\documentclass[a4paper,11pt]{article}
\usepackage{amsmath,amssymb,amsthm,MnSymbol}
\usepackage{paralist,enumitem,hyperref,url,bm,placeins}
\usepackage{color,graphicx}
\usepackage[margin=30mm]{geometry}

\newcommand{\N}{\mathbb{N}}
\newcommand{\R}{\mathbb{R}}

\newcommand{\bv}{\bm{v}}
\newcommand{\pd}{\partial}
\renewcommand{\div}{\, \mathrm{div} \,}
\newcommand{\I}{\mathbb{I}}
\newcommand{\bo}{\bm{\omega}}
\newcommand{\bnu}{\bm{n}}

\newcommand{\eps}{\varepsilon}
\newcommand{\D}{\mathrm{D}}
\newcommand{\W}{\mathrm{W}}
\newcommand{\curl}{\, \mathrm{curl}\, }

\newcommand{\pdnu}{\pd_{\bnu}}
\newcommand{\bu}{\bm{u}}
\newcommand{\nd}[1]{\pd_t^{\bullet}{#1}}
\newcommand{\inn}[2]{\langle #1, #2 \rangle}

\newcommand{\C}{\mathbb{C}}
\newcommand{\HH}{\mathbb{H}}
\newcommand{\PP}{\mathbb{P}}
\newcommand{\VV}{\mathbb{V}}

\newcommand{\bw}{\bm{w}}

\newcommand{\GL}{\mathrm{GL}}
\newcommand{\TT}{\mathbb{T}}

\theoremstyle{plain}
\newtheorem{thm}{Theorem}[section]
\newtheorem{prop}[thm]{Proposition}

\newtheorem{cor}[thm]{Corollary}
\newtheorem{remark}{Remark}[section]

\numberwithin{equation}{section}

\title{Two phase micropolar fluid flow with unmatched densities modeled by Navier--Stokes--Cahn--Hilliard systems: Local strong well-posedness and consistency estimates}
\author{Kin Shing Chan \footnotemark[1] \and Kei Fong Lam \footnotemark[1]} 
\date{ }

\begin{document}

\maketitle

\renewcommand{\thefootnote}{\fnsymbol{footnote}}
\footnotetext[1]{Department of Mathematics, Hong Kong Baptist University, Kowloon Tong, Hong Kong \tt(22482951@life.hkbu.edu.hk,akflam@hkbu.edu.hk)}

\begin{abstract}
We study a thermodynamically consistent phase field model for binary mixtures of micropolar fluids, i.e., fluids exhibiting internal rotations.  Furnishing with classical no-slip, no-spin and no-flux boundary conditions, in a smooth and bounded three-dimensional domain, we establish the well-posedness of local-in-time strong solutions. Since the model studied is a generalization of the earlier model introduced by Abels, Garcke and Gr\"un for binary Newtonian fluids with unmatched densities, we provide a consistency result between the corresponding strong solutions to both models in terms of a parameter associated to the micro-rotation viscosity. 
\end{abstract}



\section{Introduction}\label{sec:intro}
Micropolar fluids refer to a class of fluidic bodies where each material point is endowed with rigid directors \cite{Eringen1,Eringen2,Eringen3,Luka}. Their mathematical descriptions differ from conventional non-polar fluids modeled with the classical Navier--Stokes equations by possessing a non-symmetric stress tensor, as well as requiring additional balance laws for the conservation of micro-inertia. The resulting system of equations as a generalization of the classical Navier--Stokes equations read as follows for the fluid velocity $\bu$, pressure $p$, viscosity $\eta$ and micro-rotation $\bo$:
\begin{subequations}\label{microNS}
\begin{alignat}{2}
\div \bu & = 0, \label{MNS:div} \\
\pd_t \bu + (\bu \cdot \nabla) \bu & =  (\eta + \eta_r) \Delta \bu - \nabla p + 2 \eta_r \curl \bo, \label{MNS:mom} \\
\pd_t \bo + (\bu \cdot \nabla) \bo & = (c_0 + c_d - c_a) \nabla \div \bo + (c_a + c_d) \Delta \bo + 2 \eta_r(\curl \bu - 2 \bo), \label{MNS:w}
\end{alignat}
\end{subequations}
where $c_0$, $c_a$, $c_d$ are the constant coefficients of angular viscosities, and $\eta_r$ is the constant dynamic micro-rotation viscosity. We mention two interesting asymptotic limits. The first is the case $\eta_r \to 0$ where the micro-rotation equation \eqref{MNS:w} then decouples with the Navier--Stokes equations \eqref{MNS:div}-\eqref{MNS:mom}. The second is the case $\eta_r \to \infty$ where in this formal limit we obtain from \eqref{MNS:w} the relation $\bo = \frac{1}{2} \curl \bu$, which connects the micro-rotation with the well-known fluid vorticity, see e.g.~\cite{Guterres2} for some recent results in this direction.

Due to the similarity between \eqref{microNS} and the Navier--Stokes equations, various results concerning solvability, long-time behavior, global attractors, optimal time decay rates and blow-up criteria have been established, see e.g.~\cite{JLBoldMDu,Chen,FWCruz,Dong1,Dong2,GaldiRionero,Guterres,LiuZhang,Luka,Luka2D,QianChenZ} and the references cited therein. However these results concern the \emph{single-component} case where there is only one fluid occupying the domain of interest. Motivated by previous studies \cite{MDevAR,HRamkSM,PKYadSJBDS} for the two-component settings often encountered in industrial and medical applications of micro-polar fluid mixtures \cite{TAriMTNS}, a new phase field model for a binary mixture of immiscible and incompressible micro-polar fluids has been proposed in \cite{CHLS}, which reads as follows: In a bounded domain $\Omega \subset \R^3$, and for $T > 0$ an arbitrary but fixed terminal time, we consider in $\Omega \times (0,T)$
\begin{subequations}\label{bu:PF:model}
\begin{alignat}{2}
\label{model:phi} \nd{\phi} & = \Delta \mu, \\[1ex]
\label{model:mu} \mu & = \frac{\sigma}{\eps} F'(\phi) - \sigma \eps \Delta \phi,\\[1ex]
\label{model:div} \div \bu & = 0, \\[1ex]
\rho \nd{\bu} & = - \nabla p - \div( \sigma \eps \nabla \phi \otimes \nabla \phi) + \div (2 \eta(\phi) \D \bu + 2 \eta_r(\phi) \W \bu) \label{PF:lin:mom} \\[1ex] 
\notag & \quad + 2 \curl(\eta_r(\phi) \bo) + \rho'(\phi) (\nabla \mu \cdot \nabla  ) \bu, \\[1ex]
\rho \nd{\bo} & = \div (c_0(\phi) (\div \bo) \I + 2c_d(\phi) \D \bo +2 c_a(\phi) \W \bo) + 2 \eta_r(\phi)(\curl \bu - 2 \bo)  \label{PF:ang:mom} \\[1ex]
\notag  & \quad + \rho'(\phi) ( \nabla \mu \cdot \nabla  ) \bo.
\end{alignat} 
\end{subequations}
The above model comprises of a Navier--Stokes--Cahn--Hilliard system \eqref{model:phi}-\eqref{PF:lin:mom} for the mixture velocity $\bu$, phase field variable $\phi$, chemical potential $\mu$, and pressure $p$, coupled to an equation \eqref{PF:ang:mom} for the micro-rotation $\bo$.  In \eqref{model:phi} $\nd{\phi} := \pd_t \phi + \bu \cdot \nabla \phi$ is the material derivative, while in \eqref{model:mu} the positive constant $\sigma$ is associated to the surface tension. Meanwhile, $F'$ is the derivative of a double well potential $F$, which is typically chosen as a smooth quartic polynomial $F(s) = \frac{1}{4}(s^2-1)^2$ or as the Flory--Huggins logarithmic potential
\[
F(s) = \frac{\theta}{2} \lbrack (1+s)\log(1+s)+(1-s)\log(1-s)\rbrack -\frac{\theta_{0}}{2}s^2,
\]
with constants $0 < \theta < \theta_0$. By interpreting $\phi$ as the difference in volume fractions of the individual fluid constituents, we can assign the region occupied by fluid 1 as the set $\{\phi = 1\}$, the region occupied by fluid 2 as the set $\{\phi = -1\}$ and the interfacial layer separating these regions as $\{|\phi| < 1\}$ whose width scales with the parameter $\eps > 0$.

For the mixture velocity $\bu$ we choose the volume averaged velocity  \cite{AGG,FBoyer1,HDingPSCS}, which has the advantage of being solenoidal and thus renders the resulting models to be better amenable to further mathematical analysis. An alternate choice is the mass averaged velocity \cite{Aki,JLowLT,MSRGSEH} which results in models that are mathematically more difficult to treat due to strong coupling between the pressure and the chemical potential, see e.g.~\cite{2023CHNS} for further details. In \eqref{PF:lin:mom}, $\nd \bu = \pd_t \bu + (\bu \cdot \nabla) \bu$ denotes the material derivative, $\D \bu = \frac{1}{2}(\nabla \bu + (\nabla \bu)^{\top})$ is the symmetric strain tensor, while $\W \bu = \frac{1}{2}(\nabla \bu - (\nabla \bu)^{\top})$ is the anti-symmetric strain tensor.  The mass density function $\rho$ as a function of $\phi$ is defined by the algebraic relation
\begin{align}\label{rho:defn}
\rho(\phi)= \frac{\overline{\rho}_1 - \overline{\rho}_2}{2} \phi + \frac{\overline{\rho}_1 + \overline{\rho}_2}{2} \quad \text{ for } \phi \in [-1,1],
\end{align}
where $\overline{\rho}_1 > 0$, $\overline{\rho}_2> 0$ denote the constant mass densities of fluids 1 and 2, respectively.  In particular, $\rho(\phi)$ serves to interpolate between the two constant mass densities, and likewise the functions $\eta(\phi)$, $\eta_r(\phi)$, $c_0(\phi)$, $c_d(\phi)$ and $c_a(\phi)$ appearing in \eqref{PF:lin:mom} and \eqref{PF:ang:mom} interpolate between the two sets of positive viscosity values $\{\eta_i, \eta_{r,i}, c_{0,i}, c_{d,i}, c_{a,i}\}$, i.e., for $f \in \{ \eta, \eta_r, c_0, c_d, c_a\}$ we consider 
\begin{align}\label{coeff:defn}
f(\phi) := \frac{f_1 - f_2}{2} \phi + \frac{f_1 + f_2}{2} \quad \text{ for } \phi \in [-1,1].
\end{align}
Furthermore, the term $\div(\sigma \eps \nabla \phi \otimes \nabla \phi)$ accounts for capillary forces due to surface tension and the term $\rho'(\phi)(\nabla \mu \cdot \nabla ) \bu$ describes the transport mechanism by a relative flux related to diffusion of the components. In the case $\eta_r = 0$, \eqref{model:phi}-\eqref{PF:lin:mom} reduces to the model of Abels, Garcke and Gr\"un (AGG) \cite{AGG} for binary mixtures with unmatched mass densities $\overline{\rho}_1 \neq \overline{\rho}_2$.  If $\overline{\rho}_1 = \overline{\rho}_2$ also holds, then the mass density function $\rho(\phi)$ reduces to a positive constant, and we recover Model H of Hohenberg and Halperin \cite{PCHoBIH}, see also \cite{ModelH.der.}. The extension of the Navier--Stokes--Cahn--Hilliard description to micropolar fluids lies in the new equation \eqref{PF:ang:mom} as well as the coupling to \eqref{model:phi}-\eqref{PF:lin:mom} via the terms involving the dynamic micro-rotation viscosity function $\eta_r(\phi)$.

In this paper we furnish \eqref{model:phi}-\eqref{PF:lin:mom}, which we will denote as the MAGG model (for Micropolar Abels--Garcke--Gr\"un), with the standard no-slip condition for $\bu$, no-spin condition for $\bo$ and no-flux condition for $\phi$ and $\mu$:
\begin{subequations}\label{bu:PF:BC&IC}
\begin{alignat}{2}
\label{PF:BC} \bu = \bm{0}, \quad \bo = \bm{0}, \quad \pd_{\bnu} \phi = 0, \quad \pd_{\bnu} \mu & = 0 && \quad \text{ on } \partial \Omega \times (0,T), \\[1ex]
\bu(0) = \bu_0, \quad \bo(0) = \bo_0, \quad \phi(0) & = \phi_0 &&  \quad\text{ in } \Omega,
\end{alignat}
\end{subequations}
where $\bnu$ and $\pd_{\bnu}$ denotes the outward normal vector and the outer normal derivative on $\partial \Omega$, respectively.  Then, the associated energy identity for the MAGG model is given by
\begin{equation} \label{intro:Ene.ineq}
\begin{aligned}
\frac{d}{dt}E(\bu,\bo,\phi) & +\int_{\Omega} |\nabla \mu|^2 + 2 \eta(\phi) |\D \bu|^2 + 4\eta_{r}(\phi) | \tfrac{1}{2} \curl \bu - \bo|^2 \, dx \\& +\int_\Omega c_0(\phi) |\div \bo|^2 + 2 c_d(\phi) |\D \bo|^2 + 2 c_a(\phi) |\W \bo|^2 \, dx = 0,
\end{aligned}
\end{equation}
where the total energy is defined as
\begin{align}\label{defn:E}
E(\bu, \bo, \phi) := \int_\Omega \frac{\rho(\phi)}{2} |\bu|^2 + \frac{\rho(\phi)}{2} |\bo|^2 + \frac{\sigma \eps}{2} |\nabla \phi|^2 + \frac{\sigma}{\eps} F(\phi) \, dx.
\end{align}
Concerning mathematical analysis of such Navier--Stokes--Cahn--Hilliard models, we refer to \cite{HAbel.EandU.ModelH,AGAMRT} for well-posedness and regularity of solutions to Model H.  For the AGG model we mention \cite{HAbelDDHG,HAbelDDHG2} for global weak solutions, \cite{CHNS2Dstrsol,CHNS3Dstrsol} for local strong solutions, \cite{2023CHNS} for asymptotic stabilization and regularity, as well as \cite{CCGMG1,CCGMG2,CCGMG3} for long time behavior in the context of attractors, and to \cite{AbelsTera,SFrigeri1,SFrigeri2,CCGAGMG,Knopf} for similar results concerning a non-local variant of the Navier--Stokes--Cahn--Hilliard system.

For the MAGG model \eqref{model:phi}-\eqref{PF:ang:mom} the global existence of weak solutions (see Proposition \ref{prop:weaksoln} below) can be deduced from the results established in \cite{CHLS} where $F$ is taken as the Flory--Huggins logarithmic potential in order to provide a natural boundedness property for the phase field variable $\phi \in (-1,1)$ for a.e.~$(x,t) \in \Omega \times (0,T)$, and a more complicated set of boundary conditions than \eqref{bu:PF:BC&IC} accounting for moving contact lines and dynamic contact angles, see e.g.~\cite{GGM,CGGALMGHW,GLW} and the references cited therein.  

The main purpose of this paper is to establish local strong solutions for the MAGG model in three spatial dimensions, which serves to extend the results of \cite{CHNS3Dstrsol} to the micropolar setting.  Although the equations for $\bu$ and $\bo$ are similar, we encounter new difficulties arising from the fact that $\bo$ need not be solenoidal like $\bu$, as well as the coupling terms involving $\eta_r$ between \eqref{PF:lin:mom} and \eqref{PF:ang:mom}. Our strategy for establishing local strong solution takes inspiration from \cite{CHNS3Dstrsol} by employing a blend of Galerkin approxiations, Schauder fixed point theorem, suitable regularization and approximation schemes.  The main result of this paper is formulated as follows (see Section \ref{sec:pre} for the precise definitions of the function spaces):
\begin{thm}\label{thm:e&u}
Let $\Omega \subset \R^3$ be a bounded domain with $C^3$ boundary. Assuming the mass density and the viscosity functions are of the form \eqref{coeff:defn}, and the following technical assumptions on the viscosity coefficients:
\[
c_{d,i} \geq c_{a,i}, \quad 2c_{0,i} + c_{a,i} > c_{d,i} \quad \text{ for } \; i = 1, 2,
\]
along with initial data 
\begin{align*}
\bu_{0} & \in \HH^1_{\div}, \quad \bo_{0} \in \HH^1, \\
\phi_{0} & \in H^2_{n}(\Omega) \text{ s.t. } |\phi_0(x)| \leq 1 \, \forall x \in \Omega, \quad \frac{1}{|\Omega|} \int_\Omega \phi_0(x) \, dx \in (-1,1), \\
\mu_{0} & := -\sigma \eps \Delta \phi_0 + \frac{\sigma}{\eps} F'(\phi_0) \in H^{1}(\Omega),
\end{align*}
then there exists a time $T_{0} > 0$, depending on the norms of the initial data, and a quintuple of functions $(\bu,\bo, p, \phi, \mu)$ that is a strong solution to the MAGG model \eqref{bu:PF:model} on $(0,T_0)$ with the following regularities
\begin{subequations}\label{app:model1:bc}
\begin{alignat*}{2}
\bu & \in C^{0}( [0,T_{0});\HH_{\div}^{1}) \cap L^{2}( 0,T_{0};\HH_{\div}^{2}) \cap H^{1}( 0,T_{0};\HH_{\div}), \\[1ex]
 \bo & \in C^{0}( [ 0,T_{0});\HH^{1}) \cap L^{2}( 0,T_{0};\HH^{2} )  \cap H^{1}(0,T_{0};\HH), \\[1ex]
 p & \in L^{2}( 0,T_{0};H^{1}(\Omega)), \\[1ex]
 \phi & \in L^{\infty}( 0,T_{0};W^{2,6}(\Omega)) \text{ and } \left| \phi(x,t) \right| < 1 \text{ a.e. in } \Omega \times ( 0,T_{0}), \\[1ex]
 \pd_t \phi & \in L^{\infty}( 0,T_{0}; H^{1}(\Omega)^* ) \cap L^{2}( 0,T_{0};H^{1}(\Omega)), \\[1ex]
 \mu & \in L^{\infty}( 0,T_{0};H^{1}(\Omega)) \cap L^{2}( 0,T_{0};H^{3}(\Omega)), \\[1ex]
F'(\phi) & \in L^{\infty}( 0,T_{0};L^{6}(\Omega)).
\end{alignat*}
\end{subequations}
If, in addition, for some $\delta_0 > 0$, it holds that $|\phi_0(x)| \leq  1-\delta_0$ for all $x \in \Omega$, then there exists a time $T_1 \in (0,T_0)$ depending only on the norms of the initial data, such that $\Delta \phi \in L^2(0,T_1;H^2_n(\Omega))$ and the strong solution is unique in the interval $(0,T_1)$.
\end{thm}

\begin{remark}
The assertion of Theorem \ref{thm:e&u} remains valid when we consider periodic boundary conditions.
\end{remark}

\begin{remark}
The technical assumptions on $c_{0,i}$, $c_{d,i}$ and $c_{a,i}$ are not needed to establish the global existence of weak solutions to the MAGG model \eqref{bu:PF:model}, see \cite{CHLS} for more details.
\end{remark}

With the local well-posedness of strong solutions, we can  investigate the asymptotic limit as the dynamic micro-rotation viscosity $\eta_r$ tends to zero, which formally yields the AGG model from \eqref{model:phi}-\eqref{PF:lin:mom} coupled with 
\begin{align}\label{Decouple:bo}
\rho \nd{\bo} = \div (c_0(\phi) (\div \bo) \mathbb{I} + 2 c_d(\phi) \D \bo + 2 c_a(\phi) \W \bo) + \rho'(\phi) (\nabla \mu \cdot \nabla) \bo.
\end{align}
In particular, $\bo$ does not influence the evolution of $(\bu, p, \phi, \mu)$.  We mention that $\bo = \bm{0}$ is a solution to \eqref{Decouple:bo} if furnished with initial data $\bo_0 = \bm{0}$, although it is possible for non-trivial solutions to arise from non-zero initial data. Our next result details a consistency estimate between strong solutions of the MAGG model and strong solutions to the AGG model established in \cite{CHNS3Dstrsol}, which admit the same regularities as listed in Theorem~\ref{thm:e&u} for $(\bu, p, \phi, \mu)$, see Proposition \ref{thm:AGG} for more details. The strong solutions to the AGG model derived in \cite{CHNS3Dstrsol} exist on a time interval $[0,T_a]$, and without loss of generality we assume $T_1 \leq T_a$.

\begin{thm}\label{micro-AGG.consist.} Let $T_a \in [T_1,\infty)$ be a time such that on the interval $[0,T_a]$ there exists a unique strong solution $(\bu_a, p_a, \phi_a, \mu_a)$ to the AGG model subjected to the initial conditions $(\bu_a(0), \phi_a(0)) = (\bu_0, \phi_0)$. Let $R > 0$ be a fixed constant where if the standard boundary conditions \eqref{PF:BC} are considered we set $R = 1$, and if the periodic boundary conditions are considered we allow $R<\infty$ to be arbitrary. Let $\eta_r(\phi) = \eta_r \in (0,R]$ be a constant and denote by $(\bu_w, p_w, \bo_w, \phi_w, \mu_w)$ the unique strong solution to the MAGG model subjected to initial conditions $(\bu_w(0), \bo_w(0), \phi_w(0)) = (\bu_0,\bm{0}, \phi_0)$.  

Then, there exists a positive constant $C$ depending on the norms of the initial data, $T_1$, $R$ and parameters of the systems, but independent of $\eta_r \in (0,R]$ such that 
\begin{align*}
\sup_{t \in (0,T_1)} \Big ( \| \bu_w(t) - \bu_a(t) \|_{L^2}^2 + \| \bo_w(t) \|_{L^2}^2 + \| \phi_w(t) - \phi_a(t) \|_{H^2}^2 \Big ) \leq C \eta_r.
\end{align*}
\end{thm}

Similarly, setting $\overline{\rho}_1 = \overline{\rho}_2 =: \overline{\rho}$ (hence $\rho' = 0$) and $\eta_r = 0$ in \eqref{model:phi}-\eqref{PF:lin:mom} results in Model H of Hohenberg and Halperin \cite{PCHoBIH}, which admit unique strong solutions on a time interval $[0,T_h]$ with the same regularities as listed in Theorem~\ref{thm:e&u} for $(\bu, p, \phi, \mu)$, see Proposition \ref{thm:modelH} for more details. Without loss of generality we assume $T_1 \leq T_a \leq T_h$ and by leveraging the consistency estimate between the strong solutions to AGG model and the strong solutions to Model H established in \cite{CHNS3Dstrsol}, we can combine with the above result to deduce a consistency estimate between the strong solutions of MAGG model and the strong solutions of Model H in terms of $\eta_r$, $\overline{\rho}_1 - \overline{\rho}_2$ and $\frac{\overline{\rho}_1 + \overline{\rho}_2}{2} - \overline{\rho}$. 

\begin{cor}\label{micro-H.consist.} Let $T_h \in [T_1, \infty)$ be a time such that on the interval $[0,T_h]$ there exists a unique strong solution $(\bu_h, p_h, \phi_h, \mu_h)$ to Model H subjected to initial conditions $(\bu_h(0), \phi_h(0)) = (\bu_0, \phi_0)$. Let $R$ be as in Theorem \ref{micro-AGG.consist.}, $\eta_r(\phi) = \eta_r \in (0,R]$ be a constant and denote by $(\bu_w, p_w, \bo_w, \phi_w, \mu_w)$ the unique strong solution to the MAGG model subjected to initial conditions $(\bu_w(0), \bo_w(0), \phi_w(0)) = (\bu_0,\bm{0}, \phi_0)$.  

Then, there exists a positive constant $C$ depending on the norms of the initial data, $T_1$, $R$ and parameters of the systems, but independent of $\eta_r \in (0,R]$, $\overline{\rho}_1 - \overline{\rho}_2$ and $\frac{\overline{\rho}_1 + \overline{\rho}_2}{2} - \overline{\rho}$ such that 
\begin{align*}
& \sup_{t \in (0,T_1)} \Big ( \| \bu_w(t) - \bu_h(t) \|_{L^2}^2 + \| \bo_w(t) \|_{L^2}^2 + \| \phi_w(t) - \phi_h(t) \|_{H^2}^2 \Big ) \\
& \quad \leq C \Big ( \eta_r + |\overline{\rho}_1 - \overline{\rho}_2| +  \Big | \frac{\overline{\rho}_1 + \overline{\rho}_2}{2} - \overline{\rho} \Big | \Big ).
\end{align*}
\end{cor}

\begin{remark}[Two-dimensional setting]
All of the above main results are valid in the two-dimensional setting, however the Navier--Stokes and micro-rotation components of the MAGG model will need to be modified. We restrict the three-dimensional dynamics to the $(x_1,x_2)$-plane, with the flow independent of the $x_3$-coordinate, while having the axes of micro-rotation parallel to the $x_3$-axis. Namely, $\bu$ and $\bo$ reduces to $\bu = (u_1(t,x_1,x_2), u_2(t,x_1,x_2),0)^{\top}$ and $\bo = (0,0,\omega(t,x_1,x_2))^{\top}$, while $\phi$ and $\mu$ are scalar functions of $t$, $x_1$ and $x_2$. The curl of a two dimensional vector $\bm{v} = (v_1, v_2)^{\top}$ is the scalar quantity
\[
\curl_2 \bm{v} = \pd_1 v_2 - \pd_2 v_1,
\]
and we define the rotation vector associated to a scalar $w$ as
\[
\curl_1 w = \begin{pmatrix}
    \pd_2 w \\ - \pd_1 w
\end{pmatrix}.
\]
Then, in a bounded domain $\Omega \subset \R^2$, with an abuse of notation we write $\bu = (u_1, u_2)^{\top}$, and take the first and second components of \eqref{PF:lin:mom} for $\bu$, as well as the third component of \eqref{PF:ang:mom} for the scalar quantity $\omega$, so that the two-dimensional MAGG model reads as
\begin{subequations}\label{local:model:equ:2D}
\begin{alignat}{2}
\label{2D:equ:CH} \nd \phi & =\Delta \mu, \\[1ex] 
\label{2D:equ:mu} \mu & = -\sigma \eps \Delta \phi + \frac{\sigma}{\eps} F'(\phi), \\[1ex]
\label{2D:div:0} \div \bu & = 0, \\[1ex]
\label{2D:equ:bu} \rho \nd \bu & =  - \nabla p - \div (\sigma \eps \nabla \phi \otimes \nabla \phi) + \div (2 \eta(\phi) \D \bu + 2 \eta_r(\phi) \W \bu) \\[1ex]
\notag & \quad = 2\curl_1\left( \eta_{r}(\phi) \omega \right) + \rho'(\phi) (\nabla \mu \cdot  \bu), \\[1ex]
 \label{2D:equ:bo} \rho \nd \omega & = \div ((c_d(\phi) + c_a(\phi)) \nabla \omega) + 2\eta_{r}(\phi)(\curl_2 \bu - 2\omega) + \rho'(\phi) \nabla \mu \cdot \nabla \omega,
\end{alignat}
\end{subequations}
where $\nabla$ and $\div$ denote the two-dimensional spatial gradient and divergence operator.
\end{remark}

The rest of the paper is organized as follows: In section \ref{sec:pre} we provide notation, definition of relevant fuction spaces and useful preliminary results concerning the convective Cahn--Hilliard equation with logarithmic potentials and the global weak existence to the MAGG model. In section \ref{sec:exist} we prove the existence of local-in-time strong solutions employing a mixture of Galerkin approximations and fixed point arguments. Uniqueness and continuous dependence on initial data are established in section \ref{sec:unique}, while section \ref{sec:etar0} is devoted to studying the nonpolar limit $\eta_r \to 0$ and deriving consistency estimates between the MAGG model with the AGG model and Model H.

\section{Preliminaries}\label{sec:pre}
\subsection{Notation}
For the subsequent sections we often make use of the Einstein summation convention and neglect the basis vector elements. For a vector $\bv = (v_i)$ and a second order tensor $\bm{A} = (A_{ij})$, the gradient $\nabla \bv$ and divergence $\div \bm{A}$ are defined as
\[
(\nabla \bv)_{ij} = \pd_i v_j = \frac{\pd v_j}{\pd x_i}, \quad (\div \bm{A})_j = \pd_i A_{ij}.
\]
Meanwhile the vector $\bm{v} \cdot \bm{A}$ is defined as
\[
(\bm{v} \cdot \bm{A})_j = v_i A_{ij}.
\]
The Frobenius product $\bm{A}: \bm{B}$ of two second order tensors $\bm{A}$ and $\bm{B}$ is defined as $\bm{A} : \bm{B} = A_{ij} B_{ij}$. In three spatial dimensions, the entries of the third order Levi-Civita tensor $\bm{\eps} = (\eps_{ijk})$ are defined as
\[
\eps_{ijk} = \begin{cases}
1 & \text{ if } (i,j,k) \text{ is } (1,2,3), (2,3,1), \text{ or } (3,1,2), \\
-1 & \text{ if } (i,j,k) \text{ is } (3,2,1), (1,3,2), \text{ or } (2,1,3), \\
0 & \text{ if } i = j, \text{ or } j = k, \text{ or } k = i.
\end{cases}
\]
Then, in three spatial dimensions, the following properties are valid:
\begin{align}\label{eps:prop}
\eps_{ljk} \eps_{mjk} = 2 \delta_{lm}, \quad \eps_{jik} \eps_{jlm} = \delta_{il} \delta_{km} - \delta_{im} \delta_{kl}.
\end{align}
Furthermore, the cross product $\bm{a} \times \bm{b}$ between two vectors $\bm{a}$ and $\bm{b}$, as well as the curl of a vector $\bv$ are defined as
\[
(\bm{a} \times \bm{b})_j = \eps_{jkl} a_k b_l, \quad (\curl \bv)_j = \eps_{jkl} \pd_k v_l = (\nabla \times \bv)_j.
\]
Then, we have the integration by parts formula involving the curl operator:
\begin{equation}\label{IBP:curl}
\begin{aligned}
\int_\Omega \curl \bm{a} \cdot \bm{b} \, dx & = \int_{\Omega} \bm{a} \cdot \curl \bm{b} \, dx + \int_{\pd \Omega} (\bm{a} \times \bm{b}) \cdot \bnu \, dS \\
& = \int_{\Omega} \bm{a} \cdot \curl \bm{b} \, dx - \int_{\pd \Omega}  (\bm{a} \times \bnu) \cdot \bm{b} \, dS.
\end{aligned}
\end{equation}
For a real Banach space $X$ and its topological dual $X^*$, we denote by $\inn{f}{g}_X$ for $f \in X^*$ and $g \in X$ as the corresponding duality pairing.  The continuous embedding of $X$ into $Y$ is denoted by $X \subset Y$, while the compact embedding of $X$ into $Y$ is denoted by $X \Subset Y$.

For $1 \leq p \leq \infty$, the Bochner space $L^p(0,T;X)$ denotes the set of all strongly measurable $p$-integrable functions (if $p < \infty$) or essentially bounded functions (if $p = \infty$) on the time interval $[0,T]$ with values in the Banach space $X$. The space $W^{1,p}(0,T;X)$ denotes all $u \in L^p(0,T;X)$ such that its vector-valued distributional derivative $\frac{du}{dt} \in L^p(0,T;X)$. Furthermore, $C^0([0,T];X)$ denotes the Banach space of all bounded and continuous functions $u:[0,T] \to X$ equipped with the supremum norm, while $C_w^0([0,T];X)$ is the space of bounded and weakly continuous functions $u:[0,T] \to X$, i.e., $\inn{f}{u} : [0,T] \to \R$ is continuous for all $f \in X^*$. We employ the notation $C^\infty_0(0,T;X)$ to denote the vector space of all smooth and compactly support (in time) functions $u:(0,T) \to X$.

For a bounded domain $\Omega \subset \R^3$, we denote by $L^p(\Omega)$ and $W^{k,p}(\Omega)$, for $1 \leq p \leq \infty$ and $k \geq 0$ to be the Lebesgue and Sobolev spaces over $\Omega$. When $k \in \mathbb{N}$ and $p = 2$, we use the notation $H^k(\Omega) = W^{k,2}(\Omega)$. The $L^2(\Omega)$-inner product is denoted by $(\cdot,\cdot)$ with associated norm $\| \cdot \|$, while the norm of $W^{k,p}(\Omega)$ is denoted by $\| \cdot \|_{W^{k,p}}$. For vectorial functions, we use $L^p(\Omega;\R^d)$ and $W^{k,p}(\Omega;\R^d)$ to denote the corresponding Lebesgue and Sobolev spaces. For convenience we use the notation $Q := \Omega \times [0,T]$, $\Sigma := \pd \Omega \times [0,T]$, as well as $(\cdot,\cdot)_Q$ to denote the $L^2(Q)$-inner product. Due to the homogeneous Neumann boundary conditions we introduce the notation
\[
H^2_n(\Omega) := \{ f \in H^2(\Omega) \, | \, \pdnu f = 0 \text{ on } \pd \Omega \}.
\]
The generalized mean value is defined as
\[
\overline{f} = \begin{cases}
\frac{1}{|\Omega|}\int_{\Omega} f \, dx & \text{ if } f \in L^1(\Omega), \\
\frac{1}{|\Omega|}\inn{f}{1}_{H^{1}(\Omega)} & \text{ if } f \in H^1(\Omega)^*.
\end{cases}
\]
We recall some useful Sobolev inequalities valid in three spatial dimensions:
\begin{subequations}
\begin{alignat}{3}
\label{Gen.Poin.} \left\| f \right\|_{H^{1}}^{2} & \leq C\left( \left| \overline{f} \right|^{2} + \left\| \nabla f \right\|_{L^{2}}^{2} \right)  && \quad \forall f \in H^{1}(\Omega) \\ 
\label{L3estimate} \left\| f \right\|_{L^{3}} & \leq C\left\| f \right\|_{L^{2}}^{\frac{1}{2}}\left\| f \right\|_{H^{1}}^{\frac{1}{2}} && \quad \forall f \in H^{1}(\Omega) \\ 
\label{Linftyestimate} \left\| f \right\|_{L^{\infty}} & \leq C\left\| f \right\|_{H^{1}}^{\frac{1}{2}}\left\| f \right\|_{H^{2}}^{\frac{1}{2}} && \quad \forall f \in H^{2}(\Omega)\\
\label{grad:L4estimate} \left\| \nabla f \right\|_{L^{4}(\Omega)} & \leq C\left\| f \right\|_{L^{\infty}}^{\frac{1}{2}} \left\| f \right\|_{H^{2}}^{\frac{1}{2}} && \quad \forall f \in H^{2}(\Omega) \\
\label{W1,4estimate} \left\| f \right\|_{W^{1,4}} & \leq C\left\| f \right\|_{H^{1}}^{\frac{5}{8}}\left\| f \right\|_{W^{2,6}}^{\frac{3}{8}} && \quad \forall f \in W^{2,6}(\Omega),
\end{alignat}
\end{subequations}
where in the above the positive constants denoted by $C$ depend only on $\Omega$.

For the velocity field and the micro-rotation field, we consider the following function spaces:
\begin{equation*}
\begin{aligned}
& \HH_{\div} = \HH^0_{\div}:= \left\{ \bu \in L^{2}(\Omega;\R^3) \, | \, \div \bu = 0\text{ in } \Omega, \, \bu \cdot \bnu = 0 \text{ on } \partial\Omega \right \}, \\
& \HH_{\div}^{1} = \left\{ \bu \in H^{1}(\Omega;\R^3) \, | \,\div \bu = 0\text{ in } \Omega, \, \bu = \bm{0}\text{ on } \partial\Omega \right\}, \\
& \HH = \HH^0 := \left\{ \bo \in L^{2}(\Omega;\R^3) \, | \, \bo = \bm{0}\text{ on } \partial\Omega \right\}, \\
& \HH^{1} = \left\{ \bo \in H^{1}(\Omega;\R^3) \, | \, \bo = \bm{0} \text{ on } \partial\Omega \right\}, \\
& \HH_{\div}^{k} = H^{k}(\Omega;\R^3) \cap \HH_{\div}^{1}, \quad  \HH^{k} = H^{k}(\Omega;\R^3) \cap \HH^{1} \quad \text{ for } k \in \mathbb{N} \text{ and } k \geq 2.
\end{aligned}
\end{equation*}
We equip $\HH^1$ with the inner product $(\bu,\bm{v})_{\HH^{1}} = (\nabla \bu,\nabla \bm{v})$ and the norm $\left\| \bu \right\|_{\HH^{1}} = \left\| \nabla \bu \right\|_{L^{2}}$ and likewise for $\HH_{\div}^1$.  Then, upon recalling Korn's inequality
\begin{equation} \label{Korn.ineq.}
\left\| \nabla \bu \right\|_{L^{2}} \leq C\left\| \D\bu \right\|_{L^{2}} \leq C\left\| \nabla \bu \right\|_{L^{2}} \quad \forall \bu \in \HH^{1},
\end{equation}
we see that $\| \D \bu \|_{L^2}$ is an equivalent norm to the $\HH^1$-norm.

We recall the Leray projection $\PP:L^{2}(\Omega;\R^3) \rightarrow \HH_{\div}$ and denote by $\bm{A} := \PP(-\Delta)$ as the Stokes operator. Then we can equip the space $\HH_{\div}^2 := H^2(\Omega;\R^3) \cap \HH_{\div}^1$ with the inner product $(\bu, \bm{v})_{\HH_{\div}^2} := (\bm{A} \bu, \bm{A} \bm{v})$ and the norm $\| \bu \|_{\HH_{\div}^2} := \| \bm{A} \bu \|_{L^2}$. Furthermore, it holds that 
\[
\| \bu \|_{\HH^2} \leq C \| \bu \|_{\HH_{\div}^2} \quad \forall \bu \in \HH_{\div}^2,
\]
with a positive constant $C$. 

\subsection{Preliminary results}
Without loss of generality we set $\sigma = \eps = 1$ in \eqref{model:mu}.
\subsubsection{On the convective Cahn--Hilliard equation with logarithmic potentials}
We report the following collection of results from \cite[Theorem A.1]{CHNS3Dstrsol} concerning the convective Cahn--Hilliard equation with logarithmic potential:
\begin{subequations} \label{Pre.result.CH.eq}
\begin{alignat}{2}
& \partial_{t}\phi + \bu \cdot \nabla\phi = \Delta\mu && \quad \text{ in } \Omega \times (0,T), \\
& \mu = \alpha\partial_{t}\phi - \Delta\phi +F'(\phi) && \quad \text{ on } \Omega \times (0,T),
\end{alignat}
\end{subequations}
furnished with the following boundary and initial conditions:
\begin{subequations}\label{Pre.result.CH.bc&ic}
\begin{alignat}{2}
& \partial_{n}\phi = \partial_{n}\mu = 0 && \quad \text{ on } \partial\Omega \times (0,T), \\
& \phi \vert_{t = 0} = \phi_{0}  && \quad \text{ in }\Omega.
\end{alignat}
\end{subequations}

\begin{prop} \label{Pre.result}
Let $\Omega \subset \R^3$ be a bounded domain with $C^3$ boundary. Given a solenoidal velocity field $\bu \in L^\infty(0,T;\HH_{\div} \cap L^3(\Omega;\R^3))$, constant $\alpha > 0$ and $\phi_0 \in H^1(\Omega) \cap L^\infty(\Omega)$ such that $\| \phi_0 \|_{L^\infty} \leq 1$ and $|\overline{\phi_0}| < 1$, there exists a unique strong solution to \eqref{Pre.result.CH.eq}-\eqref{Pre.result.CH.bc&ic} with
\begin{subequations}
\begin{align*}
\phi &  \in L^{\infty}( 0,T;H^{1}(\Omega)) \text{ s.t. } |\phi(x,t)| < 1 \text{ a.e.~in } \Omega \times (0,T), \\
\phi & \in L^{2}( 0,T;H^{2}_n(\Omega)) \cap H^{1}( 0,T;L^{2}(\Omega)),\\
\mu &  \in L^{2}( 0,T;H^{2}_n(\Omega) ), \quad  F'(\phi) \in L^{2}( 0,T;L^{2}(\Omega)), 
\end{align*}
\end{subequations}
satisfying \eqref{Pre.result.CH.eq} a.e.~in $\Omega \times (0,T)$, \eqref{Pre.result.CH.bc&ic} a.e.~on $\partial\Omega \times (0,T)$ and $\phi \vert_{t = 0} = \phi_0$ in $\Omega$. Additionally, the following regularity results hold:
\begin{enumerate}
\item[$(\mathrm{R1})$] If $- \Delta\phi_{0} + F'( \phi_{0}) \in L^{2}(\Omega)$ and $\partial_{t}\bu \in L^{\frac{4}{3}}( 0,T;L^{1}(\Omega))$, we have
\begin{subequations}
\begin{align*}
& \partial_{t}\phi \in L^{\infty} ( 0,T;L^{2}(\Omega)) \cap L^{2}( 0,T;H^{1}(\Omega) ),\\
& \phi \in L^{\infty}( 0,T;H^{2}(\Omega) ), \\
& \mu \in L^{\infty}( 0,T;H^{2}_n(\Omega) ).
\end{align*} 
\end{subequations}
\item[$(\mathrm{R2})$] In addition to $(\mathrm{R1})$, suppose for some $\delta_0 \in (0,1)$ it holds that $\| \phi_0 \|_{L^\infty} < 1- \delta_0$, then there exists $\delta > 0$ such that 
\begin{equation*}
\max_{(x,t) \in \Omega \times [ 0,T]} \left| \phi(x,t) \right| \leq 1 - \delta \quad \text{ and } \quad \phi \in L^{2}( 0,T;H^{3}(\Omega)).
\end{equation*}
\item[$(\mathrm{R3})$] In addition to $(\mathrm{R2})$, suppose $\phi_{0} \in H^{3}(\Omega) \cap H^2_n(\Omega)$ and
$\partial_{t}\bu \in L^{2}( 0,T;L^{\frac{6}{5}}(\Omega))$, then we have
\begin{subequations}
\begin{align*}
& \phi \in L^{\infty}( 0,T;H^{3}(\Omega)), \\
& \partial_{t}\phi \in L^{\infty}( 0,T;H^{1}(\Omega)) \cap L^{2}( 0,T;H^{2}(\Omega)) \\
& \partial_{tt}\phi \in L^{2}( 0,T;L^{2}(\Omega)), \\
& \partial_{t}\mu \in L^{2}( 0,T;L^{2}(\Omega)).
\end{align*}
\end{subequations}
\end{enumerate}
\end{prop}

\subsubsection{Approximation of the initial data}
In the following we will need to consider an approximation to the initial data $\phi_0$ in order to invoke the regularity results in Proposition \ref{Pre.result}.  Thus, for $k \in \N$ we consider the Lipschitz continuous truncation operator $h_k: \R \to \R$ defined as
\[
h_{k}(z) = \begin{cases}
-k & \text{ if } z \leq - k,  \\
z & \text{ if } |z| \leq k, \\
k & \text{ if } z \geq k, 
\end{cases}
\]
and for the data $\widetilde{\mu}_0 := - \Delta \phi_0 + F'(\phi_0) \in H^1(\Omega)$ we infer that $\widetilde{\mu}_{0,k} := h_k(\widetilde{\mu}_0) \in H^1(\Omega) \cap L^\infty(\Omega)$ and satisfies
\begin{align}\label{prop1.approx.IC}
\| {\widetilde{\mu}}_{0,k} \|_{H^{1}} \leq \| {\widetilde{\mu}}_{0} \|_{H^{1}}.
\end{align}
Then, let us consider the elliptic problem
\begin{equation} \label{approx.phi0.problem}
\begin{aligned}
 -\Delta\phi_{0,k} + F'( \phi_{0,k} ) & = \widetilde{\mu}_{0,k} && \text{ in } \Omega \\
\pdnu \phi_{0,k} & = 0 && \text{ on } \partial\Omega,
\end{aligned}
\end{equation}
whose unique solution $\phi_{0,k} \in H^2_n(\Omega)$ with $F'(\phi_{0,k}) \in L^2(\Omega)$ satisfies the following properties, see e.g.~\cite{CHNS3Dstrsol,AGAMRT}:  There exists $\widetilde{m} \in (0,1)$ independent of $k$, and $\overline{k}$ sufficiently large such that for all $k > \overline{k}$,
\begin{align}\label{prop.approx.IC.phi0}
\| \phi_{0,k} \|_{H^{2}} \leq C( 1 + \| {\widetilde{\mu}}_{0} \|_{L^{2}} ), \quad \| \phi_{0,k} \|_{H^{1}} \leq 1 + \| \phi_{0} \|_{H^{1}}, \quad  | \overline{\phi_{0,k}} | \leq \widetilde{m} < 1,
\end{align}
with positive constant $C$ independent of $k$. Furthermore, it holds that 
\[
\| F'(\phi_{0,k}) \|_{L^\infty} \leq \| \widetilde{\mu}_{0,k} \|_{L^\infty} \leq k
\]
and so there exists $\delta = \delta(k) > 0$ such that 
\begin{align}\label{prop.approx.IC.delta}
\| \phi_{0,k} \|_{L^{\infty}} \leq 1 - \delta.
\end{align}
This yields further regularities $F'(\phi_{0,k}) \in H^1(\Omega)$ and $\phi_{0,k} \in H^3(\Omega)$.  Lastly, as $\widetilde{\mu}_{0,k}\rightarrow \widetilde{\mu}_{0}$ in $L^2(\Omega)$ it holds that $\phi_{0,k}\rightarrow \phi_{0}$ in $H^1 (\Omega)$.

\subsubsection{Global weak solutions to the MAGG model}

We conclude this section by providing a global weak existence result to \eqref{bu:PF:model}-\eqref{bu:PF:BC&IC} that can be deduced from \cite{CHLS}:
\begin{prop}\label{prop:weaksoln}
Let $\bu_0 \in \HH_{\div}$, $\bo_0 \in \HH$, $\phi_0 \in H^1(\Omega)$ with $\| \phi_0 \|_{L^\infty} \leq 1$, $F(\phi_0) \in L^1(\Omega)$ and $|\overline{\phi_0}| < 1$, then there exists a global weak solution $(\bu, \bo, \phi, \mu)$ to \eqref{bu:PF:model}-\eqref{bu:PF:BC&IC} on $[0,T]$ satisfying
\begin{itemize}
\item Regularity
\begin{equation*}
\begin{alignedat}{3}
\bu & \in C_w([0,T]; \HH_{\div}) \cap L^2(0,T;\HH^1_{\div}), \\
\bo & \in C_w([0,T]; \HH) \cap L^2(0,T;\HH^1), \\
\phi & \in C_w([0,T];H^1(\Omega)) \cap L^2(0,T;H^2(\Omega)) \text{ s.t. }  |\phi| < 1 \text{ a.e.~in } \Omega \times (0,T), \\
\mu & \in L^2(0,T;H^1(\Omega)).
\end{alignedat}
\end{equation*}
\item Equations
\begin{subequations}
\begin{alignat}{2}
\label{weak:ns} 0 &=  -(\rho(\phi) \bu, \pd_t \bm{v})_Q + (\div(\rho \bu \otimes \bu), \bm{v})_Q  \\
\notag & \quad + (2 \eta(\phi) \D \bu, \D \bm{v})_Q + (2 \eta_r(\phi) \W \bu, \W \bm{v})_Q - (2\eta_r(\phi) \bo, \curl \bm{v})_Q  \\
\notag & \quad - (( \bm{J} \otimes \bu), \nabla \bm{v})_Q - (\mu \nabla \phi, \bm{v})_Q \\
\label{weak:w} 0 & = -(\rho(\phi) \bo, \pd_t \bm{z})_Q + (\div (\rho \bu \otimes \bo), \bm{z})_Q \\
\notag & \quad + (c_0(\phi) \div \bo, \div \bm{z})_Q  + (2 c_d(\phi) \D \bo, \D \bm{z})_Q + (2 c_a(\phi) \W \bo, \W \bm{z})_Q\\
\notag & \quad  - (2 \eta_r(\phi) (\curl \bu - 2 \bo), \bm{z})_Q - ((\bm{J} \otimes \bo), \nabla \bm{z})_Q, \\
\label{weak:phi} 0 & = -(\phi, \pd_t \xi)_Q + (\bu \cdot \nabla \phi, \xi)_Q + (\nabla \mu, \nabla \xi)_Q, 
\end{alignat}
\end{subequations}
holding for all $\bm{v} \in C^\infty_0(0,T;\C^\infty_{0}(\overline{\Omega};\R^3))$ such that $\div \bm{v} = 0$, $\bm{z} \in C^\infty_0(0,T;\C^\infty_0(\overline{\Omega};\R^3))$, and $\xi \in C^\infty_0(0,T;C^1(\overline{\Omega}))$, along with 
\begin{subequations}
\begin{alignat}{2}
\label{weak:mu} \mu & = - \Delta \phi + F'(\phi) && \quad \text{ a.e.~in } Q, \\
\label{weak:J} \bm{J} & = -\tfrac{\overline{\rho}_1 - \overline{\rho}_2}{2} \nabla \mu && \quad \text{ a.e.~in } Q, \\
\label{weak:rho} \rho(\phi) & = \tfrac{\overline{\rho}_1 - \overline{\rho}_2}{2} \phi + \tfrac{\overline{\rho}_1+ \overline{\rho}_2}{2} && \quad \text{ a.e.~in } Q.
\end{alignat}
\end{subequations}
\item Energy inequality
\begin{equation}\label{energy:ineq}
\begin{aligned}
& E(\bu(t), \bo(t), \phi(t)) \\
& \qquad + \int_s^t \int_\Omega |\nabla \mu|^2 + 2 \eta(\phi) |\D \bu|^2 + 4 \eta_r(\phi) | \tfrac{1}{2} \curl \bu - \bo|^2 \, dx \, d\tau \\
& \qquad + \int_s^t \int_\Omega c_0(\phi) |\div \bo|^2 + 2 c_d(\phi) |\D \bo|^2 + 2 c_a(\phi) |\W \bo|^2 \, dx \, d \tau \\
& \quad \leq E(\bu(s), \bo(s), \phi(s))
\end{aligned}
\end{equation}
holds for all $t \in [s,\infty)$ and almost all $s \in [0,\infty)$ including $s = 0$, where the total energy $E$ is given by
\begin{equation}\label{total:energy}
\begin{aligned}
E(\bu, \bo, \phi) & := \int_\Omega \frac{\rho(\phi)}{2} |\bu|^2 + \frac{\rho(\phi)}{2} |\bo|^2 + \frac{1}{2}|\nabla \phi|^2 + F(\phi) \, dx.
\end{aligned}
\end{equation}
\item Initial conditions
\[
(\bu, \bo, \phi) \vert_{t = 0} = (\bu_0, \bo_0, \phi_0).
\]
\end{itemize}
\end{prop}

\section{Existence of local-in-time strong solutions}\label{sec:exist}
Following along similar ideas as in \cite{CHNS3Dstrsol}, we first approximate the initial condition $\phi_{0}$ so that we can employ Proposition \ref{Pre.result} to gain sufficient regularity for the approximated solutions. Then, we develop a suitable approximation scheme using a mixture of fixed point theorem and Faedo--Galerkin approximation to establish the existence and uniqueness of the approximate solutions.  Lastly, we pass to the limit to establish the existence of solutions to the original problem. 

\subsection{Approximate problem}
Let $\{\bm{Y}_j \}_{j = 1}^\infty$ be the set of eigenfunctions of the Stokes operator $\bm{A}$, and for any integer $m \geq 1$ we introduce the finite dimensional subspace $\VV_{\bu,m} := \text{span}\{\bm{Y}_1, \dots, \bm{Y}_m\}$ with orthogonal projection $\PP_{\bu,m}$. Let $\{\bm{Z}_j\}_{j=1}^\infty$ be the eigenfunctions of the vectorial Dirichlet Laplacian operator, and for any integer $m \geq 1$ we introduce the finite dimensional subspace $\VV_{\bo,m} := \text{span}\{\bm{Z}_1, \dots, \bm{Z}_m\}$ with orthogonal projection $\PP_{\bo,m}$. Due to the regularity of $\pd \Omega$, it follows from regularity theory that $\bm{Y}_j \in H^3(\Omega;\R^3) \cap \HH^1_{\div}$ and $\bm{Z}_j \in H^3(\Omega;\R^3)$. Furthermore, there exist positive constants $C_m$ depending on $m$ such that 
\begin{align}\label{inverse:est}
\| \bm{f} \|_{H^1} \leq C_m \| \bm{f} \|, \quad \| \bm{f} \|_{H^2} \leq C_m \| \bm{f} \|, \quad \| \bm{f} \|_{H^3} \leq C_m \| \bm{f} \|
\end{align}
for all $\bm{f} \in \VV_{\bu,m}$ or $\VV_{\bo,m}$.

Fix $m \in \mathbb{N}$ and $\alpha \in (0,1)$, we seek an approximation solution quadruple $(\bu_m, \bo_m, \phi_m, \mu_m)$ to the following system for all $t \in [0,T]$
\begin{subequations}\label{approx:problem}
\begin{alignat}{2}
\label{app:u} & (\rho(\phi_m) \pd_t \bu_m, \bm{v}) + (\rho(\phi_m) (\bu_m \cdot \nabla) \bu_m, \bm{v}) + (2 \eta(\phi_m) \D \bu_m, \D \bm{v}) \\
\notag & \qquad + (2 \eta_r(\phi_m) \W \bu_m, \W \bm{v}) - \tfrac{\overline{\rho}_1 - \overline{\rho}_2}{2} ((\nabla \mu_m \cdot \nabla) \bu_m, \bm{v}) \\
\notag & \quad = (\mu_m \nabla \phi_m, \bm{v}) + 2(\curl(\eta_r(\phi_m) \bo_m), \bm{v}) \quad \forall \bm{v} \in \VV_{\bu,m}, \\[1ex]
\label{app:w} & (\rho(\phi_m) \pd_t \bo_m, \bm{z}) + (\rho(\phi_m) (\bu_m \cdot \nabla) \bo_m, \bm{z}) + (c_0(\phi_m) \div \bo_m, \div \bm{z}) \\
\notag & \qquad + (2c_d(\phi_m) \D \bo_m, \D \bm{z})  +  (2 c_a(\phi_m) \W \bo_m, \W \bm{z}) - \tfrac{\overline{\rho}_1 - \overline{\rho}_2}{2} ((\nabla \mu_m \cdot \nabla) \bo_m, \bm{z}) \\
\notag & \quad = (2\eta_r(\phi_m) \curl \bu_m, \bm{z}) - (4 \eta_r(\phi_m) \bo_m, \bm{z}) \quad \forall \bm{z} \in \VV_{\bo,m}, \\[1ex]
\label{app:phi} & \pd_t \phi_m + \bu_m \cdot \nabla \phi_m = \Delta \mu_m \quad \text{ a.e.~in } \Omega, \\
\label{app:mu} & \mu_m = \alpha \pd_t \phi_m - \Delta \phi_m + F'(\phi_m) \quad \text{ a.e.~in } \Omega,
\end{alignat}
\end{subequations}
together with the following initial-boundary conditions:
\begin{equation*}
\begin{alignedat}{2}
\bu_m = \bm{0}, \quad \bo_m = \bm{0}, \quad \pdnu \phi_m = 0, \quad \pdnu \mu_m & = 0 && \quad \text{ on } \Sigma, \\
\bu_m(0) = \PP_{\bu,m} \bu_0, \quad \bo_m(0) = \PP_{\bo,m} \bo_0, \quad \phi_m(0) & = \phi_{0,k} && \quad \text{ in } \Omega,
\end{alignedat}
\end{equation*}
where $\phi_{0,k}$ is defined in \eqref{approx.phi0.problem}. 

\begin{prop}\label{prop:Galerkin}
For any $T > 0$ and for all $k > 0$, $\alpha \in (0,1)$ and $m \in \mathbb{N}$, there exists at least one solution quadruple $(\bu_m, \bo_m, \phi_m, \mu_m)$ to \eqref{approx:problem} and some $\delta \in (0,1)$ satisfying the following regularities:
\begin{align*}
\bu_m & \in H^1([0,T];\VV_{\bu,m}), \\
\bo_m & \in H^1([0,T];\VV_{\bo,m}), \\
\phi_m & \in L^\infty(0,T;H^3(\Omega)\cap H^2_n(\Omega)) \text{ s.t. } |\phi_m(x,t)| \leq 1-\delta \text{ a.e.~in } Q, \\
\pd_t \phi_m & \in L^\infty(0,T;H^1(\Omega)) \cap L^2(0,T;H^2_n(\Omega)), \\
\mu_m & \in L^\infty(0,T;H^2_n(\Omega)) \cap H^1(0,T;L^2(\Omega)).
\end{align*}
\end{prop}

The proof of Proposition \eqref{Galerkin} is divided in the following subsections.
\subsubsection{Galerkin approximation}
Fix $\bm{U} \in H^1(0,T;\VV_{\bu,m})$ and consider the following viscous Cahn--Hilliard equation
\begin{subequations}\label{viscousCH}
\begin{alignat}{2}
& \pd_t \phi + \bm{U} \cdot \nabla \phi = \Delta \mu && \quad \text{ in } Q, \\
& \mu = \alpha \pd_t \phi - \Delta \phi + F'(\phi) &&  \quad \text{ in } Q, \\
& \pdnu \phi = 0, \quad \pdnu \mu = 0 && \quad \text{ on } \Sigma, \\
&\phi(0) = \phi_{0,k} && \quad \text{ in } \Omega.
\end{alignat}
\end{subequations}
Invoking Proposition \ref{Pre.result} we deduce the existence of a unique solution $(\phi_{U}, \mu_U)$ to \eqref{viscousCH} such that for some $\widetilde{\delta} > 0$ depending on $\alpha$ and $k$ it holds that
\begin{align*}
\phi_U & \in L^\infty(0,T;H^3(\Omega)\cap H^2_n(\Omega)) \text{ s.t. } |\phi_U(x,t)| \leq 1-\widetilde{\delta} \text{ a.e.~in } Q, \\
\pd_t \phi_U & \in L^\infty(0,T;H^1(\Omega)) \cap L^2(0,T;H^2_n(\Omega)), \\
\mu_U & \in L^\infty(0,T;H^2_n(\Omega)) \cap H^1(0,T;L^2(\Omega)).
\end{align*}
Furthermore, we also infer that for any $T>0$, the following $L^2$ and energy identities
\begin{align*}
& \frac{1}{2} \frac{d}{dt} \Big ( \| \phi_U \|^2 + \alpha \| \nabla \phi_U \|^2 \Big )+ \| \Delta \phi_U \|^2 + \int_\Omega F''(\phi_U) |\nabla \phi_U|^2 \, dx = 0, \\
& \frac{d}{dt} \int_\Omega \frac{1}{2} |\nabla \phi_U|^2 + F(\phi_U) \, dx + \alpha \| \pd_t \phi_U \|^2 + \| \nabla \mu_U \|^2 = \int_\Omega \phi_U \bm{U} \cdot \nabla \mu_U \, dx.
\end{align*}
From the explicit formula for $F$ we note that $F''(s) \geq - \theta_0$, while using the Neumann boundary condition we infer
\[
\| \nabla \phi_U \|^2 = (-\Delta \phi_U, \phi_U) \leq \| \Delta \phi_U  \| \| \phi_U \|.
\]
Hence, by a Gr\"onwall argument we obtain the estimates 
\begin{align}
\label{vCH:1} & \sup_{t \in (0,T]} \Big ( \| \phi_U \|^2 + \alpha \| \nabla \phi_U \|^2 \Big ) + \int_0^T \| \Delta \phi_U \|^2 \, dt \leq \| \phi_{0,k} \|^2 + \alpha \| \nabla \phi_{0,k} \|^2 + 2\theta_0^2 |\Omega| T, \\
\label{vCH:2} & \sup_{t \in (0,T]} \int_\Omega \frac{1}{2} |\nabla \phi_U(t) |^2 + F(\phi_U(t)) \, dx + \int_0^T \alpha \| \pd_t \phi_U \|^2 + \frac{1}{2} \| \nabla \mu_U \|^2 \, dt \\
\notag & \quad \leq \int_\Omega \frac{1}{2} |\nabla \phi_{0,k} |^2 + F(\phi_{0,k}) \, dx + \frac{1}{2}  \int_0^T \| \bm{U} \|^2 \, dt.
\end{align}
We now report on a continuity property for the solution pair to \eqref{viscousCH}, which can be deduced from \cite[Theorem A.1]{CHNS3Dstrsol}:
\begin{prop}\label{prop:viscousCH:cts}
Let $\bm{U}_n \in H^1(0,T;\VV_{\bu,m})$ be a sequence and $\bm{U} \in H^1(0,T;\VV_{\bu,m})$ such that $\bm{U}_n \to \bm{U}$ in $L^2(0,T;\VV_{\bu,m})$. Let $(\phi_n, \mu_n)$ and $(\phi_U, \mu_U)$ be the unique strong solutions to \eqref{viscousCH} corresponding to $\bm{U}_k$ and $\bm{U}$, respectively. Then, it holds that as $n \to \infty$,
\begin{subequations}\label{viscousCH:diff:est}
\begin{alignat}{2}
\| \phi_n - \phi_U \|_{L^\infty(0,T;H^2)}& \to 0, \\
\| F'(\phi_n) - F'(\phi_U) \|_{L^\infty(0,T;L^2)} &\to 0, \\
\| \mu_n - \mu_U \|_{L^\infty(0,T;L^2)} & \to 0.
\end{alignat}
\end{subequations}
Furthermore, there exist positive constants $C$ and $\delta \in (0,1)$ independent of $n$ such that 
\begin{subequations}\label{viscousCH:ctsdep:est}
\begin{alignat}{2}
& \| \phi_n \|_{L^\infty(0,T;H^3) \cap W^{1,\infty}(0,T;H^1) \cap H^1(0,T;H^2)} + \| \mu_n \|_{L^\infty(0,T;H^2) \cap H^1(0,T;L^2)} \leq C, \\
& \| \phi_U \|_{L^\infty(0,T;H^3) \cap W^{1,\infty}(0,T;H^1) \cap H^1(0,T;H^2)} + \| \mu_U \|_{L^\infty(0,T;H^2) \cap H^1(0,T;L^2)} \leq C, \\
& \max_{(x,t) \in \Omega \times (0,T)} |\phi_k(x,t)| \leq 1-\delta, \quad  \max_{(x,t) \in \Omega \times (0,T)} |\phi_U(x,t)| \leq 1-\delta.
\end{alignat}
\end{subequations}
\end{prop}
We remark that the injection $H^1(0,T;\VV_{\bu,m}) \subset L^\infty(0,T;L^3(\Omega)) \cap H^1(0,T;L^{\frac{6}{5}}(\Omega))$ enables us to employ Theorem \ref{Pre.result} to derive the estimates \eqref{viscousCH:ctsdep:est}. Next, we consider a Faedo--Galerkin ansatz (with an abuse of notation)
\[
\bu_m(x,t) := \sum_{j=1}^m a_j^m(t) \bm{Y}_j(x), \quad \bo_m(x,t) := \sum_{j=1}^m b_j^m(t) \bm{Z}_j(x)
\]
satisfying the following micropolar Navier--Stokes system
\begin{subequations}\label{Galerkin}
\begin{alignat}{2}
\label{Gal:u} & (\rho(\phi_U) \pd_t \bu_m, \bm{Y}_k) + (\rho(\phi_U) (\bm{U} \cdot \nabla) \bu_m, \bm{Y}_k) + (2 \eta(\phi_U) \D \bu_m, \D \bm{Y}_k) \\
\notag & \qquad + (2 \eta_r(\phi_U) \W \bu_m, \W \bm{Y}_k) - \tfrac{\overline{\rho}_1 - \overline{\rho}_2}{2} ((\nabla \mu_U \cdot \nabla) \bu_m, \bm{Y}_k) \\
\notag & \quad = (2 \curl(\eta_r(\phi_U) \bm{w}_m), \bm{Y}_k) + (\mu_U \nabla \phi_U,  \bm{Y}_k), \\[1ex]
\label{Gal:w}& (\rho(\phi_U) \pd_t \bo_m, \bm{Z}_k) + (\rho(\phi_U) (\bm{U} \cdot \nabla) \bo_m, \bm{Z}_k) + (c_0(\phi_U) \div \bo_m, \div \bm{Z}_k) \\
\notag & \qquad + (2c_d(\phi_U) \D \bo_m, \D \bm{Z}_k) + (2 c_a(\phi_U) \W \bo_m, \W \bm{Z}_k) - \tfrac{\overline{\rho}_1 - \overline{\rho}_2}{2} ((\nabla \mu_U \cdot \nabla) \bo_m, \bm{Z}_k) \\
\notag & \quad = (2 \eta_r(\phi_U) \curl \bm{u}_m, \bm{Z}_k) - (4 \eta_r(\phi_U) \bo_m, \bm{Z}_k)
\end{alignat}
\end{subequations}
for all $k =1, \dots, m$ with $\bu_m(0) = \PP_{\bu,m} \bu_0$ and $\bo_m(0) = \PP_{\bo,m} \bo_0$. Let us reformulate \eqref{Galerkin} into a system of ordinary differential equations by setting 
\[
\bm{a}^m(t) = (a_1^m(t), \dots, a_m^m(t))^{\top}, \quad \bm{b}^m(t) = (b_1^m(t), \dots, b_m^m(t))^{\top}
\]
and introducing associated matrices and vectors
\begin{align*}
\bm{M}_{\bu,j,k}(t) & := \int_\Omega \rho(\phi_U) \bm{Y}_j  \cdot \bm{Y}_k \, dx, \quad \bm{M}_{\bo,j,k}(t) := \int_\Omega \rho(\phi_U) \bm{Z}_j  \cdot \bm{Z}_k \, dx, \\
\bm{L}_{\bu,j,k}(t) & :=  \int_\Omega \rho(\phi_U) (\bm{U} \cdot \nabla) \bm{Y}_j \cdot \bm{Y}_k + 2 \eta(\phi_U) \D \bm{Y}_j : \D \bm{Y}_k \, dx \\
& \quad + \int_\Omega 2 \eta_r(\phi_U) \W \bm{Y}_j : \W \bm{Y}_k - \tfrac{\overline{\rho}_1 - \overline{\rho}_2}{2} (\nabla \mu_U \cdot \nabla) \bm{Y}_j \cdot \bm{Y}_k \, dx, \\
\bm{L}_{\bo,j,k}(t) & :=  \int_\Omega \rho(\phi_U) (\bm{U} \cdot \nabla) \bm{Z}_j \cdot \bm{Z}_k + c_0(\phi_U) \div \bm{Z}_j  \div \bm{Z}_k \, dx \\
& \quad + \int_\Omega 2 c_d(\phi_U) \D \bm{Z}_j : \D \bm{Z}_k + 2 c_a(\phi_U) \W \bm{Z}_j : \W \bm{Z}_k \, dx \\
& \quad - \int_\Omega  \tfrac{\overline{\rho}_1 - \overline{\rho}_2}{2} (\nabla \mu_U \cdot \nabla) \bm{Z}_j \cdot \bm{Z}_k - 4 \eta_r(\phi_U) \bm{Z}_j \cdot \bm{Z}_k \, dx, \\
\bm{H}_{\bu,j,k}(t) & := \int_\Omega 2 \curl(\eta_r(\phi_U) \bm{Z}_j) \cdot \bm{Y}_k \, dx, \quad \bm{H}_{\bo,j,k}(t) := \int_\Omega 2 \eta_r(\phi_U) \curl \bm{Y}_j \cdot \bm{Z}_k \, dx, \\
\bm{G}_{\bu,k}(t)& := \int_\Omega \mu_U \nabla \phi_U \cdot \bm{Y}_k \, dx,
\end{align*}
so that \eqref{Galerkin} can be expressed as
\begin{align*}
\bm{M}_{\bu} \frac{d}{dt} \bm{a}^m + \bm{L}_{\bu} \bm{a}^m & =  \bm{H}_{\bu} \bm{b}^m + \bm{G}_{\bu}, \\
\bm{M}_{\bo} \frac{d}{dt} \bm{b}^m + \bm{L}_{\bo} \bm{b}^m &= \bm{H}_{\bo} \bm{a}^m,
\end{align*}
subjected to the initial conditions 
\[
a_i^m(0) = (\PP_{\bu,m} \bu_0, \bm{Y}_i), \quad b_i^m(0) = (\PP_{\bo,m} \bo_0, \bm{Z}_i) \quad \text{ for } 1 \leq i \leq m.
\]
In light of the regularity for $(\bm{U}, \phi_U, \mu_U)$ we can deduce via standard Sobolev embeddings that 
\[
\phi_U \in C^0([0,T];W^{1,4}(\Omega)), \quad \mu_U \in C^0([0,T];H^1(\Omega)), \quad \bm{U} \in C^0([0,T];\HH_{\div}).
\]
Along with the continuity of $\rho$ and the viscosity functions $\eta$, $\eta_r$, $c_0$, $c_d$ and $c_a$, it holds that $\bm{M}_{\bu}$ and $\bm{M}_{\bo}$ are continuous and uniformly positive definite in $t \in [0,T]$. Hence, invoking the standard theory for ordinary differential systems we can deduce the existence of unique solutions $\bm{a}^m \in C^1([0,T];\R^m)$ and $\bm{b}^m \in C^1([0,T];\R^m)$ on $[0,T]$, leading to the existence of unique solutions $\bu_m \in C^1([0,T];\VV_{\bu,m})$ and $\bo_m \in C^1([0,T];\VV_{\bo,m})$ to \eqref{Galerkin}. 

\subsubsection{Estimates on Galerkin solutions}
Next, we derive some estimates in preparation for a fixed point argument. For convenience we use the notations $\rho_U := \rho(\phi_U)$, $\rho_* := \min(\overline{\rho}_1, \overline{\rho}_2)$, $\rho^* := \max(\overline{\rho}_1, \overline{\rho}_2)$, and likewise for the viscosity functions $\eta$, $\eta_r$, $c_0$, $c_{d}$ and $c_a$. Furthermore, we set
\[
E_{\GL}(\phi) := \int_\Omega \frac{1}{2} |\nabla \phi|^2 + F(\phi) \, dx
\]
as the Ginzburg--Landau functional. We further assume that $\bm{U}$ satisfies
\begin{align}\label{ass:U}
\int_0^s \| \bm{U}(r) \|^2 \, dr \leq C_3 e^{C_2 s} \quad \text{ for } s \in [0,T],
\end{align}
for some positive constants $C_2$ and $C_3$ whose values will be revealed later.

Multiplying \eqref{Gal:u} by $a_k^m$ and \eqref{Gal:w} by $b_k^m$, then summing over $k$ and adding the resulting equalities yields
\begin{align*}
& \int_\Omega \rho_U \Big (\frac{d}{dt} \frac{1}{2} | \bu_m|^2 + \bm{U} \cdot \nabla \frac{|\bu_m|^2}{2} \Big ) + \rho_U \Big (\frac{d}{dt} \frac{1}{2} | \bo_m|^2 + \bm{U} \cdot \nabla \frac{|\bo_m|^2}{2} \Big ) \, dx \\
& \qquad + \int_\Omega  2 \eta_U |\D \bu_m|^2  + c_{0,U} |\div \bo_m|^2 + 2 c_{d,U} |\D \bo_m|^2 +  2 c_{a,U} |\W \bo_m|^2 \, dx \\
& \qquad + \int_\Omega 4\eta_{r,U} \big | \tfrac{1}{2}\curl \bu_m  - \bo_m \big |^2 \, dx  \\
& \quad = \int_\Omega \mu_U \nabla \phi_U \cdot \bu_m  +\tfrac{\overline{\rho}_1 - \overline{\rho}_2}{2} \nabla \mu_U \cdot \nabla \Big ( \frac{|\bu_m|^2 + |\bo_m|^2}{2} \Big ) \, dx.
\end{align*}
In the above we employed the integration by parts formula involving the curl operator \eqref{IBP:curl}, as well as the identity
\[
2\W \bm{a} : \W \bm{a} = |\curl \bm{a}|^2.
\]
Using the properties $\rho'(\phi_U) = \frac{\overline{\rho}_1 - \overline{\rho}_2}{2}$ and $\div \bm{U} = 0$, by \eqref{viscousCH} it holds that 
\begin{align*}
- \int_\Omega (\pd_t \rho_U + \div (\rho_U \bm{U}) - \tfrac{\overline{\rho}_1- \overline{\rho}_2}{2} \Delta \mu_U) \Big( \frac{|\bu_m|^2 + |\bo_m|^2}{2} \Big ) \, dx  = 0.
\end{align*}
Thus, together with the identity
\[
\int_\Omega \rho_U \frac{d}{dt} \frac{1}{2} |\bu_m|^2 + \rho_U \bm{U} \cdot \nabla \frac{|\bu_m|^2}{2} \, dx = \frac{d}{dt} \int_\Omega \frac{\rho_U}{2} |\bu_m|^2 \, dx - \int_\Omega \frac{|\bu_m|^2}{2} (\pd_t \rho_U + \bm{U} \cdot \nabla \rho_U) \,dx,
\]
and likewise for $\bo_m$, by integrating by parts in space and in time we deduce that 
\begin{equation}\label{Gal:1}
\begin{aligned}
&\frac{d}{dt}  \int_\Omega \frac{\rho_U}{2} (|\bu_m|^2 + |\bo_m|^2) \, dx + \int_\Omega 4\eta_{r,U} \big | \tfrac{1}{2}\curl \bu_m  - \bo_m \big |^2 \, dx   \\
& \qquad + \int_\Omega  2 \eta_U |\D \bu_m|^2  + c_{0,U} |\div \bo_m|^2 + 2 c_{d,U} |\D \bo_m|^2 +  2 c_{a,U} |\W \bo_m|^2 \, dx \\
& \quad =- \int_\Omega \phi_U \nabla \mu_U \cdot \bu_m \, dx  \leq \eta_* \| \D \bu_m \|^2 + \frac{1}{4 \eta_* C_p^2} \| \nabla \mu_U \|^2,
\end{aligned}
\end{equation}
where we have used that $\| \phi_U \|_{L^\infty} \leq 1$ and the Poincar\'e inequality on $\bu_m$ with Poincar\'e constant $C_p$. Applying the Gr\"onwall inequality to \eqref{Gal:1} on $[0,s]$ with $s \in [0,T]$ and the estimate \eqref{vCH:2} for $\nabla \mu_U$ we find that 
\begin{equation}\label{Gal:2}
\begin{aligned} \frac{\rho_*}{2}  \Big (\| \bu_m(s) \|^2 + \| \bo_m(s) \|^2 \Big ) & \leq \frac{\rho^*}{2} \Big ( \| \PP_{\bu,m}\bu_0 \|^2 + \| \PP_{\bo,m} \bo_0 \|^2 \Big ) \\
& \quad + \frac{1}{2 \eta_* C_p^2} E_{\GL}(\phi_{0,k}) + \frac{1}{4\eta_* C_p^2} \int_0^s \| \bm{U} \|^2 \, dt.
\end{aligned}
\end{equation}
Defining 
\begin{align*}
C_1 := \frac{\rho^*}{\rho_*} \Big ( \| \bu_0 \|^2 + \| \bo_0 \|^2 \Big ) + \frac{1}{\eta_* \rho_* C_p^2} E_{\GL}(\phi_{0,k}), \quad C_2 := \frac{1}{2 \eta_* \rho_* C_p^2}, \quad C_3 := C_1 T,
\end{align*}
we infer from \eqref{Gal:2} and \eqref{ass:U} that for all $t \in [0,T]$:
\begin{align}
\label{Gal:3} & \sup_{t \in [0,T]} \Big ( \| \bu_m(t) \|^2 + \| \bo_m(t) \|^2 \Big ) \leq (C_1 + C_3 C_2 e^{C_2 T}) =: K_0^2, \\
\label{Gal:4} & \int_0^t \Big (\| \bu_m(s) \|^2 + \| \bo_m(s) \|^2\Big ) \, ds  \leq C_1 t + C_2 \int_0^t \int_0^s \| \bm{U}(r) \|^2 \, dr \, ds \leq C_3 e^{C_2 t}.
\end{align}
Next, multiplying \eqref{Gal:u} by $(a_k^m)'(t)$ and upon summing over $k$ we obtain
\begin{align*}
\rho_* \| \pd_t \bu_m \|^2 & \leq -(\rho_U (\bm{U} \cdot \nabla) \bu_m, \pd_t \bu_m) - (2\eta_U \D \bu_m, \D \pd_t \bu_m) - (2 \eta_{r,U} \W \bu_m, \W \pd_t \bu_m) \\
& \quad + \tfrac{\overline{\rho}_1 - \overline{\rho}_2}{2} ((\nabla \mu_U \cdot \nabla) \bu_m, \pd_t \bu_m) - (\phi_U \nabla \mu_U, \pd_t \bu_m) + (2 \curl(\eta_{r,U} \bo_m), \pd_t \bu_m) \\
& \leq \rho^* \| \bm{U} \| \| \nabla \bu_m \|_{L^\infty} \| \pd_t \bu_m \| +  \big (2 \eta^* \| \D \bu_m \| + 2 \eta_{r}^* \| \W \bu_m \| \big ) \| \nabla \pd_t \bu_m \| \\
& \quad + \Big ( |\tfrac{\overline{\rho}_1 - \overline{\rho}_2}{2}| \| \nabla \bu_m \|_{L^\infty} + \| \phi_U \|_{L^\infty} \Big ) \| \nabla \mu_U \| \| \pd_t \bu_m \| + 2 \eta_r^* \| \bo_m \| \| \nabla \pd_t \bu_m \|.
\end{align*}
In light of the inverse estimates \eqref{inverse:est} and the Sobolev embedding $H^2(\Omega) \subset L^\infty(\Omega)$ we see that 
\begin{equation}\label{Gal:5}
\begin{aligned}
\rho_* \| \pd_t \bu_m \|^2  & \leq \rho^* C_m \| \bm{U} \| \| \bu_m \| \| \pd_t \bu_m \| + 2(\eta^* + \eta_r^*) C_m^2 \| \bu_m \| \| \pd_t \bu_m \| \\
& \quad + \Big ( |\tfrac{\overline{\rho}_1 - \overline{\rho}_2}{2}| C_m \| \bu_m \| + 1 \Big ) \| \nabla \mu_U \| \| \pd_t \bu_m \| + 2 \eta_r^* C_m \| \bo_m \| \| \pd_t \bu_m \|.
\end{aligned}
\end{equation}
Invoking \eqref{vCH:2} for $\nabla \mu_U$, \eqref{ass:U} for $\bm{U}$, \eqref{Gal:3} and \eqref{Gal:4} for $\bu_m$ and $\bo_m$, we deduce from \eqref{Gal:5}
\begin{equation}\label{Gal:6}
\begin{aligned}
& \int_0^T \| \pd_t \bu_m(t) \|^2 \, dt \\
& \quad \leq 4 \Big ( \frac{\rho^* C_m}{\rho_*} K_0 \Big)^2 C_3 e^{C_2 T} + 4 \Big ( \frac{2 C_m^2 (\eta^* + \eta_r^*)}{\rho_*} \Big )^2 C_3 e^{C_2 T} \\
& \qquad + 4 \Big ( \frac{|\overline{\rho}_1 - \overline{\rho}_2|}{2 \rho_*} C_m K_0 + \frac{1}{\rho_*} \Big )^2 (2 E_{\GL}(\phi_{0,k}) + C_3 e^{C_2 T}) + 4 \Big ( \frac{\eta_r^* C_m}{\rho_*}\Big)^2 C_3 e^{C_2 T} \\
& \quad =: K_1^2.
\end{aligned}
\end{equation}
Likewise, upon multiplying \eqref{Gal:w} by $(b_k^m)'(t)$ and summing over $k$ leads to 
\begin{align*}
\rho_* \| \pd_t \bo_m \|^2 & \leq -(\rho_U( \bm{U} \cdot \nabla) \bo_m, \pd_t \bo_m) - (c_{0,U} \div \bo_m, \div (\pd_t \bo_m)) - (2c_{d,U} \D \bo_m, \D \pd_t \bo_m) \\
& \quad - (2 c_{a,U} \W \bo_m, \W \pd_t \bo_m) + \tfrac{\overline{\rho}_1 - \overline{\rho}_2}{2}(( \nabla \mu_U \cdot \nabla) \bo_m, \pd_t \bo_m) \\
& \quad + (2 \eta_{r,U} \curl \bu_m, \pd_t \bo_m) - (4 \eta_{r,U} \bo_m, \pd_t \bo_m) \\
& \leq \rho^* \| \bm{U} \| \| \nabla \bo_m \|_{L^\infty} \| \pd_t \bo_m \| + (c_{0}^* + 2 c_d^* + 2 c_a^*) \| \bo_m \|_{H^1} \| \pd_t \bo_m \|_{H^1} \\
& \quad + |\tfrac{\overline{\rho}_1 - \overline{\rho}_2}{2}| \| \nabla \bo_m \|_{L^\infty} \| \nabla \mu_U \|  \| \pd_t \bo_m \| + 2 \eta_r^* \| \bo_m \|_{H^1} \| \pd_t \bo_m \| + 4 \eta_r^* \| \bo_m \| \| \pd_t \bo_m \| \\
& \leq \rho^* C_m \| \bm{U} \| \| \bo_m \| \| \pd_t \bo_m \| + (c_0^* + 2c_d^* + 2c_a^*) C_m^2 \| \bo_m \| \| \pd_t \bo_m \| \\
& \quad +  |\tfrac{\overline{\rho}_1 - \overline{\rho}_2}{2}| C_m \| \bu_m \| \| \nabla \mu_U \| \| \pd_t \bo_m \| + ( 2 \eta_r^* C_m +  4 \eta_r^*) \| \bo_m \| \| \pd_t \bo_m \|.
\end{align*}
Invoking \eqref{vCH:2} for $\nabla \mu_U$, \eqref{ass:U} for $\bm{U}$, \eqref{Gal:3} and \eqref{Gal:4} for $\bu_m$ and $\bo_m$, we obtain
\begin{equation}\label{Gal:7}
\begin{aligned}
& \int_0^T \| \pd_t \bo_m(t) \|^2 \, dt \\
& \quad \leq 4 \Big ( \frac{\rho^* C_m}{\rho_*} K_0 \Big )^2 C_3 e^{C_2 T} + 4 \Big ( \frac{(c_0^* + 2 c_d^* + 2c_a^*) C_m^2}{\rho_*} \Big )^2 C_3 e^{C_2 T} \\
& \qquad + 4 \Big ( \frac{|\overline{\rho}_1 - \overline{\rho}_2|}{2\rho_*} C_m K_0\Big)^2 (2 E_{\GL}(\phi_{0,k}) + C_3 e^{C_2T}) + 4 \Big ( \frac{(2 \eta_r^* C_m + 4 \eta_r^*}{\rho_*} \Big )^2 C_3 e^{C_2 T} \\
& \quad =: K_2^2.
\end{aligned}
\end{equation}

\subsubsection{Fixed point argument}
Based on the estimates \eqref{Gal:6} and \eqref{Gal:7} we introduce the closed and convex set
\begin{align*}
S & := \Big \{ \bu \in H^1(0,T;\VV_{\bu,m}) \, \text{ such that } \\
& \qquad \, \int_0^t \| \bu(r) \|^2 \, dt \leq C_3 e^{C_2 t}, \, t \in [0,T] \text{ and }  \| \pd_t \bu \|_{L^2(0,T;\VV_{\bu,m})} \leq K_1\Big \},
\end{align*}
that is also compact in $L^2(0,T;\VV_{\bu,m})$, see \cite[Section 3, Theorem 1]{Simon}. Furthermore, we consider the mapping $\Lambda :S \to L^2(0,T;\VV_{\bu,m})$, where $\Lambda(\bm{U})$ is the first component of the unique solution pair $(\bu_m, \bo_m)$ satisfying \eqref{Gal:u}-\eqref{Gal:w} with variables $(\phi_U, \mu_U)$, themselves as the unique solution pair to the viscous convective Cahn--Hilliard equation \eqref{viscousCH} corresponding to $\bm{U}$.

To invoke the Schauder fixed point theorem and deduce the existence of a fixed point for $\Lambda$, which then implies the solvability of the approximate problem \eqref{approx:problem}, it suffices to establish the continuity of the mapping $\Lambda$.

We consider a sequence $\{\bm{U}_n\}_{n \in \N} \subset S$ and define $(\phi_n, \mu_n)$ as the corresponding unique solution pair to the viscous Cahn--Hilliard system \eqref{viscousCH}, as well as $(\bu_n, \bo_n)$ as the unique solution pair to \eqref{Gal:u}-\eqref{Gal:w} corresponding to $(\bm{U}_n, \phi_n, \mu_n)$.  Then, $\Lambda(\bm{U}_n) = \bu_n$.  Likewise we consider $\widetilde{\bm{U}} \in S$ with solution pairs $(\phi_{\widetilde{U}}, \mu_{\widetilde{U}})$ and $(\widetilde{\bu}, \widetilde{\bo})$ such that $\Lambda(\widetilde{\bm{U}}) = \widetilde{\bu}$. The differences
\[
\bm{U} := \bm{U}_n - \widetilde{\bm{U}}, \quad \bu := \bu_n - \widetilde{\bu}, \quad \bo := \bo_n - \widetilde{\bo}, \quad \phi := \phi_{n} - \phi_{\widetilde{U}}, \quad \mu := \mu_n - \mu_{\widetilde{U}}
\]
satisfy (using the notations $\rho_n := \rho(\phi_n)$,  $\widetilde{\rho} := \rho(\phi_{\widetilde{U}})$ and $\rho := \rho_n - \widetilde{\rho}$, likewise for $\eta$, $\eta_r$, $c_0$, $c_d$ and $c_a$)
\begin{equation}\label{Lam:ctsdep:u}
\begin{aligned}
& (\rho_n \pd_t \bu, \bm{Y}) + (\rho \pd_t \widetilde{\bu}, \bm{Y}) + (\rho_n(\bm{U}_n \cdot \nabla) \bu_n - \widetilde{\rho} (\widetilde{\bm{U}} \cdot \nabla) \widetilde{\bu}, \bm{Y}) \\
& \qquad + (2 \eta_n \D \bu, \D \bm{Y}) + (2 \eta \D \widetilde{\bu}, \D \bm{Y}) + (2 \eta_{r,n} \W \bu, \W \bm{Y}) + (2 \eta_r \W \widetilde{\bu}, \W \bm{Y}) \\
& \qquad  - \tfrac{\overline{\rho}_1 - \overline{\rho}_2}{2} ((\nabla \mu_n \cdot \nabla) \bu_n - (\nabla \mu_{\widetilde{U}} \cdot \nabla) \widetilde{\bu}, \bm{Y}) \\
& \quad = (\mu_n \nabla \phi_n - \mu_{\widetilde{U}} \nabla \phi_{\widetilde{U}} , \bm{Y}) + 2( \curl(\eta_{r,n} \bo_n) - \curl( \widetilde{\eta_r} \widetilde{\bo}), \bm{Y}),
\end{aligned}
\end{equation}
for all $\bm{Y} \in \VV_{\bu,m}$ and for all $t \in [0,T]$, as well as
\begin{equation}\label{Lam:ctsdep:w}
\begin{aligned}
& (\rho_n \pd_t \bo, \bm{Z}) + (\rho \pd_t \widetilde{\bo}, \bm{Z}) + (\rho_n(\bm{U}_n \cdot \nabla) \bo_n - \widetilde{\rho} (\widetilde{\bm{U}} \cdot \nabla) \widetilde{\bo}, \bm{Z}) \\
& \qquad + (c_{0,n} \div \bo, \div \bm{Z}) + (c_0 \div \widetilde{\bo}, \div \bm{Z}) + (2 c_{d,n} \D \bo, \D \bm{Z}) + (2 c_d \D \widetilde{\bo}, \D \bm{Z}) \\
& \qquad + (2 c_{a,n} \W \bo, \W \bm{Z}) + (2 c_a \W \widetilde{\bo}, \W \bm{Z}) - \tfrac{\overline{\rho}_1 - \overline{\rho}_2}{2} ((\nabla \mu_n \cdot \nabla) \bo_n - (\nabla \mu_{\widetilde{U}} \cdot \nabla) \widetilde{\bo}, \bm{Z}) \\
& \quad = 2((\eta_{r,n} \curl \bu_n - \widetilde{\eta_{r}} \curl \bu), \bm{Z}) - 4(\eta_{r,n} \bo + \eta_r \widetilde{\bo}, \bm{Z})
\end{aligned}
\end{equation}
for all $\bm{Z} \in \VV_{\bo,m}$ and for all $t \in [0,T]$. Choosing $\bm{Y} = \bu \in \VV_{\bu,m}$ in \eqref{Lam:ctsdep:u} leads to 
\begin{equation}\label{Lam:cts:u:equ}
\begin{aligned}
& \frac{1}{2} \frac{d}{dt} \int_\Omega \rho_n |\bu|^2 \, dx + \int_\Omega 2 \eta_n |\D \bu|^2 + 2 \eta_{r,n} |\W \bu|^2 \, dx \\
& \quad = \int_\Omega \frac{\overline{\rho}_1 - \overline{\rho}_2}{4} \pd_t \phi_n |\bu|^2 - \frac{\overline{\rho}_1 - \overline{\rho}_2}{2} \phi \pd_t \widetilde{\bu} \cdot \bu - (\rho_n (\bm{U}_n \cdot \nabla) \bu_n - \widetilde{\rho}( \widetilde{\bm{U}} \cdot \nabla) \widetilde{\bu}) \cdot \bu \, dx \\
& \qquad - \int_\Omega 2 \eta \D \widetilde{\bu} : \D \bu + 2 \eta_r \W \widetilde{\bu} : \W \bu \, dx \\
& \qquad +  \int_\Omega \frac{\overline{\rho}_1 - \overline{\rho}_2}{2} ((\nabla \mu_n \cdot \nabla) \bu_n - (\nabla \mu_{\widetilde{U}} \cdot \nabla) \widetilde{\bu}) \cdot \bu + (\mu_n \nabla \phi_n - \mu_{\widetilde{U}} \nabla \phi_{\widetilde{U}}) \cdot \bu \, dx \\
& \qquad + \int_\Omega 2 \curl(\eta_r \bo_n - \widetilde{\eta_r} \widetilde{\bo})  \cdot \bu \, dx \\
& =: A_1 + \cdots + A_8.
\end{aligned}
\end{equation}
Employing the Korn inequality \eqref{Korn.ineq.}, the inverse estimates \eqref{inverse:est}, the estimates \eqref{viscousCH:ctsdep:est} and \eqref{Gal:3} and the fact that $\widetilde{\bu} \in S$, we see that 
\begin{align*}
|A_1|& \leq C \| \pd_t \phi_n \|_{L^6} \| \bu \| \| \bu \|_{L^3} \leq \delta \| \D \bu \|^2 + C_\delta \| \bu \|^2, \\
|A_2| & \leq C \| \phi \|_{L^\infty} \| \pd_t \widetilde{\bu} \| \| \bu \| \leq C \| \bu \|^2 + C \| \pd_t \widetilde{\bu} \|^2 \| \phi \|_{H^2}^2, \\
|A_4|& \leq C \| \phi \|_{L^\infty} \| \D \widetilde{\bu} \| \| \D \bu \| \leq \delta \| \D \bu \|^2 + C_\delta \| \D \widetilde{\bu} \|^2 \| \phi \|_{H^2}^2 \leq \delta \| \D \bu \|^2 + C_{\delta,m} \| \phi \|_{H^2}^2,  \\
|A_5|& \leq C \| \phi \|_{L^\infty} \| \W \widetilde{\bu} \| \| \W \bu \| \leq \delta \| \W \bu \|^2 + C_\delta \| \W \widetilde{\bu} \|^2 \| \phi \|_{H^2}^2 \leq \delta \| \W \bu \|^2 + C_{\delta,m} \| \phi \|_{H^2}^2, \\
|A_7| & \leq (\| \mu \| \| \nabla \phi_n \|_{L^6} + \| \mu_{\widetilde{U}} \| \| \nabla \phi \|_{L^6}) \| \bu \|_{L^3} \leq \delta \| \D \bu \|^2 + C_\delta( \| \mu \|^2 + \| \phi \|_{H^2}^2),  \\
|A_8| & \leq C (\| \bo \| +  \| \phi \|_{L^\infty} \| \widetilde{\bo} \|) \| \curl \bu \| \leq \delta \| \W \bu \|^2 + C_\delta ( \| \bo \|^2 + \| \phi \|_{H^2}^2).
\end{align*}
On the other hand for $A_3$ and $A_6$ we estimate with the help of the inverse estimates \eqref{inverse:est} and the fact that $\bm{U}_n, \bu_n, \widetilde{\bu} \in S$ as follows:
\begin{align*}
|A_3| & \leq |(\tfrac{\overline{\rho}_1 - \overline{\rho}_2}{2} \phi (\bm{U}_n \cdot \nabla) \bu_n, \bu)| + |(\widetilde{\rho}( \bm{U} \cdot \nabla) \bu_n, \bu)| + |(\widetilde{\rho}( \widetilde{\bm{U}} \cdot \nabla) \bu, \bu)| \\
& \leq C \| \phi \|_{L^\infty} \| \bm{U}_n \|_{L^\infty} \| \nabla \bu_n \| \| \bu \| + C \| \bm{U} \| \| \nabla \bu_n \|_{L^\infty} \| \bu \| + C \| \widetilde{\bm{U}} \|_{L^\infty} \| \nabla \bu \| \| \bu \| \\
& \leq C_m \| \phi \|_{H^2} \| \bu_n \| \| \bu \| + C_m \| \bm{U} \| \| \bu_n \| \| \bu \| + C \| \nabla \bu \| \| \bu \| \\
& \leq \delta \| \D \bu \|^2 + C_{\delta,m}( \| \bu \|^2 + \| \phi \|_{H^2}^2 + \| \bm{U} \|^2 ), 
\end{align*}
and
\begin{align*}
|A_6| & \leq |\tfrac{\overline{\rho}_1 - \overline{\rho}_2}{2}| |(\mu_n \Delta \bu_n - \mu_{\widetilde{U}} \Delta \widetilde{\bu}) , \bu)| + |\tfrac{\overline{\rho}_1 - \overline{\rho}_2}{2}| |(\mu_n \nabla \bu_n - \mu_{\widetilde{U}} \nabla \widetilde{\bu}, \nabla \bu)| \\
& = |\tfrac{\overline{\rho}_1 - \overline{\rho}_2}{2}| |(\mu \Delta \bu_n - \mu_{\widetilde{U}} \Delta \bu) , \bu)| + |\tfrac{\overline{\rho}_1 - \overline{\rho}_2}{2}| |(\mu \nabla \bu_n - \mu_{\widetilde{U}} \nabla \bu, \nabla \bu)| \\
& \leq C \| \mu \| \| \Delta \bu_n \| \| \bu \|_{L^\infty} + C \| \mu_{\widetilde{U}} \|_{L^6} \| \Delta \bu \| \| \bu \|_{L^3} \\
& \quad + C \| \mu \| \| \nabla \bu_n \|_{L^6} \| \nabla \bu \|_{L^3} + C \| \mu_{\widetilde{U}} \|_{L^6} \| \nabla \bu \| \| \nabla \bu \|_{L^3} \\
& \leq C_m \| \mu \| \| \nabla \bu \| + C_m \| \nabla \bu \| \| \bu \| \\
& \leq \delta \| \D \bu \|^2 + C_{\delta,m} (\| \mu \|^2 + \| \bu \|^2).
\end{align*}
Meanwhile, choosing $\bm{Z} = \bo \in \VV_{\bo,m}$ in \eqref{Lam:ctsdep:w} yields
\begin{equation}\label{Lam:cts:w:equ}
\begin{aligned}
& \frac{1}{2} \frac{d}{dt} \int_\Omega \rho_n |\bo|^2 \, dx + \int_\Omega c_{0,n} |\div \bo|^2 + 2c_{d,n} | \D \bo|^2 + 2 c_{a,n} |\W \bo|^2  + 4 \eta_{r,n} |\bo|^2 \, dx \\
& \quad = \int_\Omega   \frac{\overline{\rho}_1 - \overline{\rho}_2}{4} \pd_t \phi |\bo|^2 - \frac{\overline{\rho}_1 - \overline{\rho}_2}{2} \phi \pd_t \widetilde{\bo} \cdot \bo - (\rho_n (\bm{U}_n \cdot \nabla) \bo_n - \widetilde{\rho}(\widetilde{\bm{U}} \cdot \nabla) \widetilde{\bo}) \cdot \bo \, dx \\
& \qquad - \int_\Omega [c_0 (\div \widetilde{\bo}) (\div \bo) + 2 c_d \D \widetilde{\bo} : \D \bo + 2 c_a \W \widetilde{\bo} : \W \bo - 4 \eta_r \widetilde{\bo} \cdot \bo] \, dx \\
& \qquad +  \int_\Omega \frac{\overline{\rho}_1 - \overline{\rho}_2}{2} ((\nabla \mu_{U_n} \cdot \nabla) \bo_n - (\nabla \mu_{\widetilde{U}} \cdot \nabla) \widetilde{\bo}) \cdot \bo \, dx \\
& \qquad + \int_\Omega 2 (\eta_{r,n} \curl \bu_n - \widetilde{\eta_r} \curl \widetilde{\bu}) \cdot \bo \, dx \\
& \quad =: A_9 + \cdots + A_{14}.
\end{aligned}
\end{equation} 
Via similar calculations we find that 
\begin{align*}
|A_9|& \leq \delta \| \D \bo \|^2 + C_{\delta} \| \bo \|^2, \\
|A_{10}| & \leq C \| \bo \|^2 + C \| \pd_t \widetilde{\bo} \|^2 \| \phi \|_{H^2}^2, \\
|A_{11}| & \leq \delta \| \D \bo \|^2 + C_{\delta,m} (\| \bo \|^2 + \| \phi \|_{H^2}^2 + \| \bm{U} \|^2), \\
|A_{12}|& \leq \delta (\| \div \bo \|^2 + \| \D \bo \|^2 + \| \W \bo \|^2 + \| \bo \|^2) + C_{\delta,m} \|  \phi \|_{H^2}^2, \\
|A_{13}| & \leq \delta \| \D \bo \|^2 + C_{\delta,m} (\| \mu \|^2 + \| \bo \|^2), \\
|A_{14}|& \leq  C(\| \curl \bu \| + \| \phi \|_{H^2} \| \curl \widetilde{\bu} \|) \| \bo \| \leq \delta \| \W \bu \|^2 + C_{\delta,m}( \| \phi \|_{H^2}^2 + \| \bo \|^2).
\end{align*}
Combining \eqref{Lam:cts:u:equ} and \eqref{Lam:cts:w:equ} while employing the above estimates with $\delta$ sufficiently small, we arrive at the differential inequality
\begin{align*}
& \frac{d}{dt} \int_\Omega \rho_n (|\bu|^2 + |\bo|^2) \, dx \\
& \quad \leq C \int_\Omega \rho_n(|\bu|^2 + |\bo|^2) \, dx + C (1 +  \| \pd_t \widetilde{\bu} \|^2 + \| \pd_t \widetilde{\bo} \|^2) \| \phi \|_{H^2}^2 + C \| \mu \|^2 + C \| \bm{U} \|^2.
\end{align*}
Application of Gr\"onwall's inequality then yields
\begin{align*}
\sup_{t \in (0,T]}( \| \bu(t) \|^2 + \| \bo(t) \|^2)& \leq \frac{C\exp(T)}{\rho_*}  \int_0^T(1 +  \| \pd_t \widetilde{\bu} \|^2 + \| \pd_t \widetilde{\bo} \|^2) \| \phi \|_{H^2}^2 + \| \mu \|^2 + \| \bm{U} \|^2 \, dr \\
& \leq C \big ( \| \phi \|_{L^\infty(0,T;H^2)}^2 + \| \mu \|_{L^\infty(0,T;L^2)}^2 + \| \bm{U} \|_{L^2(0,T;L^2)}^2 \big ).
\end{align*}
On account of the stability estimates \eqref{viscousCH:diff:est} we see that if $\bm{U}_n \to \widetilde{\bm{U}}$ in $L^2(0,T;\VV_{\bu,m})$, then $\bu_n \to \widetilde{\bu}$ in $L^\infty(0,T;\VV_{\bu,m})$ and $\bo_n \to \widetilde{\bo}$ in $L^\infty(0,T;\VV_{\bo,m})$.  This in turn implies the continuity of the mapping $\Lambda$, and hence by Schauder's fixed point theorem there exists a fixed point of $\Lambda$ in $S$, leading to the existence of approximate solutions $(\bu_m, \bo_m, \phi_m, \mu_m)$ on $[0,T]$ to \eqref{approx:problem} for any $m \in \N$.

\subsection{Uniform estimates}
In this section the symbol $C$ denotes positive constants independent of $m \in \N$, $\alpha \in (0,1]$ and $k > \overline{k}$, whose values may change line by line and also within the same line.
\subsubsection{Standard estimates}
Integrating \eqref{app:phi} over $\Omega$ and employing $\div \bu_m = 0$ and the Neumann condition for $\mu_m$ yields the conservation of mass property
\begin{align}\label{uni:mass}
\int_\Omega \phi_m(t) \, dx = \int_\Omega \phi_{0,k} \, dx \quad \forall t \in (0,T].
\end{align}
Next, we take $\bm{v}= \bu_m$ in \eqref{app:u}, $\bm{z} = \bo_m$ in \eqref{app:w} and testing \eqref{app:phi} with $\mu_m$ and \eqref{app:mu} with $\pd_t \phi_m$, upon summing we obtain the energy identity analogous to \eqref{intro:Ene.ineq}:
\begin{equation}\label{uni:energy}
\begin{aligned}
\frac{d}{dt} E(\bu_m, \bo_m, \phi_m) &+ \int_\Omega |\nabla \mu_m|^2 + 2 \eta(\phi_m) |\D \bu_m|^2 + 4 \eta_r(\phi_m) |\tfrac{1}{2} \curl \bu_m - \bo_m|^2 \, dx \\
&  + \int_\Omega c_0(\phi_m) |\div \bo_m|^2 + 2 c_d(\phi_m) |\D \bo_m|^2 + 2 c_a(\phi_m) |\W \bo_m|^2 \, dx \\
&  + \int_\Omega \alpha |\pd_t \phi_m|^2 \, dx = 0,
\end{aligned}
\end{equation}
where $E$ is defined as in \eqref{defn:E}. Due to \eqref{prop.approx.IC.delta} for $\phi_{0,k}$, it holds that 
\[
\int_\Omega F(\phi_{0,k}) \, dx \leq |\Omega| \max_{s \in [-1,1]} |F(s)|,
\]
using the continuity of $F$ over $[-1,1]$. Hence, upon integrating \eqref{uni:energy} over $[0,T]$ and employing the properties $\| \PP_{\bu,m} \bu_0 \| \leq \| \bu_0 \|$ and $\| \PP_{\bo,m} \bo_0 \| \leq \| \bo_0 \|$ we find there exists a positive constant $C$ independent of $m$, $\alpha$ and $k$ such that 
\begin{subequations}\label{uni:est:1}
\begin{alignat}{2}
\| \bu_m \|_{L^\infty(0,T;\HH_{\div})} + \| \bu_m \|_{L^2(0,T;\HH^1)} & \leq C, \\
\| \bo_m \|_{L^\infty(0,T;\HH)} + \| \bo_m \|_{L^2(0,T;\HH^1)} & \leq C, \\
\| \phi_m \|_{L^\infty(0,T;H^1)} & \leq C, \\
\| \nabla \mu_m \|_{L^2(0,T;L^2)} & \leq C, \\
\| F(\phi_m) \|_{L^\infty(0,T;L^1)} & \leq C, \\
\sqrt{\alpha} \| \pd_t \phi_m \|_{L^2(0,T;L^2)} & \leq C.
\end{alignat}
\end{subequations}
Next, testing \eqref{app:phi} with $- \Delta \phi_m$ and integrating by parts, then using the property $F''(s) \geq - \theta_0$ obtained from the convexity of the logarithmic part, we find that
\begin{equation}\label{Lapphi:nablamu}
\begin{aligned}
\frac{1}{2}\| \Delta \phi_m \|^2  & \leq \frac{\alpha^2}{2} \| \pd_t \phi_m \|^2 +  \theta_0 \| \nabla \phi_m \|^2 + \| \nabla \mu_m \| \| \nabla \phi_m \|\\
& \leq C (1 + \| \nabla \mu_m \|^2 + \alpha^2 \| \pd_t \phi_m \|^2),
\end{aligned}
\end{equation}
and this yields after integrating and then employing elliptic regularity theory
\begin{align}\label{uni:est:2}
\| \phi_m \|_{L^2(0,T;H^2)} \leq C.
\end{align}
We now invoke the following inequality valid for the convex part of the logarithmic potential $\Psi(s) = (1+s) \log(1+s) + (1-s) \log(1-s)$: For any $u \in (-1,1)$ there exist constants $c_u > 0$ and $C_u \geq 0$ such that 
\[
|\Psi'(s)| \leq c_u \Psi'(s)(s-u) + C_u \quad \forall s \in (-1,1).
\]
Choosing $u = \overline{\phi_{0,k}} \in (-1,1)$ and testing \eqref{app:mu} with $\phi_{m} - \overline{\phi_{0,k}}$ yields after applying the Poincar\'e inequality:
\begin{align*}
\|\nabla \phi_m \|^2 + \frac{1}{c_u} \| \Psi'(\phi_m) \|_{L^1} \leq C\alpha \| \pd_t \phi_m \| \|\nabla \phi_m \| + C\theta_0 \| \phi_m \| \| \nabla \phi_m \| + C \| \nabla \mu_m \| \| \nabla \phi_m \| + \frac{C_u}{c_u},
\end{align*}
and so on account of \eqref{uni:est:1}
\begin{align*}
\| F'(\phi_m) \|_{L^1} \leq C (1 + \sqrt{\alpha} \| \pd_t \phi_m \| + \| \nabla \mu_m \|)
\end{align*}
and then by integrating \eqref{app:mu} we infer the estimate for the mean value:
\begin{align}\label{mean:est}
|\overline{\mu_m}| \leq C (1 + \sqrt{\alpha} \| \pd_t \phi_m \| + \| \nabla \mu_m \|).
\end{align}
Hence, by the Poincar'e inequality we deduce
\begin{align}\label{uni:est:3}
\| \mu_m \| \leq C (1 + \sqrt{\alpha} \| \pd_t \phi_m \| + \| \nabla \mu_m \|) \quad \implies \quad 
\| \mu_m \|_{L^2(0,T;L^2)} \leq C,
\end{align}
and after testing \eqref{app:mu} with $F'(\phi_m)$:
\[
\frac{d}{dt} \int_\Omega \alpha F(\phi_m) \, dx + \| F'(\phi_m) \|^2 \leq \| \mu_m \| \| F'(\phi_m) \| + \theta_0 \| \phi_m \| \| F'(\phi_m) \|
\]
we also obtain
\begin{align}\label{uni:est:F'}
\| F'(\phi_m) \|_{L^2(0,T;L^2)} \leq C.
\end{align}
Lastly, testing \eqref{app:phi} with an arbitrary test function $\zeta \in H^1(\Omega)$ it follows that 
\begin{align*}
\| \pd_t \phi_m\|_{(H^1)^*} \leq \| \bu_m \| \| \nabla \phi_m \|_{L^3} + \| \nabla \mu_m \| \leq C(\| \phi_m \|_{H^2} + \| \nabla \mu_m \|),
\end{align*}
which leads to 
\begin{align}\label{uni:est:4}
\| \pd_t \phi_m \|_{L^2(0,T;(H^1)^*)} \leq C.
\end{align}

\subsubsection{High order estimates}
Choosing $\bm{v} = \pd_t \bu_m$ in \eqref{app:u} gives
\begin{align*}
& \frac{d}{dt} \int_\Omega  \eta(\phi_m) |\D \bu_m|^2 + \eta_r(\phi_m) |\W \bu_m|^2  \, dx  + \int_\Omega \rho(\phi_m) |\pd_t \bu_m|^2 \, dx \\
& \quad = - (\rho(\phi_m) (\bu_m \cdot \nabla) \bu_m, \pd_t \bu_m) + (\eta'(\phi_m) \pd_t \phi_m, |\D \bu_m|^2) \\
& \qquad + ( \eta_r'(\phi_m) \pd_t \phi_m, |\W \bu_m|^2) + \tfrac{\overline{\rho}_1 -\overline{\rho}_2}{2} ((\nabla \mu_m \cdot \nabla) \bu_m, \pd_t \bu_m) \\
& \qquad + (\mu_m \nabla \phi_m, \pd_t \bu_m) + (2 \curl(\eta_r(\phi_m) \bo_m), \pd_t \bu_m),
\end{align*}
while choosing $\bm{z} = \pd_t \bo_m$ in \eqref{app:w} gives
\begin{align*}
&\frac{d}{dt} \int_\Omega \tfrac{1}{2} c_0(\phi_m) |\div \bo_m|^2 + c_d(\phi_m) |\D \bo_m|^2 + c_a(\phi_m) |\W \bo_m|^2 \, dx  + \int_\Omega \rho(\phi_m) |\pd_t \bo_m|^2 \, dx \\
& \quad = (\tfrac{1}{2} c_0'(\phi_m) |\div \bo_m|^2 + c_d'(\phi_m) |\D \bo_m|^2 + c_a'(\phi_m) |\W \bo_m|^2, \pd_t \phi_m) \\
& \qquad  -(\rho(\phi_m) (( \bu_m \cdot \nabla) \bo_m, \pd_t \bo_m)  + \tfrac{\overline{\rho}_1 - \overline{\rho}_2}{2} ((\nabla \mu_m \cdot \nabla) \bo_m, \pd_t \bo_m) \\
& \qquad +  (2\eta_r(\phi_m) \curl \bu_m, \pd_t \bo_m) - (4\eta_r(\phi_m) \bo_m, \pd_t \bo_m).
\end{align*}
We note the last term of the former identity and the last two terms of the latter identity can be expressed together as
\begin{align*}
& (2 \eta_r(\phi_m) \bo_m, \curl \pd_t \bu_m) + (2 \eta_r(\phi_m) \curl \bu_m, \pd_t \bo_m) - (4 \eta_r(\phi_m) \bo_m, \pd_t \bo_m) \\
& \quad = \frac{d}{dt} \int_\Omega  2 \eta_r(\phi_m) \curl \bu_m \cdot \bo_m - 2 \eta_r(\phi_m) |\bo_m|^2 \, dx \\
& \qquad - (2 \eta_r'(\phi_m) \pd_t \phi_m, \curl \bu_m \cdot \bo_m) + (2 \eta_r'(\phi_m) \pd_t \phi_m, |\bo_m|^2).
\end{align*} 
Using the following identity deduce from the second identity in \eqref{eps:prop} of the Levi-Civita tensor:
\begin{align}\label{W:W:curl}
\W \bm{a} : \W \bm{b} = \frac{1}{2} \curl \bm{a} \cdot \curl \bm{b},
\end{align}
we find that upon summing 
\begin{equation}\label{reg:main}
\begin{aligned}
& \frac{d}{dt} \int_\Omega  \eta(\phi_m) |\D \bu_m|^2 + 2\eta_r(\phi_m) |\tfrac{1}{2} \curl \bu_m - \bo_m|^2 \, dx \\
& +\frac{d}{dt} \int_\Omega \tfrac{1}{2} c_0(\phi_m) |\div \bo_m|^2 + c_d(\phi_m) |\D \bo_m|^2 + c_a(\phi_m) |\W \bo_m|^2 \, dx \\
& \qquad  + \int_\Omega \rho(\phi_m) |\pd_t \bu_m|^2 + \rho(\phi_m) |\pd_t \bo_m|^2  \, dx \\
& \quad = - (\rho(\phi_m) (\bu_m \cdot \nabla) \bu_m, \pd_t \bu_m) + (\eta'(\phi_m) \pd_t \phi_m, |\D \bu_m|^2) \\
& \qquad + ( \eta_r'(\phi_m) \pd_t \phi_m, |\W \bu_m|^2) + \tfrac{\overline{\rho}_1 -\overline{\rho}_2}{2} ((\nabla \mu_m \cdot \nabla) \bu_m, \pd_t \bu_m) + (\mu_m \nabla \phi_m, \pd_t \bu_m) \\
& \qquad +  (\tfrac{1}{2} c_0'(\phi_m) |\div \bo_m|^2 + c_d'(\phi_m) |\D \bo_m|^2 + c_a'(\phi_m) |\W \bo_m|^2, \pd_t \phi_m) \\
& \qquad  -(\rho(\phi_m) (( \bu_m \cdot \nabla) \bo_m, \pd_t \bo_m)  + \tfrac{\overline{\rho}_1 - \overline{\rho}_2}{2} ((\nabla \mu_m \cdot \nabla) \bo_m, \pd_t \bo_m) \\
& \qquad  - (2 \eta_r'(\phi_m) \pd_t \phi_m, \curl \bu_m \cdot \bo_m) + (2 \eta_r'(\phi_m) \pd_t \phi_m, |\bo_m|^2) \\
& \quad =: B_1 + \cdots + B_{10}.
\end{aligned}
\end{equation}
From Proposition \ref{prop:Galerkin} we have $\mu_m \in L^\infty(0,T;H^2_n(\Omega)) \cap H^1(0,T;L^2(\Omega))$, so that upon testing \eqref{app:phi} by $\pd_t \mu_m$ and simplifying using the time derivative of \eqref{app:mu} yields
\begin{align*}
& \frac{d}{dt} \Big ( \int_\Omega \frac{1}{2} |\nabla \mu_m|^2 + \frac{\alpha}{2} |\pd_t \phi_m|^2 + \mu_m \bu_m \cdot \nabla \phi_m \, dx \Big ) + \| \nabla \pd_t \phi_m \|^2 \\
& \quad \leq \theta_0 \| \pd_t \phi_m \|^2 + (\mu_m \pd_t \bu_m , \nabla \phi_m) + (\mu_m \bu_m, \nabla \pd_t \phi_m ) \\
& \quad =: B_{11} + B_{12} + B_{13}.
\end{align*}
We define
\begin{equation}\label{defn:H}
\begin{aligned}
H_m(t) & = \int_\Omega \eta(\phi_m) |\D \bu_m|^2 + 2\eta_r(\phi_m) |\tfrac{1}{2} \curl \bu_m - \bo_m |^2 \, dx \\
& \quad + \int_\Omega  \tfrac{1}{2} c_0(\phi_m) |\div \bo_m|^2 + c_d(\phi_m) |\D \bo_m|^2 + c_a(\phi_m) |\W \bo_m|^2 \, dx \\
& \quad + \int_\Omega \frac{1}{2} |\nabla \mu_m|^2 + \frac{\alpha}{2}|\pd_t \phi_m|^2 + \mu_m \bu_m \cdot \nabla \phi_m \, dx,
\end{aligned}
\end{equation}
where by the Sobolev inequality \eqref{L3estimate} and the uniform estimates \eqref{uni:est:1} and \eqref{mean:est},
\begin{align*}
|(\mu_m \bu_m, \nabla \phi_m)| & \leq \| \mu_m - \overline{\mu_m} \|_{L^6} \| \bu_m \|_{L^3} \| \nabla \phi_m \| + |\overline{\mu_m}| \| \bu_m \| \| \nabla \phi_m \| \\
& \leq C\| \mu_m - \overline{\mu_m} \|_{H^1} \| \bu_m \|^{1/2} \| \nabla \bu_m \|^{1/2} + C(1 + \sqrt{\alpha} \| \pd_t \phi_m \| + \| \nabla \mu_m \|) \\
& \leq C \| \nabla \mu_m \| \| \D \bu_m \|^{1/2} + C(1 + \sqrt{\alpha} \| \pd_t \phi_m \| + \| \nabla \mu_m \|)  \\
& \leq  \int_\Omega \frac{1}{2}\eta(\phi_m) | \D \bu_m|^2 +  \frac{1}{4} |\nabla \mu_m|^2 + \frac{\alpha}{4} |\pd_t \phi_m|^2 \, dx + C,
\end{align*}
and so we have the lower bound for $H_m$:
\begin{equation}\label{Hm:lb}
\begin{aligned}
H_m & \geq  \int_\Omega \frac{1}{2}\eta(\phi_m) |\D \bu_m|^2 + 2\eta_r(\phi_m) |\tfrac{1}{2} \curl \bu_m - \bo_m |^2 \, dx \\
& \quad + \int_\Omega  \tfrac{1}{2} c_0(\phi_m) |\div \bo_m|^2 + c_d(\phi_m) |\D \bo_m|^2 + c_a(\phi_m) |\W \bo_m|^2 \, dx \\
& \quad + \int_\Omega \frac{1}{4} |\nabla \mu_m|^2 + \frac{\alpha}{4}|\pd_t \phi_m|^2 - C_0,
\end{aligned}
\end{equation}
for some positive constant $C_0$ independent of $m$. On the other hand we have the upper bound
\begin{equation}\label{Hm:ub}
\begin{aligned}
H_m & \leq  \int_\Omega \eta(\phi_m) |\D \bu_m|^2 + 2\eta_r(\phi_m) |\tfrac{1}{2} \curl \bu_m - \bo_m |^2 \, dx \\
& \quad + \int_\Omega  \tfrac{1}{2} c_0(\phi_m) |\div \bo_m|^2 + c_d(\phi_m) |\D \bo_m|^2 + c_a(\phi_m) |\W \bo_m|^2 \, dx \\
& \quad + \int_\Omega  |\nabla \mu_m|^2 + \alpha |\pd_t \phi_m|^2  \, dx + C.
\end{aligned}
\end{equation}
Recalling the Stokes operator $\bm{A}$, there exists $p_m \in C^0([0,T];H^1(\Omega))$ such that $- \Delta \bu_m + \nabla p_m = \bm{A} \bu_m$ a.e.~in $\Omega \times (0,T)$ and satisfies (see \cite{CHNS2Dstrsol,AGAMRT})
\begin{align}\label{pi:est}
\| p_m \| \leq C \| \nabla \bu_m \|^{1/2}\| \bm{A} \bu_m \|^{1/2}, \quad \| p_m \|_{H^1} \leq C \| \bm{A} \bu_m \|.
\end{align}
Then, when considering $\bm{v} = \bm{A} \bu_m$ in \eqref{app:u} we have
\begin{equation}\label{uni:est:6}
\begin{aligned}
& ((\eta(\phi_m) + \eta_r(\phi_m)) \bm{A} \bu_m, \bm{A} \bu_m) \\
& \quad = -(\rho(\phi_m) \pd_t \bu_m, \bm{A} \bu_m) - (\rho(\phi_m) (\bu_m \cdot \nabla) \bu_m, \bm{A}\bu_m) \\
& \qquad + \tfrac{\overline{\rho}_1 - \overline{\rho}_2}{2} ((\nabla \mu_m \cdot \nabla) \bu_m, \bm{A} \bu_m) + (\mu_m \nabla \phi_m, \bm{A} \bu_m) \\
& \qquad + ( 2 \eta'(\phi_m) \nabla \phi_m \D \bu_m , \bm{A} \bu_m) + (2\eta_r'(\phi_m) \nabla \phi_m \W \bu_m, \bm{A} \bu_m) \\
& \qquad + 2(\curl(\eta_r(\phi_m) \bo_m), \bm{A} \bu_m) - (p_m \eta'(\phi_m) \nabla \phi_m,  \bm{A} \bu_m) \\
& \qquad  - (p_m \eta_r'(\phi_m) \nabla \phi_m,  \bm{A} \bu_m) \\
& \quad =: B_{14} + \cdots + B_{22}.
\end{aligned}
\end{equation}
Meanwhile, choosing $\bm{z} = - \Delta \bo_m$ in \eqref{app:w} results in
\begin{equation}\label{uni:est:10}
\begin{aligned}
&\int_\Omega (c_d(\phi_m) + c_a(\phi_m)) |\Delta \bo_m|^2 + c_0(\phi_m) |\nabla \div \bo_m|^2 \, dx \\
& \quad = (\rho(\phi_m) \pd_t \bo_m, \Delta \bo_m) + (\rho(\phi_m) (\bu_m \cdot \nabla) \bo_m, \Delta \bo_m) \\
& \qquad - \tfrac{\overline{\rho}_1 - \overline{\rho}_2}{2} ((\nabla \mu_m \cdot \nabla) \bo_m, \Delta \bo_m)  -  (c_0'(\phi_m) \nabla \phi_m \div \bo_m, \nabla \div \bo_m)  \\
& \qquad- (2 c_d'(\phi_m) \nabla \phi_m \D \bo_m + 2 c_a'(\phi_m) \nabla \phi_m \W \bo_m, \Delta \bo_m)\\
&\qquad  - ( (c_d(\phi_m) - c_a(\phi_m)) \nabla \div \bo_m, \Delta \bo_m) - (2 \eta_r(\phi_m) \curl \bu_m, \Delta \bo_m) \\
& \qquad + (4 \eta_r(\phi_m) \bo_m , \Delta \bo_m) \\
& \quad =: B_{23} + \cdots + B_{30}.
\end{aligned}
\end{equation}
On the other hand, taking the gradient of \eqref{app:phi} and testing with $\nabla \Delta \mu_m$ leads to 
\begin{align}\label{uni:est:7}
\| \nabla \Delta \mu_m \|^2 = (\nabla \pd_t \phi_m, \nabla \Delta \mu_m) + (\nabla (\bu_m \cdot \nabla \phi_m), \nabla \Delta \mu_m).
\end{align}
By elliptic regularity theory it holds that 
\[
\| \mu_m \|_{H^3} \leq C (\| \mu_m \|_{H^1} + \| \nabla \Delta \mu_m \|),
\]
and hence from \eqref{mean:est} and the upper bound \eqref{Hm:ub}, we infer that
\begin{equation}\label{uni:est:8}
\begin{aligned}
\| \mu_m \|_{H^3}^2 & \leq C ( 1 + \| \nabla \mu_m \|^2 + \alpha \| \pd_t \phi_m \|^2) \\
& \quad + (\nabla \pd_t \phi_m, \nabla \Delta \mu_m)  + ( \nabla (\bu_m \cdot \nabla \phi_m), \nabla \Delta \mu_m ) \\
& \leq C(1 + C_0 + H_m) + B_{31} + B_{32},
\end{aligned}
\end{equation}
where $C_0$ is the positive constant appearing in \eqref{Hm:lb} that ensures $C_0 + H_m$ is non-negative.
Let $\beta_1$, $\beta_2$ and $\beta_3$ be positive constants yet to be determined.  Multiplying \eqref{uni:est:6} by $\beta_1$, \eqref{uni:est:10} by $\beta_2$ and \eqref{uni:est:8} by $\beta_3$, and adding the results to \eqref{reg:main} leads to the inequality
\begin{equation}\label{reg:main2}
\begin{aligned}
& \frac{d}{dt} H_m(t) + \rho_*(\|\pd_t \bu_m \|^2 + \|\pd_t \bo_m\|^2 ) + \| \nabla \pd_t \phi_m \|^2 + \beta_1 (\eta_* + \eta_{r,*}) \| \bm{A} \bu_m \|^2 \\
& \qquad + \beta_2 \int_\Omega (c_d(\phi_m) + c_a(\phi_m)) |\Delta \bo_m|^2 + c_0(\phi_m) |\nabla \div \bo_m|^2 \, dx + \beta_3 \| \mu_m \|_{H^3}^2 \\
& \quad \leq C(1 + C_0 + H_m) + (B_1 + \cdots + B_{13}) + \beta_1 (B_{14} + \cdots + B_{22}) \\
& \qquad + \beta_2 (B_{23} + \cdots + B_{30}) + \beta_3(B_{31} + B_{32}).
\end{aligned}
\end{equation}
We now estimate the right-hand side, and while many of them are similar to those in \cite{CHNS3Dstrsol}, there are some new contributions coming from terms involving $\W \bu_m$ and $\bo_m$. In the following we often make use of the Poincar\'e inequality applied $\bu_m$, $\bo_m$ and $\pd_t \phi_m$, the Sobolev inequalities \eqref{L3estimate} and \eqref{Linftyestimate}, Korn's inequality \eqref{Korn.ineq.}, previous uniform estimates \eqref{uni:est:1}-\eqref{uni:est:4}, particularly \eqref{Lapphi:nablamu} and \eqref{mean:est} to control $\| \phi_m \|_{H^2}$ and $\| \mu_m \|_{H^1}$, as well as the lower bound \eqref{Hm:lb} for $H_m$:
\begin{align*}
|B_1|& \leq C \|\bu_m \|_{L^6} \| \nabla \bu_m \|_{L^3} \| \pd_t \bu_m \| \leq \frac{\rho_*}{10} \| \pd_t \bu_m \|^2 + C \| \D \bu_m \|^3 \| \bm{A} \bu_m \| \\
& \leq \frac{\rho_*}{10} \| \pd_t \bu_m \|^2 + \frac{\eta_* \beta_1}{44} \| \bm{A} \bu_m \|^2 + C(1 + C_0 + H_m)^3,  \\
|B_2|& \leq C\| \pd_t \phi_m \|_{L^6} \| \D \bu_m \|_{L^3} \| \D \bu_m \| \leq \frac{1}{16} \| \nabla \pd_t \phi_m \|^2 + C \| \bm{A} \bu_m \| \| \D \bu_m \|^3 \\
& \leq \frac{1}{14} \| \nabla \pd_t \phi_m \|^2 + \frac{\eta_* \beta_1}{44} \| \bm{A} \bu_m \|^2 + C(1+ C_0 + H_m)^3, \\
|B_3| & \leq C\| \pd_t \phi_m \|_{L^6} \| \W \bu_m \|_{L^3} \| \W \bu_m \| \leq \frac{1}{16} \| \nabla \pd_t \phi_m \|^2 + C \| \bm{A} \bu_m \| \| \W \bu_m \|^3 \\
& \leq \frac{1}{14} \| \nabla \pd_t \phi_m \|^2 + \frac{\eta_* \beta_1}{44} \| \bm{A} \bu_m \|^2 + C(1+ C_0 + H_m)^3, \\
|B_4|& \leq C \| \nabla \mu_m \|_{L^\infty} \| \nabla \bu_m \| \| \pd_t \bu_m \| \leq C \| \nabla \mu_m \|_{H^1}^{1/2} \| \mu_m \|_{H^3}^{1/2} \| \D \bu_m \| \| \pd_t \bu_m \| \\
& \leq \frac{\rho_*}{10} \| \pd_t \bu_m \|^2 +C \| \nabla \mu_m \|^{1/2} \| \mu_m \|_{H^3}^{3/2} \| \D \bu_m \|^2 \\
&  \leq \frac{\rho_*}{10} \| \pd_t \bu_m \|^2 + \frac{\beta_3}{10} \| \mu_m \|_{H^3}^{2} + C \| \nabla \mu_m \|^2 \| \D \bu_m \|^8  \\
& \leq \frac{\rho_*}{10} \| \pd_t \bu_m \|^2 + \frac{\beta_3}{10} \| \mu_m \|_{H^3}^2 + C(1 +C_0 + H_m)^5.
\end{align*}
Notice that $B_{12}$ is exactly the same as $B_{5}$, and so
\begin{align*}
|B_5 + B_{12}| & \leq 2 \| \mu_m \|_{L^6} \| \nabla \phi_m \|_{L^3} \| \pd_t \bu_m \| \leq \frac{\rho_*}{5} \| \pd_t \bu_m \|^2 + C \| \phi_m \|_{H^2}^2 \| \mu_m \|_{H^1}^2 \\
& \leq \frac{\rho_*}{5} \| \pd_t \bu_m \|^2 + C(1 + C_0 + H_m)^2.
\end{align*}
Similar to the estimation for $B_2$ and $B_3$ we have
\begin{align*}
|B_{6}| & \leq C\| \pd_t \phi_m \|_{L^6} ( \|\div  \bo_m \|_{L^3} \|\div \bo_m \| + \| \D \bo_m \|_{L^3} \| \D \bo_m \| + \| \W \bo_m \|_{L^3} \| \W \bo_m \|) \\
& \leq C \| \pd_t \phi_m \|_{L^6} \| \nabla \bo_m \|_{L^3} \| \nabla \bo_m \| \leq C \| \nabla \pd_t \phi_m \| \| \nabla \bo_m \|^{3/2} \| \Delta \bo_m \|^{1/2} \\
& \leq \frac{\beta_2 c_*}{32} \| \Delta \bo_m \|^2 + \frac{1}{14} \| \nabla \pd_t \phi_m \|^{2} + C \| \D \bo_m \|^6 \\
& \leq \frac{\beta_2 c_*}{32} \| \Delta \bo_m \|^2 + \frac{1}{14} \| \nabla \pd_t \phi_m \|^2 + C (1 + C_0 + H_m)^3,
\end{align*}
while analogous to $B_1$ and $B_4$, respectively:
\begin{align*}
|B_{7}| & \leq \rho^* \| \bu_m \|_{L^6} \| \nabla \bo_m \|_{L^3} \| \pd_t \bo_m \| \leq C \| \nabla \bu_m \| \| \nabla \bo_m \|^{1/2} \| \Delta \bo_m \|^{1/2} \| \pd_t \bo_m \| \\
& \leq \frac{\rho_*}{8} \| \pd_t  \bo_m \|^2 + \frac{\beta_2 c_*}{32} \| \Delta \bo_m \|^2 + C \|\D \bu_m  \|^4 \| \D \bo_m \|^2, \\
& \leq \frac{\rho_*}{8} \| \pd_t  \bo_m \|^2 + \frac{\beta_2 c_*}{32} \| \Delta \bo_m \|^2 + C (1+ C_0  + H_m)^3, \\
|B_{8}| & \leq C \| \nabla \mu_m \|_{L^\infty} \| \nabla \bo_m \| \| \pd_t \bo_m \| \leq C \| \nabla \mu_m \|_{H^1}^{1/2} \| \mu_m \|_{H^3}^{1/2} \| \D \bo_m \| \| \pd_t \bo_m \| \\
& \leq \frac{\rho_*}{8} \| \pd_t \bo_m \|^2 + \frac{\beta_3}{10} \| \mu_m \|_{H^3}^2 + C \| \nabla \mu_m \|^2 \| \D \bo_m \|^8, \\
& \leq \frac{\rho_*}{8} \| \pd_t \bo_m \|^2 + \frac{\beta_3}{10} \| \mu_m \|_{H^3}^2 + C(1+ C_0 +H_m)^5.
\end{align*}
Meanwhile, 
\begin{align*}
|B_9 + B_{10}| & \leq C \| \pd_t \phi_m \|_{L^6} \| \| \bo_m \|_{L^3} ( \| \nabla \bu_m \| + \| \bo_m \|)\\
& \leq C \| \nabla \pd_t \phi_m \| \| \D \bo_m \| ( \| \D \bu_m \| + \| \D \bo_m \|) \\
& \leq \frac{1}{14} \| \nabla \pd_t \phi_m \|^2 + C (1 + C_0 + H_m)^2, \\
|B_{11}| &  \leq C \| \pd_t \phi_m \|_{(H^1)^*} \| \nabla \pd_t \phi_m \| \leq \frac{1}{14} \| \nabla \pd_t \phi_m \|^2 + C(1 + C_0 + H_m), \\
|B_{13}| & \leq \| \mu_m \|_{L^6} \| \bu_m \|_{L^3} \| \nabla \pd_t \phi_m \| \leq \frac{1}{10} \| \nabla \pd_t \phi_m \|^2 + C \| \D \bu_m \|^2 \| \mu_m \|_{H^1}^2 \\
& \leq  \frac{1}{14} \| \nabla \pd_t \phi_m \|^2 + C(1+ C_0 + H_m)^2,\\
|B_{14}| & \leq \rho^* \| \pd_t \bu_m \| \| \bm{A} \bu_m \| \leq \frac{\rho_*}{10 \beta_1} \| \pd_t \bu_m \|^2 + \frac{4 (\rho^*)^2 \beta_1}{\rho_*} \| \bm{A} \bu_m \|^2.
\end{align*}
Similar to $B_{1}$, $B_{4}$ and $B_{5}$, respectively, we have
\begin{align*}
|B_{15}| & \leq \rho^* \| \bu_m \|_{L^6} \| \nabla \bu_m \|_{L^3} \| \bm{A} \bu_m \| \leq C \| \D \bu_m \|^{3/2} \| \bm{A} \bu_m \|^{3/2} \\
&  \leq \frac{\eta_*}{44} \| \bm{A} \bu_m \|^2 + C(1 + C_0 + H_m)^3, \\
|B_{16}| & \leq C \| \nabla \mu_m \|_{L^\infty} \| \nabla \bu_m \| \| \bm{A} \bu_m\| \leq C \| \nabla \mu_m \|_{H^1}^{1/2} \| \mu_m \|_{H^3}^{1/2} \| \D \bu_m \| \| \bm{A} \bu_m \| \\
& \leq \frac{\eta_*}{44}\| \bm{A} \bu_m \|^2 + C \| \nabla \mu_m \|^{1/2} \| \mu_m \|_{H^3}^{3/2} \| \D \bu_m \|^2 \\
& \leq \frac{\eta_*}{44} \| \bm{A} \bu_m \|^2 + \frac{\beta_3}{10 \beta_1} \|\mu_m \|_{H^3}^2 + C(1 + C_0 + H_m)^5, \\
|B_{17}| & \leq \| \mu_m \|_{L^6} \| \nabla \phi_m \|_{L^3} \| \bm{A} \bu_m \| \\
&  \leq \frac{\eta_*}{44} \| \bm{A} \bu_m \|^2 + C \| \mu_m \|_{H^1}^2 \| \phi_m \|_{H^2}^2 \leq \frac{\eta_*}{44} \| \bm{A} \bu_m \|^2 + C(1 + C_0 + H_m)^2,
\end{align*}
while recalling the estimations for $B_{2}$ and $B_{3}$,
\begin{align*}
|B_{18}| & \leq C \| \D \bu_m \|_{L^3} \| \nabla \phi_m \|_{L^6} \| \bm{A} \bu_m \| \leq C \| \D \bu_m \|^{1/2} \| \bm{A} \bu_m \|^{3/2} \| \phi_m \|_{H^2} \\
& \leq \frac{\eta_*}{44} \| \bm{A} \bu_m \|^2 + C \| \D \bu_m \|^2 \| \| \phi_m \|_{H^2}^4 \leq \frac{\eta_*}{44} \| \bm{A} \bu_m \|^2 + C (1 + C_0 + H_m)^3, \\
|B_{19}| & \leq C\| \nabla \bu_m \|_{L^3} \| \nabla \phi_m \|_{L^6} \| \bm{A} \bu_m \| \leq C\| \D \bu_m \|^{1/2} \| \bm{A} \bu_m \|^{3/2} \| \phi_m \|_{H^2} \\
& \leq \frac{\eta_*}{44} \| \bm{A} \bu_m \|^2 + C(1 + C_0 + H_m)^3.
\end{align*}
For $B_{20}$ we employ the product rule for the curl operator:
\begin{align*}
|B_{20}| & \leq C \| \bm{A} \bu_m \| (\| \curl \bo_m \| + \| \nabla \phi_m \|_{L^6} \|  \bo_m\|_{L^3}) \\
& \leq \frac{\eta_*}{44} \| \bm{A} \bu_m \|^2 + C (1 + \| \phi_m \|_{H^2}^2) \| \nabla \bo_m \|^2 \\
& \leq \frac{\eta_*}{44} \| \bm{A} \bu_m \|^2 + C( 1 + C_0 + H_m)^2.
\end{align*}
Invoking the estimate \eqref{pi:est} we find that 
\begin{align*}
|B_{21} +B_{22}| & \leq C \| p_m \|_{L^3} \| \nabla \phi_m \|_{L^6} \| \bm{A} \bu_m \| \leq C \| p_m \|^{1/2} \| p_m \|_{H^1}^{1/2} \| \phi_m \|_{H^2} \| \bm{A} \bu_m \| \\
& \leq C \| \D \bu_m \|^{1/4} \| \bm{A} \bu_m \|^{7/4} \| \phi_m \|_{H^2} \leq \frac{\eta_*}{44} \| \bm{A} \bu_m \|^2 + C (1 + C_0 + H_m)^5, 
\end{align*}
while recalling $B_{14}$, $B_{15}$ and $B_{16}$, respectively:
\begin{align*}
|B_{23}| & \leq \frac{c_*}{32} \| \Delta \bo_m \|^2 + \frac{8(\rho^*)^2}{c_*} \| \pd_t \bo_m \|^2, \\
|B_{24}| & \leq \rho^* \| \bu_m \|_{L^6} \| \nabla \bo_m \|_{L^3} \| \Delta \bo_m \| \leq C \| \nabla \bu_m \| \| \nabla \bo_m \|^{1/2} \| \Delta \bo_m \|^{3/2} \\
& \leq \frac{c_*}{32} \| \Delta \bo_m \|^2 + C\| \D \bu_m \|^4 \| \D \bo_m \|^2 \\
& \leq   \frac{c_*}{32} \| \Delta \bo_m \|^2 + C(1+ C_0 + H_m)^3, \\
|B_{25}| & \leq C \| \nabla \mu \|_{L^\infty} \| \nabla \bo_m\| \| \Delta \bo_m \| \leq C \| \nabla \mu_m \|_{H^1}^{1/2} \| \mu_m \|_{H^3}^{1/2} \| \D \bo_m \| \| \Delta \bo_m \| \\
& \leq \frac{c_*}{32} \| \Delta \bo_m \|^2 + \frac{\beta_3}{10\beta_2} \| \mu_m \|_{H^3}^2 + C \| \nabla \mu_m \|^2 \| \D \bo_m \|^8, \\
& \leq \frac{c_*}{32} \| \Delta \bo_m \|^2 + \frac{\beta_3}{10\beta_2} \| \mu_m \|_{H^3}^2 + C (1+ C_0 + H_m)^5.
\end{align*}
Similar to $B_{18}$ and $B_{19}$ we have
\begin{align*}
|B_{26}| & \leq C \| \nabla \phi_m \|_{L^6} \| \div \bo_m \|_{L^3} \| \nabla \div \bo_m \| \leq C \| \phi_m \|_{H^2} \| \div \bo_m \|^{1/2} \| \nabla \div \bo_m \|^{3/2} \\
& \leq \frac{c_*}{4} \| \nabla \div \bo_m \|^2 + C \| \phi_m \|_{H^2}^4 \| \div \bo_m \|^2, \\
& \leq \frac{c_*}{4} \| \nabla \div \bo_m \|^2 + C(1+ C_0 + H_m)^3, \\
|B_{27}| & \leq C \| \nabla \phi_m \|_{L^6} \| \nabla \bo_m \|_{L^3} \| \Delta \bo_m \| \leq C \| \phi_m \|{H^2} \| \D \bo_m \|^{1/2} \| \Delta \bo_m \|^{3/2} \\
& \leq \frac{c_*}{32} \| \Delta \bo_m \|^2 + C \| \phi_m \|_{H^2}^4 \| \D \bo_m \|^2, \\
& \leq \frac{c_*}{32} \| \Delta \bo_m \|^2 + C(1+ C_0 + H_m)^3,
\end{align*}
while by Young's inequality 
\begin{align*}
|B_{29} + B_{30}| & \leq 4 \eta_r^* \| \Delta \bo_m \| \| \tfrac{1}{2} \curl \bu_m - \bo_m \| \leq \frac{c_*}{16} \| \Delta \bo_m \|^2 + C \| \tfrac{1}{2} \curl \bu_m - \bo_m \|^2 \\
& \leq \frac{c_*}{16} \| \Delta \bo_m \|^2 + C(1 + C_0 + H_m).
\end{align*}
Concerning $B_{28}$ we make use of the technical assumptions $|c_d - c_a| \leq (c_d + c_a)$ and $2 c_0 > |c_d - c_a|$ to deduce that
\begin{align*}
|B_{28}| & \leq \frac{1}{2} \int_\Omega |c_d (\phi_m) - c_a(\phi_m)| \big ( |\nabla \div \bo_m|^2 + |\Delta \bo_m|^2 \big ) \, dx \\
& \leq \frac{1}{2} \int_\Omega (c_d(\phi_m) + c_a(\phi_m)) |\Delta \bo_m|^2 + c_0(\phi_m) |\nabla \div \bo_m|^2 \, dx.
\end{align*}
Lastly,
\begin{align*}
|B_{31} | & \leq \frac{1}{14 \beta_3} \| \nabla \pd_t \phi_m \|^2 + \frac{7 \beta_3}{2} \| \mu_m \|_{H^3}^2 , \\
|B_{32}| & \leq  C (\| \D \bu_m \|_{L^3} \| \nabla \phi_m \|_{L^6} + \| \phi_m \|_{H^2} \| \bu_m \|_{L^\infty}) \| \mu_m \|_{H^3}  \\
& \leq C \| \D \bu_m \|^{1/2} \| \bm{A} \bu_m \|^{1/2} \| \phi_m \|_{H^2} \| \mu_m \|_{H^3} \\
& \leq \frac{\eta_* \beta_1}{44 \beta_3} \| \bm{A} \bu_m \|^2 + \frac{1}{10} \| \mu_m \|_{H^3}^2 + C \| \D \bu_m \|^2 \| \phi_m \|_{H^2}^4 \\
& \leq \frac{\eta_* \beta_1}{44 \beta_3} \| \bm{A} \bu_m \|^2 + \frac{1}{10} \| \mu_m \|_{H^3}^2 + C(1 + C_0 + H_m)^3.
\end{align*}
Substituting all the above estimates into \eqref{reg:main2} leads to the following differential inequality
\begin{align*}
& \frac{d}{dt} H_m(t) + \frac{\rho_*}{2} \|\pd_t \bu_m \|^2 + \Big ( \frac{3\rho_*}{4} -\beta_2 \frac{8(\rho^*)^2}{c_*}  \Big ) \|\pd_t \bo_m\|^2 + \frac{\beta_3}{2} (1 - 7 \beta_3) \| \mu_m \|_{H^3}^2\\
& \qquad + \frac{1}{2} \| \nabla \pd_t \phi_m \|^2 + \beta_1 \eta_* \Big ( \frac{3}{4} - \beta_1 \frac{4 (\rho^*)^2 }{\rho_*} \Big ) \| \bm{A} \bu_m \|^2 + \frac{\beta_2 c_*}{4} (\|\Delta \bo_m\|^2 + \|\nabla \div \bo_m\|^2)  \\
& \quad \leq C(1 + C_0 + H_m)^5.
\end{align*}
Choosing $\beta_1$, $\beta_2$ and $\beta_3$ sufficiently small, we infer that
\begin{equation}\label{reg:main3}
\begin{aligned}
& \frac{d}{dt} H_m(t) + \frac{\rho_*}{2} (\|\pd_t \bu_m \|^2 +  \|\pd_t \bo_m\|^2) + \frac{\beta_3}{4}\| \mu_m \|_{H^3}^2\\
& \qquad + \frac{1}{2} \| \nabla \pd_t \phi_m \|^2 + \frac{\beta_1 \eta_*}{2} \| \bm{A} \bu_m \|^2 + \frac{\beta_2 c_*}{4} (\|\Delta \bo_m\|^2 + \|\nabla \div \bo_m\|^2)  \\
& \quad \leq C(1 + C_0 + H_m)^5.
\end{aligned}
\end{equation}
where the positive constant $C$ is independent of the approximation parameters $m$, $\alpha$ and $k$. From \cite[(3.63)]{CHNS3Dstrsol} we have
\[
\alpha \| \pd_t \phi_m(0) \|^2 + \| \nabla \mu_m(0) \|^2 \leq C \| \nabla (-\Delta \phi_{0,k} + F'(\phi_{0,k})) \|^2 + C \| \bu_m(0) \|^2,
\]
so that upon recalling the properties \eqref{prop.approx.IC.phi0} and the upper bound \eqref{Hm:ub} for $H_m$, we see
\begin{align*}
H_m(0) \leq C \Big ( 1 + \| \bu_0 \|_{H^1}^2 + \| \bo_0 \|_{H^1}^2 + \| - \Delta \phi_0 + F'(\phi_0) \|_{H^1}^2 + \| \phi_0 \|_{H^1}^2 \Big ) =: L_0.
\end{align*}
Then, after neglecting the non-negative terms on the left-hand side in \eqref{reg:main3} we deduce from the resulting Bernoulli-type differential inequality 
\[
\frac{1}{(1 + C_0 + H_m)^5} \frac{d}{dt} (1 + C_0 + H_m) \leq C
\]
so that for all $t \in [0, \widetilde{T}_0)$ with $\widetilde{T}_0 := \frac{1}{4C(1 + C_0 + L_0)^4}$ independent of $m$, $\alpha$ and $k$,
\[
1 + C_0 + H_m(t) \leq \frac{1 + C_0 + L_0}{\Big ( 1 - 4Ct \big ( 1 + C_0 + L_0 \big )^4 \Big)^{1/4}}.
\]
Hence, for fixed $T_0 \in (0, \widetilde{T}_0)$, we infer from \eqref{Hm:lb} that 
\begin{align}\label{uni:est:13}
\sup_{t \in [0,T_0]} \Big ( \| \nabla \bu_m(t) \|^2 + \| \nabla \bo_m(t) \|^2 + \| \nabla \mu_m(t) \|^2 + \alpha \| \pd_t \phi_m(t) \|^2 \Big ) \leq C.
\end{align}
Then, invoking \eqref{Lapphi:nablamu} and \eqref{uni:est:3} we also have
\begin{align}\label{uni:est:14}
\sup_{t \in [0,T_0]} \Big ( \| \phi_m(t) \|_{H^2}^2 + \| \mu_m(t) \|_{H^1}^2 + \| F'(\phi_m(t)) \|^2 \Big ) \leq C.
\end{align}
Integrating \eqref{reg:main3} on $[0,T_0]$ provides
\begin{equation}\label{uni:est:15}
\begin{aligned}
& \int_0^{T_0} \Big (\| \pd_t \bu_m \|^2 + \| \pd_t \bo_m \|^2 + \| \nabla \pd_t \phi_m \|^2 + \| \bm{A} \bu_m \|^2 \Big ) \, dt \\
& \quad +  \int_0^{T_0} \Big ( \| \Delta \bo_m \|^2 + \| \nabla \div \bo_m \|^2 + \| \mu_m \|_{H^3}^2 \Big ) \, dt \leq C,
\end{aligned}
\end{equation}
and by revisiting Proposition \ref{Pre.result} applied to the convective Cahn--Hilliard subsystem \eqref{app:phi}-\eqref{app:mu}, we obtain the strict separation property where there exists $\tilde{\delta} = \tilde{\delta}(\alpha, k)$ independent of $m$ such that 
\begin{align}\label{uni:est:16}
|\phi_m(x,t)| \leq 1 - \tilde{\delta} \quad \text{ a.e.~in } \Omega \times (0,T_0).
\end{align}

\subsection{Passing to the limit}
The sequence in passing to the limit proceeds as $m \to \infty$, followed by $\alpha \to 0$, then followed by $k \to \infty$. Thanks to the uniform estimates \eqref{uni:est:1}-\eqref{uni:est:4} and higher order uniform estimates \eqref{uni:est:13}-\eqref{uni:est:16}, we deduce that along a non-relabelled subsequence $m \to \infty$, 
\begin{equation}\label{compactness}
\begin{aligned}
\bu_m \to \bu_{\alpha,k} & \text{ weakly* in } L^\infty(0,T_0;\HH^1_{\div}) \cap L^2(0,T_0;\HH^2_{\div}) \cap H^1(0,T_0;\HH_{\div}), \\
\bu_m \to \bu_{\alpha,k} & \text{ strongly in } L^2(0,T_0;\HH^1_{\div}), \\
\bo_m \to \bo_{\alpha,k} & \text{ weakly* in } L^\infty(0,T_0;\HH^1) \cap L^2(0,T_0;\HH^2) \cap H^1(0,T_0;\HH), \\
\bo_m \to \bo_{\alpha,k} & \text{ strongly in } L^2(0,T_0;\HH^1), \\
\phi_m \to \phi_{\alpha,k} & \text{ weakly* in } L^\infty(0,T_0;H^2(\Omega)) \cap H^1(0,T_0;H^1(\Omega)), \\
\phi_m \to \phi_{\alpha,k} & \text{ strongly in } C^0([0,T_0];W^{1,p}(\Omega)) \quad \forall p \in [2,6), \\
\mu_m \to \mu_{\alpha,k} & \text{ weakly* in } L^\infty(0,T_0;H^1(\Omega)) \cap L^2(0,T_0;H^3(\Omega)).
\end{aligned}
\end{equation}
Additionally we have that the limit function $\phi_{\alpha,k}$ satisfies the strict separation property:
\[
\phi_{\alpha,k} \in L^\infty(\Omega \times (0,T_0)) \quad \text{ with } |\phi_{\alpha,k}(x,t)| \leq 1 - \delta \text{ a.e.~in } \Omega \times (0,T_0),
\]
for some $\delta = \delta(\alpha,k)$, while for $f \in \{\rho, \eta, \eta_r, c_0, c_d, c_a\}$ we have
\begin{align}\label{lim:coeff}
f(\phi_m) \to f(\phi_{\alpha,k}) \text{ strongly in } C^0([0,T_0];W^{1,p}(\Omega)) \quad \forall p \in [2,6).
\end{align}
The above weak and strong convergences are sufficient to pass to the limit in the Galerkin approximation \eqref{approx:problem} as $m \to \infty$ to deduce that the limit quadruple $(\bu_{\alpha,k}, \bo_{\alpha,k}, \phi_{\alpha,k}, \mu_{\alpha,k})$ satisfies for a.e.~$t \in [0,T_0)$,
\begin{subequations}\label{lim:m:problem}
\begin{alignat}{2}
 & (\rho(\phi_{\alpha,k}) \pd_t \bu_{\alpha,k}, \bm{v}) + (\rho(\phi_{\alpha,k}) (\bu_{\alpha,k} \cdot \nabla) \bu_{\alpha,k}, \bm{v}) - (\div (2 \eta(\phi_{\alpha,k}) \D \bu_{\alpha,k}), \bm{v}) \\
\notag & \qquad  -(\div (2 \eta_r(\phi_{\alpha,k}) \W \bu_{\alpha,k}), \bm{v}) - \tfrac{\overline{\rho}_1 - \overline{\rho}_2}{2} ((\nabla \mu_{\alpha,k} \cdot \nabla) \bu_{\alpha,k}, \bm{v}) \\
\notag & \quad = (\mu_{\alpha,k} \nabla \phi_{\alpha,k}, \bm{v}) + 2(\curl(\eta_r(\phi_{\alpha,k}) \bo_{\alpha,k}), \bm{v}) \quad \forall \bm{v} \in \HH_{\div}, \\[1ex]
 & (\rho(\phi_{\alpha,k}) \pd_t \bo_{\alpha,k}, \bm{z}) + (\rho(\phi_{\alpha,k}) (\bu_{\alpha,k} \cdot \nabla) \bo_{\alpha,k}, \bm{z}) - (\div (c_0(\phi_{\alpha,k}) \div \bo_{\alpha,k} \mathbb{I}),  \bm{z}) \\
\notag & \qquad - (\div (2c_d(\phi_{\alpha,k}) \D \bo_{\alpha,k}) + 2 c_a(\phi_{\alpha,k}) \W \bo_{\alpha,k}), \bm{z}) - \tfrac{\overline{\rho}_1 - \overline{\rho}_2}{2} ((\nabla \mu_{\alpha,k} \cdot \nabla) \bo_{\alpha,k}, \bm{z}) \\
\notag & \quad = (2\eta_r(\phi_m) \curl \bu_{\alpha,k}, \bm{z}) - (4 \eta_r(\phi_{\alpha,k}) \bo_{\alpha,k}, \bm{z}) \quad \forall \bm{z} \in \HH, \\[1ex]
 & \pd_t \phi_{\alpha,k} + \bu_{\alpha,k} \cdot \nabla \phi_{\alpha,k} = \Delta \mu_{\alpha,k} \quad \text{ a.e.~in } \Omega, \\
& \mu_{\alpha,k} = \alpha \pd_t \phi_{\alpha,k} - \Delta \phi_{\alpha,k} + F'(\phi_{\alpha,k}) \quad \text{ a.e.~in } \Omega,
\end{alignat}
\end{subequations}
together with the following initial-boundary conditions:
\begin{align*}
\begin{cases}
\bu_{\alpha,k} = \bm{0}, \quad \bo_{\alpha,k} = \bm{0}, \quad \pdnu \phi_{\alpha,k} = 0, \quad \pdnu \mu_{\alpha,k} = 0 & \quad \text{ on } \Sigma, \\
\bu_{\alpha,k}(0) = \bu_0, \quad \bo_{\alpha,k}(0) = \bo_0, \quad \phi_{\alpha,k}(0) = \phi_{0,k} & \quad \text{ in } \Omega,
\end{cases}
\end{align*}
where $\phi_{0,k}$ is defined in \eqref{approx.phi0.problem}. 

Next we consider the vanishing viscosity limit $\alpha \to 0$.  By the weak lower semicontinuity of the norms, we obtain from \eqref{uni:est:1}-\eqref{uni:est:4}, \eqref{uni:est:13}-\eqref{uni:est:16} that along a non-relabelled subsequence $\alpha \to 0$, there exists a limit quadruple $(\bu_k, \bo_k, \phi_k, \mu_k)$ such that analogous compactness assertions to \eqref{compactness} holds, along with the strong convergence \eqref{lim:coeff} for the density function $\rho$ and the viscosity functions $\eta$, $\eta_r$, $c_0$, $c_d$ and $c_a$. From the a.e.~convergence of $\phi_{\alpha,k}$ to $\phi_k$ in $\Omega \times (0,T_0)$ arising as a consequence of the strong convergence in $C^0([0,T_0];W^{1,p}(\Omega))$ and the fact that $|\phi_{\alpha,k}(x,t)| < 1$ a.e.~in $\Omega \times (0,T_0)$, we first deduce that $\phi_{k}$ also satisfies the property $|\phi_k(x,t)| < 1$ a.e.~in $\Omega \times (0,T_0)$. Then, continuity of $F'$ provides that $F'(\phi_{\alpha,k}) \to F'(\phi_k)$ a.e.~in $\Omega \times (0,T_0)$, and with the help of Fatou's lemma and the analogue of \eqref{uni:est:15} for $F'(\phi_{\alpha,k})$ we deduce that 
\[
F'(\phi_{\alpha,k}) \to F'(\phi_k) \text{ weakly* in } L^\infty(0,T_0;L^2(\Omega))
\]
and $F'(\phi_k) \in L^2(0,T;L^2(\Omega))$. Thus, passing to the limit $\alpha \to 0$ in \eqref{lim:m:problem}, we find that the limit quadruple $(\bu_k, \bo_k, \phi_k, \mu_k)$ satisfies
\begin{subequations}\label{lim:alpha:problem}
\begin{alignat}{2}
& (\rho(\phi_{k}) \pd_t \bu_{k}, \bm{v}) + (\rho(\phi_{k}) (\bu_{k} \cdot \nabla) \bu_{k}, \bm{v}) - (\div (2 \eta(\phi_{k}) \D \bu_{k}), \bm{v}) \\
\notag & \qquad  -(\div (2 \eta_r(\phi_{k}) \W \bu_{k}), \bm{v}) - \tfrac{\overline{\rho}_1 - \overline{\rho}_2}{2} ((\nabla \mu_{k} \cdot \nabla) \bu_{k}, \bm{v}) \\
\notag & \quad = (\mu_{k} \nabla \phi_{k}, \bm{v}) + 2(\curl(\eta_r(\phi_{k}) \bo_{k}), \bm{v}) \quad \forall \bm{v} \in \HH_{\div}, \\[1ex]
 & (\rho(\phi_{k}) \pd_t \bo_{k}, \bm{z}) + (\rho(\phi_{k}) (\bu_{k} \cdot \nabla) \bo_{k}, \bm{z}) - (\div (c_0(\phi_{k}) \div \bo_{k} \mathbb{I}),  \bm{z}) \\
\notag & \qquad - (\div (2c_d(\phi_{k}) \D \bo_{k} + 2 c_a(\phi_{k}) \W \bo_{k}), \bm{z}) - \tfrac{\overline{\rho}_1 - \overline{\rho}_2}{2} ((\nabla \mu_{k} \cdot \nabla) \bo_{k}, \bm{z}) \\
\notag & \quad = (2\eta_r(\phi_{k}) \curl \bu_{k}, \bm{z}) - (4 \eta_r(\phi_{k}) \bo_{k}, \bm{z}) \quad \forall \bm{z} \in \HH, \\[1ex]
 & \pd_t \phi_{k} + \bu_{k} \cdot \nabla \phi_{k} = \Delta \mu_{k} \quad \text{ a.e.~in } \Omega, \\
& \mu_{k} = - \Delta \phi_{k} + F'(\phi_{k}) \quad \text{ a.e.~in } \Omega,
\end{alignat}
\end{subequations}
together with the following initial-boundary conditions:
\begin{align*}
\begin{cases}
\bu_{k} = \bm{0}, \quad \bo_{k} = \bm{0}, \quad \pdnu \phi_{k} = 0, \quad \pdnu \mu_{k} = 0 & \quad \text{ on } \Sigma, \\
\bu_{k}(0) = \bu_0, \quad \bo_{k}(0) = \bo_0, \quad \phi_{k}(0) = \phi_{0,k} & \quad \text{ in } \Omega.
\end{cases}
\end{align*}
The passage to the limit $k \to \infty$ follows analogously as the above arguments for the limit $\alpha \to 0$ and so we omit the details. From this we infer the existence of a limit quadruple $(\bu, \bo, \phi, \mu)$ satisfying analogous compactness assertion as \eqref{compactness} and 
\begin{subequations}\label{lim:k:problem}
\begin{alignat}{2}
 & (\rho(\phi) \pd_t \bu, \bm{v}) + (\rho(\phi) (\bu \cdot \nabla) \bu, \bm{v}) - (\div (2 \eta(\phi) \D \bu), \bm{v}) \\
\notag & \qquad  -(\div (2 \eta_r(\phi) \W \bu), \bm{v}) - \tfrac{\overline{\rho}_1 - \overline{\rho}_2}{2} ((\nabla \mu \cdot \nabla) \bu, \bm{v}) \\
\notag & \quad = (\mu \nabla \phi, \bm{v}) + 2(\curl(\eta_r(\phi) \bo), \bm{v}) \quad \forall \bm{v} \in \HH_{\div}, \\[1ex]
& (\rho(\phi) \pd_t \bo, \bm{z}) + (\rho(\phi) (\bu \cdot \nabla) \bo, \bm{z}) - (\div (c_0(\phi) \div \bo \mathbb{I}),  \bm{z}) \\
\notag & \qquad - (\div (2c_d(\phi) \D \bo) + 2 c_a(\phi) \W \bo), \bm{z})  - \tfrac{\overline{\rho}_1 - \overline{\rho}_2}{2} ((\nabla \mu \cdot \nabla) \bo, \bm{z}) \\
\notag & \quad = (2\eta_r(\phi) \curl \bu, \bm{z}) - (4 \eta_r(\phi) \bo, \bm{z}) \quad \forall \bm{z} \in \HH, \\[1ex]
 & \pd_t \phi + \bu \cdot \nabla \phi = \Delta \mu \quad \text{ a.e.~in } \Omega, \\
& \mu = - \Delta \phi + F'(\phi) \quad \text{ a.e.~in } \Omega,
\end{alignat}
\end{subequations}
together with the following initial-boundary conditions:
\begin{align*}
\begin{cases}
\bu = \bm{0}, \quad \bo = \bm{0}, \quad \pdnu \phi = 0, \quad \pdnu \mu = 0 & \quad \text{ on } \Sigma, \\
\bu(0) = \bu_0, \quad \bo(0) = \bo_0, \quad \phi(0) = \phi_{0} & \quad \text{ in } \Omega.
\end{cases}
\end{align*}
In order to recover the pressure variable, we first express 
\[
\mu \nabla \phi = - \div (\nabla \phi \otimes \nabla \phi) + \nabla \Big ( \frac{1}{2} |\nabla \phi|^2 + F(\phi) \Big ),
\]
and by invoking classical arguments, see e.g.~\cite[Lemma III.1.1]{Galdi} or \cite[Lemma IV.1.4.1]{Sohr}, to deduce the existence of $p \in L^2(0,T_0;L^2(\Omega))$ with $\overline{p}(t) = 0$ such that 
\begin{align*}
\nabla p & = - \rho(\phi) \pd_t \bu - \rho(\phi)( \bu \cdot \nabla) \bu + \div (2 \eta(\phi) \D \bu + 2 \eta_r(\phi) \W \bu) \\
& \quad +\rho'(\phi) (\nabla \mu \cdot \nabla) \bu - \div (\nabla \phi \otimes \nabla \phi) + 2 \curl(\eta_r(\phi) \bo)
\end{align*}
holds in the sense of distributions for a.e.~$t \in (0,T_0)$. But since the right-hand side belongs to $L^2(\Omega)$ for a.e.~$t \in (0,T_0)$, we infer that $p \in L^2(0,T_0;H^1(\Omega))$. Furthermore, employing the regularity theory of the Cahn--Hilliard equation with logarithmic potential, see e.g.~\cite[Lemma 2]{HAbel.EandU.ModelH} or \cite[Theorem A.2]{AGAMRT}, we find that $\phi \in L^\infty(0,T_0;W^{2,6}(\Omega))$ and $F'(\phi) \in L^\infty(0,T_0;L^6(\Omega))$. This completes the proof of the existence of local-in-time strong solutions.

Using the continuous embedding $W^{1,4}(\Omega) \subset C^0(\overline{\Omega})$ and the interpolation of $W^{1,4}(\Omega)$ via $H^1(\Omega)$ and $W^{2,6}(\Omega)$, we can infer in a similar fashion to \cite{CHNS3Dstrsol} that for all $t_1, t_2 \in [0,T_0]$, 
\begin{align*}
\| \phi(t_1) - \phi(t_2) \|_{C^0(\overline{\Omega})} & \leq C \| \phi(t_1) - \phi(t_2) \|_{H^1}^{5/8} \| \phi(t_1) - \phi(t_2) \|_{W^{2,6}}^{3/8} \\
& \leq C\| \pd_t \phi \|_{L^2(0,T_0;H^1(\Omega))}^{5/8} |t_1 - t_2|^{5/16},
\end{align*}
i.e., $\phi \in C^{0,\frac{5}{16}}([0,T_0];C^0(\overline{\Omega}))$. Along with the assumption on the initial condition $|\phi_{0}(x)| \leq 1- \delta_0$ for some $\delta_0 > 0$, we then deduce the existence of a time $T_1 \in (0,T_0]$ such that
\begin{align}\label{strict:bdd}
|\phi(x,t)| \leq 1- \frac{\delta_0}{2} \quad \forall (x,t) \in \overline{\Omega} \times [0,T_1].
\end{align}
This further implies that $F'(\phi) \in L^\infty(0,T_1;H^2_n(\Omega))$, and by comparison of terms yields that $\Delta \phi \in L^2(0,T_1;H^2_n(\Omega))$.

\section{Continuous dependence and uniqueness of strong solutions}\label{sec:unique}
We now consider strong solutions $\{(\bu_i, \bo_i,  p_i, \phi_i, \mu_i)\}_{i=1,2}$ on the time interval $[0,T_0]$ originating from the initial conditions $\{(\bu_{0,i}, \bo_{0,i}, \phi_{0,i})\}_{i=1,2}$. Setting $f := f_1 - f_2$ for $f \in \{ \bu, p, \bo, \phi, \mu, \rho, \eta, \eta_r, c_0, c_d, c_a\}$ and $g_i = g(\phi_i)$ for $g \in \{\rho, \eta, \eta_r, c_0, c_d, c_a\}$, we find that the differences between the two strong solutions satisfy
\begin{subequations}
\begin{alignat}{2}
\label{uniq:u} & \rho_1 \pd_t \bu + \rho \pd_t \bu_2 + \rho_1( \bu_1 \cdot \nabla) \bu + \rho_1 (\bu \cdot \nabla)  \bu_2 + \rho (\bu_2 \cdot \nabla) \bu_2 \\
\notag & \qquad - \tfrac{\overline{\rho}_1 - \overline{\rho}_2}{2} (\nabla \mu_1 \cdot \nabla) \bu - \tfrac{\overline{\rho}_1 - \overline{\rho}_2}{2} (\nabla \mu \cdot \nabla) \bu_2  - \div (2 \eta_1 \D \bu + 2 \eta \D \bu_2) \\
\notag & \qquad - \div (2 \eta_{r,1} \W \bu + 2 \eta_{r} \W \bu_2) + \nabla p \\
\notag & \quad = 2 \curl(\eta_{r,1} \bo + \eta_{r} \bo_2) - \div (\nabla \phi_1 \otimes \nabla \phi) - \div (\nabla \phi \otimes \nabla \phi_2), \\
\label{uniq:w} & \rho_1 \pd_t \bo + \rho \pd_t \bo_2 + \rho_1( \bu_1 \cdot \nabla) \bo + \rho_1 (\bu \cdot \nabla)  \bo_2 + \rho (\bu_2 \cdot \nabla) \bo_2 \\
\notag & \qquad - \div (c_{0,1} (\div \bo) \mathbb{I}) + c_0 (\div \bo_2) \mathbb{I} + 2 c_{d,1} \D \bo + 2 c_{d} \D \bo_2) \\
\notag & \qquad - \div ( 2 c_{a,1} \W \bo + 2 c_{a} \W \bo_2) - \tfrac{\overline{\rho}_1 - \overline{\rho}_2}{2} (\nabla \mu_1 \cdot \nabla) \bo - \tfrac{\overline{\rho}_1 - \overline{\rho}_2}{2} (\nabla \mu \cdot \nabla) \bo_2 \\
\notag & \quad = 2 \eta_{r,1} (\curl \bu - 2 \bo) + 2 \eta_r (\curl \bu_2 - 2 \bo_2), \\
\label{uniq:phi} & \pd_t \phi + \bu_1 \cdot \nabla \phi + \bu \cdot \nabla \phi_2 = \Delta \mu, \\
\label{uniq:mu} & \mu = - \Delta \phi + F'(\phi_1) - F'(\phi_2),
\end{alignat}
\end{subequations}
almost everywhere in $\Omega \times (0,T_0)$, along with the initial-boundary conditions
\begin{align*}
\begin{cases}
\bu = \bm{0}, \quad \bo = \bm{0}, \quad \pdnu \phi = 0, \quad \pdnu \mu = 0 & \quad \text{ on } \pd \Omega \times (0,T_0), \\
\bu(0) = \bu_{0,1} - \bu_{0,2}, \quad \bo(0) = \bo_{0,1} - \bo_{0,2}, \quad \phi(0) = \phi_{0,1} - \phi_{0,2} & \quad \text{ on } \Omega.
\end{cases}
\end{align*}
Note that when we multiply \eqref{model:phi} for $(\phi_1, \bu_1, \mu_1)$ by $\rho'(\phi_1) = \tfrac{\overline{\rho}_1 - \overline{\rho}_2}{2}$ and then test with $\frac{1}{2} |\bu|^2$, we have the following identity
\begin{align}\label{rho:cancel}
- (\pd_t \rho_1, \tfrac{1}{2} |\bu|^2) + (\rho_1 \bu_1, \nabla \tfrac{1}{2}|\bu|^2) - \tfrac{\overline{\rho}_1 - \overline{\rho}_2}{2} (\nabla \mu_1, \nabla \tfrac{1}{2}|\bu|^2) = 0,
\end{align}
so that upon testing \eqref{uniq:u} by $\bu$, we note a simplification involving the first, third and sixth terms, leading to 
\begin{equation}\label{uniq:1}
\begin{aligned}
& \frac{1}{2} \frac{d}{dt} \int_\Omega \rho_1 |\bu|^2 \, dx + \int_\Omega 2\eta_1 |\D \bu|^2 + 2 \eta_{r,1} |\W \bu|^2 \,dx -   \int_\Omega 2\eta_{r,1} \bo \cdot \curl \bu \, dx \\
& \quad = -(\rho \pd_t \bu_2 , \bu) - (\rho_1 (\bu \cdot \nabla) \bu_2, \bu) - (\rho ( \bu_2 \cdot \nabla) \bu_2, \bu) + \tfrac{\overline{\rho}_1 - \overline{\rho}_2}{2} ( (\nabla \mu \cdot \nabla) \bu_2, \bu)  \\
& \qquad - (2\eta \D \bu_2, \D \bu) + (\nabla \phi_1 \otimes \nabla \phi + \nabla \phi \otimes \nabla \phi_2, \nabla \bu) \\
& \qquad - (2\eta_r \W \bu_2, \W \bu) + (2\eta_r \bo_2, \curl \bu) \\
& \quad =: D_1 + \cdots + D_8.
\end{aligned}
\end{equation}
Similarly, testing \eqref{uniq:w} by $\bo$ we obtain
\begin{equation}\label{uniq:2}
\begin{aligned}
& \frac{1}{2} \frac{d}{dt} \int_\Omega \rho_1 |\bo|^2 \,dx + \int_\Omega c_{0,1} |\div \bo|^2 + 2 c_{d,1} |\D \bo|^2 + 2 c_{a,1} |\W \bo|^2 \, dx \\
& \qquad + \int_\Omega 2 \eta_{r,1}(2 |\bo|^2 - \curl \bu \cdot \bo) \, dx \\
& \quad = - (\rho \pd_t \bo_2, \bo) - (\rho_1 (\bu \cdot \nabla) \bo_2, \bo) - (\rho(\bu_2 \cdot \nabla) \bo_2, \bo)  + \tfrac{\overline{\rho}_1 - \overline{\rho}_2}{2} ( (\nabla \mu \cdot \nabla) \bo_2, \bo)  \\
& \qquad - (c_{0} \div \bo_2, \div \bo) - (2 c_{d} \D \bo_2, \D \bo) - (2 c_a \W \bo_2, \W \bo) \\
& \qquad + (2 \eta_r (\curl \bu_2 - 2 \bo_2), \bo) \\
& \quad =: D_9 + \cdots + D_{16}.
\end{aligned}
\end{equation}
Meanwhile, taking gradient of \eqref{uniq:phi} and then testing by $\nabla \Delta \phi$ and $\nabla \phi$, respectively yields
\begin{equation}\label{uniq:3}
\begin{aligned}
& \frac{1}{2} \frac{d}{dt} \| \Delta \phi \|^2 + \| \Delta^2 \phi \|^2 \\
& \quad = (\bu_1 \cdot \nabla \phi, \Delta^2 \phi) + (\bu \cdot \nabla \phi_2, \Delta^2 \phi) + (\Delta (F'(\phi_1) - F'(\phi_2)), \Delta^2 \phi) \\
& \quad =: D_{17} + D_{18} + D_{19},
\end{aligned}
\end{equation} 
and
\begin{equation}\label{uniq:3b}
\begin{aligned}
& \frac{1}{2} \frac{d}{dt} \| \nabla \phi \|^2 + \| \nabla \Delta \phi \|^2 \\
& \quad =  -(\nabla (\bu_1 \cdot \nabla \phi), \nabla \phi) - (\nabla (\bu \cdot \nabla \phi_2), \nabla \phi) + (\nabla (F'(\phi_1) - F'(\phi_2)), \nabla \Delta \phi) \\
& \quad =: D_{20} + D_{21} + D_{22}.
\end{aligned}
\end{equation}
Lastly, integrating \eqref{uniq:phi} over $\Omega$ yields
\begin{align}\label{uniq:mean}
\frac{d}{dt} \overline{\phi}(t) = 0 \quad \implies \quad  \frac{d}{dt} \frac{1}{2} |\overline{\phi}(t)|^2 = 0.
\end{align}
Then, upon summing \eqref{uniq:1}-\eqref{uniq:mean} we obtain 
\begin{equation}\label{uniq:4}
\begin{aligned}
& \frac{1}{2}  \frac{d}{dt}\Big (  \int_\Omega \rho_1 |\bo|^2 + \rho_1 |\bo|^2 + |\Delta \phi|^2 + |\nabla \phi |^2  \,dx + |\overline{\phi}|^2 \Big ) + \int_\Omega |\Delta^2 \phi|^2 + |\nabla \Delta \phi |^2 + 2\eta_1 |\D \bu|^2 \, dx \\
& \qquad + \int_\Omega c_{0,1} |\div \bo|^2 + 2 c_{d,1} |\D \bo|^2 + 2 c_{a,1} |\W \bo|^2 + 4 \eta_{r,1} \big | \tfrac{1}{2} \curl \bu - \bo|^2 \, dx \\
& \quad = D_1 + \cdots + D_{22}.
\end{aligned}
\end{equation}
Arguing similarly as in \cite[Section 6]{CHNS2Dstrsol} we have
\begin{align*}
& |D_1 + D_2 + D_3 + D_5 + D_6| \\
& \quad \leq \frac{\eta_*}{2} \| \D \bu \|^2 + C(1 + \| \bu_2 \|_{H^2}^2 + \| \pd_t \bu_2 \|^2) ( \| \bu \|^2 + \| \Delta \phi \|^2),\\
& |D_9 + D_{10} + D_{11} + D_{13} + D_{14} + D_{15}| \\
& \quad \leq \frac{c_*}{2} \| \div \bo \|^2 + c_* (\| \D \bo \|^2 + \| \W \bo \|^2) + \frac{\eta_*}{2} \| \D \bu \|^2 \\
& \qquad + C(1 + \| \bu \|_{H^2}^2 + \| \bo_2 \|_{H^2}^2 + \| \pd_t \bo_2 \|^2) (\| \bu \|^2 + \| \bo \|^2 + \| \Delta \phi \|^2)
\end{align*}
while by \eqref{strict:bdd}, the regularity $\phi_i \in L^\infty(0,T_0;W^{2,6}(\Omega))$ and Sobolev embeddings it holds that 
\begin{align*}
|D_4| & \leq | ((\nabla \Delta \phi \cdot \nabla) \bu_2, \bu)| + |(\nabla (F'(\phi_1) - F'(\phi_2)) \cdot \nabla) \bu_2, \bu)| \\
& \leq \| \nabla \Delta \phi \|_{L^6} \| \nabla \bu_2 \|_{L^3} \| \bu \| + \| F''(\phi_1) \|_{L^\infty} \| \nabla \phi \|_{L^6} \| \nabla \bu_2 \|_{L^3} \| \bu \| \\
& \quad + ( \| F'''(\phi_1) \|_{L^\infty} + \| F'''(\phi_2) \|_{L^\infty}) \| \phi \|_{L^\infty} \| \nabla \phi_2 \|_{L^\infty} \| \nabla \bu_2 \| \| \bu \| \\
& \leq \frac{1}{8} \| \Delta^2 \phi \|^2 + C \| \nabla \bu_2 \|_{L^3}^2 \| \bu \|^2 + C( 1 + \| \nabla \bu_2 \|_{L^3}) (\| \bu \|^2 + \| \Delta \phi \|^2 + \| \phi \|_{H^1}^2), \\
|D_{12}| & \leq \frac{1}{8} \| \Delta^2 \phi \|^2 + C \| \nabla \bo_2 \|_{L^3}^2 \| \bo \|^2 + C( 1 + \| \nabla \bo_2 \|_{L^3}) (\| \bo \|^2 + \| \Delta \phi \|^2  + \| \phi \|_{H^1}^2), \\
|D_{17}+D_{18}| & \leq \| \bu_1 \|_{L^3} \| \nabla \phi \|_{L^6} \| \Delta^2 \phi \| + \| \bu \| \| \nabla \phi_2 \|_{L^\infty} \| \Delta^2 \phi \| \\
& \leq \frac{1}{8} \| \Delta^2 \phi \|^2 + C( \| \bu \|^2 + \| \Delta \phi \|^2), \\
|D_{19}| & \leq |(F''(\phi_1) \Delta \phi + (F''(\phi_1) - F''(\phi_2)) \Delta \phi_2, \Delta^2 \phi)| \\
& \quad + |( F'''(\phi_1)(|\nabla \phi_1|^2 - |\nabla \phi_2|^2) + (F'''(\phi_1) - F'''(\phi_2)) |\nabla \phi_2|^2, \Delta^2 \phi)| \\
&  \leq C \| \Delta \phi \| \| \Delta^2 \phi \| + C( \| F'''(\phi_1) \|_{L^\infty} + \| F'''(\phi_2) \|_{L^\infty}) \| \phi \|_{L^\infty} \| \Delta \phi_2 \| \| \Delta^2 \phi \| \\
& \quad + C( \| \nabla \phi_1 \|_{L^\infty} + \| \nabla \phi_2 \|_{L^\infty}) \| \nabla \phi \| \| \Delta^2 \phi \| \\
& \quad + C(\| F''''(\phi_1) \|_{L^\infty} + \| F''''(\phi_2) \|_{L^\infty}) \| \phi \|_{L^\infty} \| \nabla \phi_2 \|_{L^\infty}^2 \| \Delta^2 \phi \| \\
& \leq \frac{1}{8} \| \Delta^2 \phi \|^2 + C \| \Delta \phi \|^2 + C \| \phi \|_{H^1}^2, \\
|D_{20}| & \leq (\| \nabla \bu_1 \| \| \nabla \phi \|_{L^\infty} + \| \bu_1 \|_{L^6} \|  \phi \|_{W^{2,3}}) \| \nabla \phi \| \leq C(\| \phi \|_{H^1} + \| \nabla \Delta \phi \|) \| \nabla \phi \| \\
& \leq \frac{1}{4} \| \nabla \Delta \phi \|^2 + C\| \phi \|_{H^1}^2, \\
|D_{21}| & \leq (\| \nabla \bu \| \| \nabla \phi_2 \|_{L^\infty} + \| \bu \|_{L^6} \| \phi_2 \|_{W^{2,3}}) \| \nabla \phi \| \\
& \leq C \| \nabla \phi \|^2 + \frac{\eta_*}{2} \| \D \bu \|^2, \\
|D_{22}| & \leq \| F''(\phi_1) \nabla \phi \| \| \nabla \Delta \phi \| + \| (F''(\phi_1) - F''(\phi_2)) \nabla \phi_2 \| \| \nabla \Delta \phi \| \\
& \leq \big [ \| F''(\phi_1) \|_{L^6} \| \nabla \phi \|_{L^3} + (\| F'''(\phi_1) \|_{L^6} + \| F'''(\phi_2) \|_{L^6}) \| \phi \|_{L^3} \| \nabla \phi_2 \|_{L^\infty} \big ] \| \nabla \Delta \phi \| \\
& \leq C \| \phi \|_{H^2} \| \nabla \Delta \phi \| \leq C \| \phi \|_{H^1}^{1/2} ( \| \phi \|_{H^1} + \| \nabla \Delta \phi \|)^{1/2} \| \nabla \Delta \phi \| \\
& \leq \frac{1}{4} \| \nabla \Delta \phi \|^2 + C \| \phi \|_{H^1}^2.
\end{align*}
On the other hand, using the identity \eqref{W:W:curl} with $\bm{a} = \bu$ and $\bm{b} = \bu_2$, we deduce
\begin{align*}
D_7 + D_8 + D_{16} & = \int_\Omega 2 \eta_r (\curl \bu \cdot \bo_2 - \W \bu : \W \bu_2 + \bo \cdot \curl \bu_2 - 2 \bo \cdot \bo_2) \, dx \\
& =  \int_\Omega 4 \eta_r ( \tfrac{1}{2} \curl \bu -  \bo) \cdot ( \bo_2 - \tfrac{1}{2} \curl \bu_2) \, dx \\
& \leq 2 \eta_{r,*} \| \tfrac{1}{2} \curl \bu - \bo \|^2 + C \| \phi \|_{L^\infty}^2 \| \tfrac{1}{2} \curl \bu_2 - \bo_2 \|^2 \\
& \leq 2 \eta_{r,*} \| \tfrac{1}{2} \curl \bu - \bo \|^2 + C \| \Delta \phi \|^2 + C \| \phi \|_{H^1}^2.
\end{align*}
Hence, we deduce from \eqref{uniq:4} the following differential inequality:
\begin{align*}
&\frac{1}{2}  \frac{d}{dt} \Big ( \int_\Omega \rho_1 |\bu|^2 + \rho_1 |\bo|^2 + |\Delta \phi|^2 + |\nabla \phi|^2 \,dx  + |\overline{\phi}|^2 \Big )  +\int_\Omega \frac{1}{2} |\Delta^2 \phi|^2 + \frac{\eta_*}{2} |\D \bu|^2 \, dx \\
& \qquad + \int_\Omega \frac{c_{0,*}}{2} |\div \bo|^2 + 2 c_{d,*} |\D \bo|^2 + 2 c_{a,*} |\W \bo|^2 + 2 \eta_{r,*} \big | \tfrac{1}{2} \curl \bu - \bo|^2 \, dx \\
& \quad \leq C(1 + \| \bu_2 \|_{H^2}^2 + \| \bo_2 \|_{H^2}^2 + \| \pd_t \bu_2 \|^2 + \| \pd_t \bo \|^2) \\
& \qquad \times ( \| \bu \|^2 + \| \bo \|^2 + \| \Delta \phi \|^2 + \| \nabla \phi \|^2 + |\overline{\phi}|^2),
\end{align*}
and application of the Gr\"onwall inequality yields 
\begin{align*}
& \sup_{t \in (0,T_1]} \Big ( \| \bu \|^2 + \| \bo \|^2 + \| \phi \|_{H^2}^2 \Big ) + \int_0^{T_1} \| \Delta^2 \phi \|^2 + \| \bu \|_{H^1}^2 + \| \bo \|_{H^1}^2 \, dt \\
& \quad \leq C ( \| \bu(0) \|^2 + \| \bo(0) \|^2 + \| \phi(0) \|_{H^2}^2 ),
\end{align*}
which also ensures the uniqueness of strong solutions in the time interval $[0,T_1]$.

\section{Nonpolar limit}\label{sec:etar0}

We recall the main result of \cite{CHNS3Dstrsol} concerning the local well-posedness of strong solution to the model of Abels, Garcke and Gr\"un \cite{AGG}:
\begin{prop}[Theorem 1.1 of \cite{CHNS3Dstrsol}]\label{thm:AGG}
Let $\Omega$ be a bounded domain of class $C^3$ in $\R^3$. Assume that $\bu_0 \in \HH^1_{\div}$ and $\phi_0 \in H^2(\Omega)$ such that $\| \phi_0\|_{L^\infty} \leq 1 - \delta_0$ for some $\delta_0 > 0$, $|\overline{\phi_0}| < 1$, $\mu_0 = - \Delta \phi_0 + F'(\phi_0) \in H^1(\Omega)$, and $\pdnu \phi_0 = 0$ on $\pd \Omega$. Then, there exists $T_a > 0$, depending on the norms of the initial data and $\delta_0$, and a quadruple of functions $(\bu_a, p_a, \phi_a, \mu_a)$ with the regularities
\begin{align*}
\bu_a & \in C^0([0,T_a);\HH^1_{\div}) \cap L^2(0,T_a; \HH^2_{\div}) \cap H^1(0,T_a; \HH_{\div}), \\
p_a & \in L^2(0,T_a;H^1(\Omega)), \\
\phi_a & \in L^\infty(0,T_a;W^{2,6}(\Omega)) \cap H^1(0,T_a;H^1(\Omega)) \cap W^{1,\infty}(0,T_a;H^1(\Omega)^*), \\
\phi_a & \in L^\infty(\Omega \times (0,T_a)) \text{ s.t. } |\phi_a(x,t)| < 1 \text{ a.e.~in } \Omega \times (0,T_a), \\
\Delta \phi_a & \in L^2(0,T_a; H^2(\Omega)), \\
\mu_a & \in L^\infty(0,T_a;H^1(\Omega)) \cap L^2(0,T_a;H^3(\Omega)),
\end{align*}
that is the unique strong solution to the Abels--Garcke--Gr\"un model in the time interval $[0,T_a)$, i.e.,
\begin{subequations}\label{agg:model}
\begin{alignat}{2}
\label{agg:phi} \nd{\phi_a} & = \Delta \mu_a, \\[1ex]
\label{agg:mu} \mu_a & =  F'(\phi_a) - \Delta \phi_a,\\[1ex]
\label{agg:div} \div \bu_a & = 0, \\[1ex]
\label{agg:lin:mom} \rho(\phi_a) \nd{\bu_a} & = - \nabla p_a - \div(  \nabla \phi_a \otimes \nabla \phi_a) + \div (2 \eta(\phi_a) \D \bu_a)  + \rho'(\phi_a) \big ( \nabla \mu_a \cdot \nabla \big ) \bu_a,
\end{alignat} 
\end{subequations}
holding a.e.~in $\Omega \times (0,T_a)$ with 
\begin{align}\label{agg:bc}
\begin{cases}
 \bu_a = \bm{0}, \quad \pd_{\bnu} \phi_a = 0, \quad \pd_{\bnu} \mu_a = 0 & \text{ on } \partial \Omega \times (0,T_a), \\
\bu_a(0) = \bu_0, \quad \phi_a(0)= \phi_0 & \text{ in } \Omega.
\end{cases}
\end{align}
\end{prop}

Similarly, we recall the main result of \cite{AGAMRT} concerning the local well-posedness of strong solution to the model H of Hohenberg and Halperin \cite{PCHoBIH}:
\begin{prop}[Theorem 5.1 of \cite{AGAMRT}]\label{thm:modelH}
Let $\Omega$ be a bounded domain of class $C^3$ in $\R^3$. Assume that $\bu_0 \in \HH^1_{\div}$ and $\phi_0 \in H^2(\Omega)$ such that $\| \phi_0\|_{L^\infty} \leq 1$, $|\overline{\phi_0}| < 1$, $\mu_0 = - \Delta \phi_0 + F'(\phi_0) \in H^1(\Omega)$, and $\pdnu \phi_0 = 0$ on $\pd \Omega$. Then, there exists $T_h > 0$, depending on the norms of the initial data, and a quadruple of functions $(\bu_h, p_h, \phi_h, \mu_h)$ with the regularities
\begin{align*}
\bu_h & \in C^0([0,T_h);\HH^1_{\div}) \cap L^2(0,T_h; \HH^2_{\div}) \cap H^1(0,T_h; \HH_{\div}), \\
p_h & \in L^2(0,T_h;H^1(\Omega)), \\
\phi_h & \in L^\infty(0,T_h;W^{2,6}(\Omega)) \cap H^1(0,T_h;H^1(\Omega)) \cap W^{1,\infty}(0,T_h;H^1(\Omega)^*), \\
\phi_h & \in L^\infty(\Omega \times (0,T_h)) \text{ s.t. } |\phi_h(x,t)| < 1 \text{ a.e.~in } \Omega \times (0,T_h), \\
\mu_h & \in L^\infty(0,T_h;H^1(\Omega)) \cap L^2(0,T_h;H^3(\Omega)),
\end{align*}
that is the unique strong solution to Model H in the time interval $[0,T_h)$, i.e.,
\begin{subequations}\label{modelh:model}
\begin{alignat}{2}
\label{modelh:phi} \nd{\phi_h} & = \Delta \mu_h, \\[1ex]
\label{modelh:mu} \mu_h & =  F'(\phi_h) - \Delta \phi_h,\\[1ex]
\label{modelh:div} \div \bu_h & = 0, \\[1ex]
\label{modelh:lin:mom} \overline{\rho}\nd{\bu_h} & = - \nabla p_h - \div(  \nabla \phi_h \otimes \nabla \phi_h) + \div (2 \eta(\phi_h) \D \bu_h) ,
\end{alignat} 
\end{subequations}
holding a.e.~in $\Omega \times (0,T_h)$ with constant mass density $\overline{\rho}$ and
\begin{align}\label{modelh:bc}
\begin{cases}
\bu_h = \bm{0}, \quad \pd_{\bnu} \phi_h = 0, \quad \pd_{\bnu} \mu_h = 0 & \text{ on } \partial \Omega \times (0,T_h), \\
 \bu_h(0) = \bu_0, \quad \phi_h(0) = \phi_0 & \text{ in } \Omega.
\end{cases}
\end{align}
\end{prop}

\begin{remark}
Although it is not explicitly pointed out in \cite{AGAMRT}, but via similar arguments we can deduce for initial data $\phi_0$ satisfying $\| \phi_0 \|_{L^{\infty}} \leq 1 - \delta_0$ for some $\delta_0 > 0$, it holds that $\Delta \phi_h \in L^2(0,\hat{T}_h;H^2(\Omega))$ for some $\hat{T}_h \in (0,T_h)$.
\end{remark}

\subsection{Uniformality of the well-posedness time interval}\label{sec:uniformtime}
In this section we set $\eta_r(\cdot)$ as a constant $\eta_r > 0$ and investigate the difference between the strong solutions to the MAGG model \eqref{bu:PF:model}, denoted as $(\bu_w, \bo_w, p_w, \phi_w, \mu_w)$, and the strong solution to the AGG model \eqref{agg:model}, denoted as $(\bu_a, p_a, \phi_a, \mu_a)$, originating from the same initial conditions $\bu_w(0) = \bu_a(0) = \bu_0$, $\phi_w(0) = \phi_a(0) = \phi_0$, while we consider $\bo_0 = \bm{0}$ for the micro-rotation. 

We note that the initial energy $E(\bu_w(0), \bo_w(0), \phi_w(0))$ is independent of $\eta_r$, and thus the analogous (standard) estimates to \eqref{uni:est:1}-\eqref{uni:est:4} are uniform in $\eta_r$. For the analogous higher order estimates, we first note that due to $\eta_r(\cdot)$ being a constant, the terms $B_3$, $B_{9}$, $B_{10}$, $B_{19}$ and $B_{22}$ vanish, while $B_{20}$ simplifies to $(2 \eta_r \curl \bo_m, \bm{A} \bu_m)$. Furthermore, let us point out the following identity (where the Einstein notation for summation is used, and $\bm{n}$ denotes the outward unit normal on $\pd \Omega$):
\begin{equation}\label{curlu_w_grad}
\begin{aligned}
& \| \nabla (\tfrac{1}{2} \curl \bu - \bo) \|^2 = (\pd_i (\tfrac{1}{2} \eps_{jkl} \pd_k u_l  - w_j), \pd_i (\tfrac{1}{2} \eps_{jmn} \pd_m u_n - w_j)) \\
& \quad = (\pd_i w_j, \pd_i w_j) + \tfrac{1}{4} (\eps_{jmn} \eps_{jkl}  \pd_{i} \pd_k u_l, \pd_i \pd_m u_n) \\
& \qquad - \tfrac{1}{2} \eps_{jmn} (\pd_i \pd_m u_n, \pd_i w_j) - \tfrac{1}{2}  (\eps_{jkl} \pd_i \pd_k u_l, \pd_i w_j) \\
&  \quad= (-\bo, \Delta \bo) + (\pd_i w_j, w_j n_i)_{\pd \Omega} + \tfrac{1}{4} (\pd_{i} \pd_k u_l, \pd_i \pd_k u_l - \pd_i \pd_l u_k) \\
& \qquad - \tfrac{1}{2}  (\eps_{jmn} \pd_i \pd_m u_n, \pd_i w_j) - \tfrac{1}{2}  (\eps_{jkl} \pd_i \pd_k u_l, \pd_i w_j) \\
& \quad = (-\bo, \Delta \bo) +  (\pd_i w_j, w_j n_i)_{\pd \Omega} - \tfrac{1}{2} (\pd_i \pd_i \pd_k u_l, (\W \bu)_{kl}) + \tfrac{1}{2} (\pd_i \pd_k u_l, (\W \bu)_{kl} n_i)_{\pd \Omega}  \\
& \qquad + \tfrac{1}{2} \eps_{jmn} (\pd_i \pd_i \pd_m u_n, w_j) - \tfrac{1}{2} (\eps_{jmn} \pd_{i} \pd_m u_n, w_j n_i)_{\pd \Omega} \\
& \qquad +  \tfrac{1}{2} (\eps_{jkl} \pd_k u_l, \pd_i \pd_i w_j) - \tfrac{1}{2}(\eps_{jkl} \pd_i w_j, \pd_k u_l n_i)_{\pd \Omega} \\
& \quad = (-\bo, \Delta \bo) + (\pd_i w_j, w_j n_i)_{\pd \Omega} \\
& \qquad + \tfrac{1}{2} (\Delta \bu, \div (\W \bu)) - \tfrac{1}{2} (\pd_i \pd_i u_l, (\W \bu)_{kl} n_k)_{\pd \Omega} + \tfrac{1}{2}(\pd_i \pd_k u_l, (\W \bu)_{kl} n_i)_{\pd \Omega} \\
& \qquad - \tfrac{1}{2} (\eps_{jmn} \pd_i \pd_i u_n, \pd_m w_j) + \tfrac{1}{2}(\eps_{jmn} \pd_i \pd_i u_n, w_j n_m)_{\pd \Omega} - \tfrac{1}{2} (\eps_{jmn} \pd_{i} \pd_m u_n, w_j n_i)_{\pd \Omega} \\
& \qquad +  \tfrac{1}{2} (\eps_{jkl} \pd_k u_l, \pd_i \pd_i w_j) - \tfrac{1}{2}(\eps_{jkl} \pd_i w_j, \pd_k u_l n_i)_{\pd \Omega} \\
& \quad = (-\bo, \Delta \bo) + (\pd_i w_j, w_j n_i)_{\pd \Omega} \\
& \qquad + \tfrac{1}{2} (\Delta \bu, \div (\W \bu)) - \tfrac{1}{2} (\pd_i \pd_i u_l, (\W \bu)_{kl} n_k)_{\pd \Omega} + \tfrac{1}{2}(\pd_i \pd_k u_l, (\W \bu)_{kl} n_i)_{\pd \Omega} \\
& \qquad + \tfrac{1}{2} (\curl \bo, \Delta \bu) + \tfrac{1}{2}(\eps_{jmn} \pd_i \pd_i u_n, w_j n_m)_{\pd \Omega} - \tfrac{1}{2} (\eps_{jmn} \pd_{i} \pd_m u_n, w_j n_i)_{\pd \Omega} \\
& \qquad +  \tfrac{1}{2} (\curl \bu, \Delta \bo) - \tfrac{1}{2}(\eps_{jkl} \pd_i w_j, \pd_k u_l n_i)_{\pd \Omega}.
\end{aligned}
\end{equation}
The no-slip and no-spin boundary conditions enable only the first, fourth and fifth boundary terms on the right-hand side to vanish, and thus the sum of the analogous of $B_{20}$, $B_{29}$ and $B_{30}$ along with $\frac{1}{2} \eta_r (\bm{A} \bu_m, \bm{A} \bu_m)$ on the left-hand side of \eqref{uni:est:6} cannot be fully expressed as $\eta_r \| \nabla (\tfrac{1}{2} \curl \bu_m - \bo_m) \|^2$ under the no-slip and no-spin boundary conditions.

Hence, with the way $B_{20}$, $B_{29}$ and $B_{30}$ are currently estimated, the positive constant $C$ in \eqref{reg:main3} is dependent on $\eta_r$, leading to the dependence of the existence time on the value of $\eta_r$. To study the nonpolar limit $\eta_r \to 0$, we restrict ourselves to the range $\eta_r \in (0,1]$ when considering the standard boundary conditions \eqref{PF:BC}, and perform estimates for $B_{20}$, $B_{29}$ and $B_{30}$ in a different manner:
\begin{align*}
|B_{20}| & \leq 2 \eta_r \| \bm{A} \bu_m \| \| \curl \bo_m \| \leq \frac{\eta_r^2}{2} \| \bm{A} \bu_m \|^2 + C (1 + C_0 + H_m)^2 \\
& \leq \frac{\eta_r}{2} \| \bm{A} \bu_m \|^2 + C (1 + C_0 + H_m)^2 , \\
|B_{29} + B_{30}| & \leq 4 \eta_r \| \Delta \bo_m \| \| \tfrac{1}{2} \curl \bu_m - \bo_m \| \leq \frac{c_*}{16} \| \Delta \bo_m \|^2 + C \eta_r^2 \| \tfrac{1}{2} \curl \bu_m - \bo_m \|^2 \\
& \leq \frac{c_*}{16} \| \Delta \bo_m \|^2 + C(1 + C_0 + H_m)^2,
\end{align*}
with positive constants independent of $\eta_r$.  Then, choosing $\beta_1, \beta_2$ and $\beta_3$ sufficiently small, we arrive at the differential inequality \eqref{reg:main3} with positive constants independent of $\eta_r \in (0,1]$. Furthermore, due to $\eta_r \in (0,1]$, we see that $H_m(0)$ can be bounded above by a constant $L_0$ that is independent of $\eta_r$. Hence, repeating the argument yields a time interval $[0,T_1]$, with $T_1$ independent of $\eta_r \in (0,1]$ for which there exists a unique strong solution to the MAGG model \eqref{bu:PF:model}.

\subsection{Under periodic boundary conditions}
Let us now consider $\Omega = \TT^3$ as the 3-torus and replace \eqref{PF:BC} with periodic boundary conditions.  As mentioned in \cite[Remark 4.1]{CHNS2Dstrsol}, the proof of Theorem \ref{thm:e&u} holds true in the case of periodic boundary conditions, where $\{\bm{Y}_j\}_{j=1}^\infty$ can be chosen as the eigenfunctions of the Stokes operator augmented by the constant function, and instead of choosing $\bm{v} = \bm{A} \bu_m$ in \eqref{app:u} to obtain the identity \eqref{uni:est:6}, it is sufficient to take $\bm{v}= - \Delta \bu_m$, which in turn leads to the absence of the terms $B_{21}$ and $B_{22}$ involving $p_m$.

Then, returning to the identity \eqref{curlu_w_grad}, we infer the validity of the following relation
\begin{align}\label{nabla:curl:u:w:periodic}
4 \| \nabla (\tfrac{1}{2} \curl \bu - \bo) \|^2 = (2 \curl \bu - 4\bo, \Delta \bo)  + (2 \curl \bo + 2 \div (\W \bu), \Delta \bu),
\end{align}
and thus choosing $\bm{v} = - \Delta \bu_m$ in \eqref{app:u} yields the alternate identity to \eqref{uni:est:6}:
\begin{equation}\label{uni:est:6:periodic}
\begin{aligned}
& (\eta(\phi_m) \Delta \bu_m, \Delta \bu_m) + (2 \eta_r \div (\W \bu_m), \Delta \bu_m) + (2 \eta_r \curl \bo_m, \Delta \bu_m) \\
& \quad = (\rho(\phi_m) \pd_t \bu_m, \Delta \bu_m) + (\rho(\phi_m) (\bu_m \cdot \nabla ) \bu_m, \Delta \bu_m) \\
& \qquad - \tfrac{\overline{\rho}_1 - \overline{\rho}_2}{2} (( \nabla \mu_m \cdot \nabla) \bu_m, \Delta \bu_m) - (\mu_m \nabla \phi_m, \Delta \bu_m) \\
& \qquad - (2 \eta'(\phi_m) \nabla \phi_m \D \bu_m, \Delta \bu_m) \\
& \quad =: \widetilde{B}_{14} + \cdots + \widetilde{B}_{18}.
\end{aligned}
\end{equation}
Adding the above to \eqref{uni:est:10} and employing \eqref{nabla:curl:u:w:periodic}, we arrive at
\begin{equation}
\begin{aligned}
& \eta_* \| \Delta \bu_m \|^2 + 4 \eta_r \| \nabla (\tfrac{1}{2} \curl \bu_m - \bo_m) \|^2 \\
& \qquad + \int_\Omega (c_d(\phi_m) + c_a(\phi_m)) |\Delta \bu_m|^2 + c_0(\phi_m) |\nabla \div \bo_m|^2 \, dx \\
& \quad \leq \widetilde{B}_{14} + \cdots \widetilde{B}_{18} + B_{23} + \cdots + B_{28}.
\end{aligned}
\end{equation}
Hence, when choosing $\beta_1 = \beta_2$, we arrive at the following analogue of \eqref{reg:main2}:
\begin{equation}
\begin{aligned}
& \frac{d}{dt} H_m(t) + \rho_*(\|\pd_t \bu_m \|^2 + \|\pd_t \bo_m\|^2 ) + \| \nabla \pd_t \phi_m \|^2 + \beta_1 \eta_* \| \Delta \bu_m \|^2  \\
& \qquad + 4 \beta_1 \eta_{r} \|\nabla (\tfrac{1}{2} \curl \bu_m - \bo_m) \|^2 + \beta_3 \| \mu_m \|_{H^3}^2  \\
& \qquad + \beta_1 \int_\Omega (c_d(\phi_m) + c_a(\phi_m)) |\Delta \bo_m|^2 + c_0(\phi_m) |\nabla \div \bo_m|^2 \, dx \\
& \leq C(1 + C_0 + H_m) + (B_1 + \cdots + B_{13}) + \beta_1 (\widetilde{B}_{14} + \cdots + \widetilde{B}_{18}) \\
& \quad + \beta_1 (B_{23} + \cdots + B_{28}) + \beta_3(B_{31} + B_{32}).
\end{aligned}
\end{equation}
We note that the estimates for $\widetilde{B}_{14}, \dots, \widetilde{B}_{18}$ proceed analogously as for $B_{14}, \dots, B_{18}$, which yields with sufficiently small $\beta_1$ and $\beta_3$ the differential inequality
\begin{equation}\label{reg:main3:alt}
\begin{aligned}
& \frac{d}{dt} H_m(t) + \frac{\rho_*}{2} (\|\pd_t \bu_m \|^2 +  \|\pd_t \bo_m\|^2) + \frac{\beta_3}{4}\| \mu_m \|_{H^3}^2 + 4 \eta_r \| \nabla ( \tfrac{1}{2} \curl \bu_m - \bo_m) \|^2 \\
& \qquad + \frac{1}{2} \| \nabla \pd_t \phi_m \|^2 + \frac{\beta_1 \eta_*}{2} \| \bm{A} \bu_m \|^2 + \frac{\beta_1 c_*}{4} (\|\Delta \bo_m\|^2 + \|\nabla \div \bo_m\|^2)  \\
& \quad \leq C(1 + C_0 + H_m)^5,
\end{aligned}
\end{equation}
with a positive constant $C$ independent of $\eta_r > 0$ and the approximation parameters $m$, $\alpha$ and $k$. Then, for arbitrary but fixed constant $R > 0$ and for any $\eta_r \in (0,R]$, we can argue as in Section \ref{sec:uniformtime} to infer that $H_m(0)$ can be bounded above by a constant $L_0$ that is independent of $\eta_r$. Hence, there exists a time instance $T_1$ independent of $\eta_r \in (0,R]$ for which there exists a unique strong solution to the MAGG model \eqref{bu:PF:model} under periodic boundary conditions.

\subsection{Consistency estimates with AGG solutions}
Without loss of generality, we assume $T_1 \leq T_a$ and work on the common time interval $[0,T_1]$. Returning to the setting of the uniqueness proof, we set $(\bu_1, \bo_1, p_1, \phi_1, \mu_1)$ as $(\bu_w, \bo_w, p_w, \phi_w, \mu_w)$ and set $(\bu_2, \bo_2, p_2, \phi_2, \mu_2)$ as $(\bu_a, \bm{0}, p_a, \phi_a, \mu_a)$. Then, the differences $(\bu, \bo, p, \phi, \mu)$ satisfy
\begin{subequations}
\begin{alignat}{2}
\label{stab:u} & \rho_w \pd_t \bu + \rho \pd_t \bu_a + \rho_w( \bu_w \cdot \nabla) \bu + \rho_w (\bu \cdot \nabla)  \bu_a + \rho (\bu_a \cdot \nabla) \bu_a \\
\notag & \qquad - \tfrac{\overline{\rho}_1 - \overline{\rho}_2}{2} (\nabla \mu_w \cdot \nabla) \bu - \tfrac{\overline{\rho}_1 - \overline{\rho}_2}{2} (\nabla \mu \cdot \nabla) \bu_a  - \div (2 \eta_w \D \bu + 2 \eta \D \bu_a) \\
\notag & \qquad - \div (2 \eta_w \W \bu_w) + \nabla p \\
\notag & \quad = 2 \curl( \eta_{r} \bo_w) - \div (\nabla \phi_w \otimes \nabla \phi) - \div (\nabla \phi \otimes \nabla \phi_a), \\
\label{stab:w} & \rho_w \pd_t \bo_w + \rho_w( \bu_w \cdot \nabla) \bo_w  - \div (c_{0,m} (\div \bo_w) \mathbb{I}) + 2 c_{d,m} \D \bo_w + 2 c_{a,m} \W \bo_w) \\
\notag & \quad = \tfrac{\overline{\rho}_1 - \overline{\rho}_2}{2} (\nabla \mu_w \cdot \nabla) \bo_w  + 2 \eta_{r} (\curl \bu_w - 2 \bo_w), \\
\label{stab:phi} & \pd_t \phi + \bu_w \cdot \nabla \phi + \bu \cdot \nabla \phi_a = \Delta \mu, \\
\label{stab:mu} & \mu = - \Delta \phi + F'(\phi_w) - F'(\phi_a),
\end{alignat}
\end{subequations}
almost everywhere in $\Omega \times (0,T_1)$, along with the initial-boundary conditions
\begin{align*}
\begin{cases}
\bu = \bm{0}, \quad \bo = \bm{0}, \quad \pdnu \phi = 0, \quad \pdnu \mu = 0 & \quad \text{ on } \pd \Omega \times (0,T_1), \\
\bu(0) = \bm{0}, \quad \bo(0) = \bm{0}, \quad \phi(0) = 0 & \quad \text{ on } \Omega.
\end{cases}
\end{align*}
Note that \eqref{uniq:mean} holds simply from integrating \eqref{stab:phi} over $\Omega$. Recalling the identity \eqref{rho:cancel} and the following analogue
\[
-(\pd_t \rho_w, \tfrac{1}{2} |\bo_w|^2)+ (\rho_w \bo_w , \nabla \tfrac{1}{2} |\bo_w|^2) - \tfrac{\overline{\rho}_1 - \overline{\rho}_2}{2} (\nabla \mu_w , \nabla \tfrac{1}{2} |\bo_w|^2) = 0,
\]
when testing \eqref{stab:u} with $\bu$, \eqref{stab:w} with $\bo_w$, taking the gradient of \eqref{stab:phi} and testing with $\nabla \Delta \phi$ and $\nabla \phi$, respectively, and summing all resulting identities including \eqref{uniq:mean} yields
\begin{equation}\label{stab:1}
\begin{aligned}
& \frac{1}{2} \frac{d}{dt} \Big ( \int_\Omega \rho_w |\bu|^2  + \rho_w |\bo_w|^2 + |\Delta \phi|^2 + |\nabla \phi|^2 \, dx + |\overline{\phi}|^2 \Big ) \\
& \qquad + \int_\Omega |\Delta^2 \phi|^2 + |\nabla \Delta \phi |^2 + 2 \eta_w |\D \bu|^2 \, dx \\
& \qquad + \int_\Omega c_{0,w} |\div \bo_w|^2 + 2 c_{d,w} |\D \bo_w|^2 + 2 c_{a,w} |\W \bo_w|^2 + 4 \eta_w \big | \tfrac{1}{2} \curl \bu_w - \bo_w \big |^2 \, dx \\
& \quad = -(\rho \pd_t \bu_a, \bu) - (\rho_w(\bu \cdot \nabla) \bu_a, \bu) - (\rho(\bu_a \cdot \nabla) \bu_a, \bu) + \tfrac{\overline{\rho}_1 - \overline{\rho}_2}{2} ((\nabla \mu \cdot \nabla) \bu_a, \bu) \\
& \qquad - (2 \eta \D \bu_a, \D \bu) + (\nabla \phi_w \otimes \nabla \phi + \nabla \phi \otimes \nabla \phi_a, \nabla \bu) - (2 \eta_r \W \bu_w, \W \bu_a)  \\
& \qquad + (2 \eta_r \bo_w, \curl \bu_a) + (\bu_w \cdot \nabla \phi, \Delta^2 \phi) + (\bu \cdot \nabla \phi_a, \Delta^2 \phi) \\
& \qquad + (\Delta (F'(\phi_w) - F'(\phi_a)), \Delta^2 \phi) - (\nabla (\bu_w \cdot \nabla \phi), \nabla \phi) - (\nabla (\bu \cdot \nabla \phi_a), \nabla \phi) \\
& \qquad + (\nabla (F'(\phi_w) - F'(\phi_a)), \nabla \Delta \phi) \\
& \quad =: E_1 + \cdots + E_{14}.
\end{aligned}
\end{equation}
We note that $E_1, \dots, E_6$ are analogous to $D_1, \dots, D_6$, while $E_{9}, \dots, E_{14}$ are analogous to $D_{17}, \dots, D_{22}$, and hence
\begin{align*}
|E_{1} + \cdots + E_6| & \leq \eta_* \| \D \bu \|^2  + \frac{1}{4} \| \Delta^2 \phi \|^2 + C(1 + \| \bu_a \|_{H^2}^2 + \| \pd_t \bu_a \|^2)( \| \bu \|^2 + \| \Delta \phi \|^2 + \| \phi \|_{H^1}^2), \\
|E_{9} + \cdots + E_{11}| & \leq \frac{1}{4} \| \Delta^2 \phi \|^2 + C( \| \bu \|^2 + \| \Delta \phi \|^2 + \| \phi \|_{H^1}^2), \\
|E_{12} + \cdots + E_{14}| & \leq \frac{1}{2} \| \nabla \Delta \phi \|^2 + C \| \phi \|_{H^1}^2 + \frac{\eta_*}{2} \| \D \bu \|^2.
\end{align*}
On the other hand, we have
\[
|E_{7} + E_8| = \eta_r |(\curl \bu_a, 2 \bo_w -\curl \bu_w)| \leq 2 \eta_r \| \tfrac{1}{2} \curl \bu_w - \bo_w \|^2 + C \eta_r \| \bu_a \|_{H^1}^2,
\]
and thus we infer from \eqref{stab:1} the following differential inequality
\begin{align*}
& \frac{1}{2} \frac{d}{dt} \Big (\int_\Omega \rho_w |\bu|^2  + \rho_w |\bo_w|^2 + |\Delta \phi|^2 + |\nabla \phi|^2 \, dx + |\overline{\phi}|^2 \Big ) + \int_\Omega \frac{1}{2} |\Delta^2 \phi|^2 + \frac{1}{2} |\nabla \Delta \phi |^2 +  \eta_* |\D \bu|^2 \, dx \\
& \qquad + \int_\Omega c_{0,m} |\div \bo_w|^2 + 2 c_{d,m} |\D \bo_w|^2 + 2 c_{a,m} |\W \bo_w|^2 + 2 \eta_r \big | \tfrac{1}{2} \curl \bu_w - \bo_w \big |^2 \, dx \\
& \quad \leq C ( 1 + \| \bu_a \|_{H^2}^2 + \| \pd_t \bu_a \|^2) ( \| \bu \|^2 + \| \Delta \phi \|^2) + C \eta_r \| \bu_a \|_{H^1}^2.
\end{align*}
Using that $\bu_a \in L^\infty(0,T_1;\HH^1_{\div})$ from Proposition \ref{thm:AGG}, by applying Gr\"onwall's inequality and then the lower bound for $\rho_w = \rho(\phi_w)$, we deduce that
\begin{equation}\label{stab:AGG:micro}
\begin{aligned}
& \sup_{t \in (0,T_1)} \Big (\| (\bu_w - \bu_a)(t) \|^2 + \| \bo_w(t) \|^2 + \|  (\phi_w - \phi_a)(t) \|_{H^2}^2 \Big ) \\
& \qquad + \int_0^{T_1} \|\Delta (\phi_w - \phi_a) \|_{H^2}^2 + \| \nabla (\bu_w - \bu_a) \|^2 + \| \bo_w \|_{H^1}^2 \, dt  \leq C \eta_r,
\end{aligned}
\end{equation}
with positive constant $C$ independent of $\eta_r \in (0,R]$.

\subsection{Consistency estimates with Model H solutions}
Without loss of generality, we assume $T_1 \leq T_a \leq \hat{T}_h$ and work on the common time interval $[0,T_1]$. We derive a consistency estimate between strong solutions $(\bu_a, p_a, \phi_a, \mu_a)$ of the AGG model and $(\bu_h, p_h, \phi_h, \mu_h)$ of Model H that is analogous to the estimate in \cite[Theorem 1.8]{CHNS3Dstrsol} but under stronger norms. Returning to the setting of the uniqueness proof, with $(\bu_a, \bm{0}, p_a, \phi_a, \mu_a)$ playing the role of $(\bu_1, \bo_1, p_1, \phi_1, \mu_1)$, and likewise $(\bu_h, \bm{0}, p_h, \phi_h, \mu_h)$ playing the role of $(\bu_2, \bo_2, p_2, \phi_2, \mu_2)$, we find that the differences $(\bu, p, \phi, \mu)$ satisfy
\begin{subequations}
\begin{alignat}{2}
\label{stab:h:u} & \rho_a \pd_t \bu + \rho \pd_t \bu_h + \rho_a( \bu_a \cdot \nabla) \bu + \rho_a (\bu \cdot \nabla)  \bu_h + \rho (\bu_h \cdot \nabla) \bu_h \\
\notag & \qquad - \tfrac{\overline{\rho}_1 - \overline{\rho}_2}{2} (\nabla \mu_a \cdot \nabla) \bu_a  - \div (2 \eta_a \D \bu + 2 \eta \D \bu_h) + \nabla p \\
\notag & \quad = - \div (\nabla \phi_a \otimes \nabla \phi) - \div (\nabla \phi \otimes \nabla \phi_h), \\
\label{stab:h:phi} & \pd_t \phi + \bu_a \cdot \nabla \phi + \bu \cdot \nabla \phi_h = \Delta \mu, \\
\label{stab:h:mu} & \mu = - \Delta \phi + F'(\phi_a) - F'(\phi_h),
\end{alignat}
\end{subequations}
almost everywhere in $\Omega \times (0,T_1)$, along with the initial-boundary conditions
\begin{align*}
\begin{cases}
\bu = \bm{0}, \quad \pdnu \phi = 0, \quad \pdnu \mu = 0 & \quad \text{ on } \pd \Omega \times (0,T_1), \\
\bu(0) = \bm{0}, \quad \phi(0) = 0 & \quad \text{ on } \Omega.
\end{cases}
\end{align*}
Testing \eqref{stab:h:u} with $\bu$, taking the gradient of \eqref{stab:h:phi} and testing with $\nabla \Delta \phi$ and $\nabla \phi$, respectively, and summing all resulting identities including \eqref{mean:est} yields
\begin{equation}\label{stab:2}
\begin{aligned}
& \frac{1}{2} \frac{d}{dt} \Big ( \int_\Omega \rho_a |\bu|^2  + |\Delta \phi|^2 + |\nabla \phi|^2 \, dx + |\overline{\phi}|^2 \Big ) + \int_\Omega |\Delta^2 \phi|^2 + |\nabla \Delta \phi |^2 + 2 \eta_a |\D \bu|^2 \, dx \\
& \quad = -(\rho \pd_t \bu_h, \bu) - (\rho_a(\bu \cdot \nabla) \bu_h, \bu) - (\rho(\bu_h \cdot \nabla) \bu_h, \bu) + \tfrac{\overline{\rho}_1 - \overline{\rho}_2}{2} ((\nabla \mu_a \cdot \nabla) \bu_a, \bu) \\
& \qquad - (2 \eta \D \bu_h, \D \bu) + (\nabla \phi_a \otimes \nabla \phi + \nabla \phi \otimes \nabla \phi_h, \nabla \bu)  + (\bu_a \cdot \nabla \phi, \Delta^2 \phi)\\
& \qquad + (\bu \cdot \nabla \phi_h, \Delta^2 \phi)  + (\Delta (F'(\phi_a) - F'(\phi_h)), \Delta^2 \phi) - (\nabla (\bu_a \cdot \nabla \phi), \nabla \phi)  \\
& \qquad - (\nabla (\bu \cdot \nabla \phi_h), \nabla \phi)+ (\nabla (F'(\phi_a) - F'(\phi_h)), \nabla \Delta \phi) \\
& \quad =: G_1 + \cdots + G_{12}.
\end{aligned}
\end{equation}
Once again we note that $G_{1}, \dots, G_{6}$ are analogous to $D_{1}, \dots, D_{6}$, while $G_{7}, \dots, G_{12}$ are analogous to $D_{17}, \dots, D_{22}$. However, $G_1$, $G_3$ and $G_4$ have to be treated differently.  Hence, we estimate as follows:
\begin{align*}
|G_{2}| & \leq C \| \bu \|_{L^3} \| \nabla \bu_a \|_{L^6} \| \bu \| \leq \frac{\eta_*}{4} \| \D \bu \|^2 + C \| \bu_a \|_{H^2}^2 \| \bu \|^2, \\
|G_{5} + G_{6}| & \leq \frac{\eta_*}{4} \| \D \bu \|^2 + C( 1 + \| \bu_h \|_{H^2}^2) \| \phi \|_{H^1}^2, \\
|G_{7} + \cdots + G_{12}| & \leq \frac{1}{2} \| \Delta^2 \phi \|^2 + \frac{\eta_*}{2} \| \D \bu \|^2 + \frac{1}{2} \| \nabla \Delta \phi \|^2 + C( \| \bu \|^2 + \| \Delta \phi \|^2 + \| \phi \|_{H^1}^2).
\end{align*}
For $G_{4}$ we follow similar ideas as for $D_{4}$: using the regularity $\phi_a \in L^\infty(0,T_1;W^{2,6}(\Omega))$,
\begin{align*}
|G_{4}| & \leq \frac{|\overline{\rho}_1 - \overline{\rho}_2|}{2} \big ( |((\nabla \Delta \phi_a \cdot \nabla) \bu_h, \bu)| + |((\nabla F'(\phi_a) \cdot \nabla) \bu_h, \bu)| \big ) \\
& \leq C|\overline{\rho}_1 - \overline{\rho}_2| \big ( \| \nabla \Delta \phi_a \|_{L^6} \| \nabla \bu_h \|_{L^3} \| \bu \| + \| F''(\phi_a) \|_{L^\infty} \| \nabla \phi_a \|_{L^6} \| \nabla \bu_h \|_{L^3} \| \bu \| \big ) \\
& \leq C ( 1 + \|\Delta \phi_a \|_{H^2}^2) \| \bu \|^2 + C |\overline{\rho}_1 - \overline{\rho}_2|^2 \| \nabla \bu_h \|_{L^3}^2.
\end{align*}
On the other hand, using the explicit form \eqref{rho:defn} for $\rho(\phi_a)$, we estimate $G_{1}$ and $G_{3}$ as follows:
\begin{align*}
|G_{1}| & \leq \Big (\Big |\frac{ \overline{\rho}_1 + \overline{\rho}_2}{2} - \overline{\rho}\Big | + \frac{|\overline{\rho}_1 - \overline{\rho}_2|}{2} \| \phi_a \|_{L^\infty} \Big ) \| \pd_t \bu_h \| \| \bu \| \\
& \leq C \| \pd_t \bu_h \|^2 \| \bu \|^2 + C\Big (\Big |\frac{ \overline{\rho}_1 + \overline{\rho}_2}{2} - \overline{\rho}\Big | + |\overline{\rho}_1 - \overline{\rho}_2| \Big ),\\
|G_{3}| & \leq \Big (\Big |\frac{ \overline{\rho}_1 + \overline{\rho}_2}{2} - \overline{\rho}\Big | + \frac{|\overline{\rho}_1 - \overline{\rho}_2|}{2} \| \phi_a \|_{L^\infty} \Big ) \| \bu_h \|_{L^6} \| \nabla \bu_h \|_{L^3} \| \bu \| \\
& \leq C \| \nabla \bu_h \|_{L^3}^2 \| \bu \|^2 + C\Big (\Big |\frac{ \overline{\rho}_1 + \overline{\rho}_2}{2} - \overline{\rho}\Big | + |\overline{\rho}_1 - \overline{\rho}_2| \Big ).
\end{align*}
Hence, we derive from \eqref{stab:2} the following differential inequality
\begin{align*}
& \frac{1}{2} \frac{d}{dt} \Big ( \int_\Omega \rho_a |\bu|^2  + |\Delta \phi|^2 + |\nabla \phi|^2 \, dx + |\overline{\phi}|^2 \Big ) + \int_\Omega \frac{1}{2} |\Delta^2 \phi|^2 + \frac{1}{2} |\nabla \Delta \phi |^2 + \eta_* |\D \bu|^2 \, dx \\
& \quad \leq C(1 + \| \bu_a \|_{H^2}^2 + \| \bu_h \|_{H^2}^2 + \| \Delta \phi_a \|_{H^2}^2)( \| \bu \|^2 + \| \Delta \phi \|^2 + \| \phi \|_{H^1}^2) \\
& \qquad + C(1 + \| \bu_h \|_{H^2}^2) \Big (\Big |\frac{ \overline{\rho}_1 + \overline{\rho}_2}{2} - \overline{\rho}\Big | + |\overline{\rho}_1 - \overline{\rho}_2| \Big ).
\end{align*}
Integrating over $(0,T_1)$ and using Propositions \ref{thm:AGG} and \ref{thm:modelH}, we deduce that 
\begin{align*}
& \sup_{t \in (0,T_1)} \Big ( \| (\bu_h - \bu_a)(t) \|^2 + \| (\phi_h - \phi_a)(t) \|_{H^2}^2 \Big )  + \int_0^{T_1} \| \Delta (\phi_h - \phi_a) \|_{H^2}^2 + \| \nabla (\bu_h - \bu_a) \|^2 \, dt \\
& \quad \leq C\Big (\Big |\frac{ \overline{\rho}_1 + \overline{\rho}_2}{2} - \overline{\rho}\Big | + |\overline{\rho}_1 - \overline{\rho}_2| \Big ).
\end{align*}
Together with \eqref{stab:AGG:micro} and the triangle inequality we infer the stability estimate between the strong solution $(\bu_w, \bw_w, p_w, \phi_w, \mu_w)$ to the MAGG model \eqref{bu:PF:model} and the strong solution $(\bu_h, p_h, \phi_h, \mu_h)$ to Model H:
\begin{align*}
& \sup_{t \in (0,T_1)} \Big ( \| (\bu_h - \bu_w)(t) \|^2 + \| \bo_w(t) \|^2 + \| (\phi_h - \phi_w)(t) \|_{H^2}^2 \Big ) \\
& \qquad + \int_0^{T_1} \| \Delta (\phi_h - \phi_w) \|_{H^2}^2 + \| \nabla (\bu_h - \bu_w) \|^2 + \| \bo_w \|_{H^1}^2 \, dt \\
& \quad \leq C\Big (\Big |\frac{ \overline{\rho}_1 + \overline{\rho}_2}{2} - \overline{\rho}\Big | + |\overline{\rho}_1 - \overline{\rho}_2| \Big ) + C \eta_r.
\end{align*}
with positive constant $C$ independent of $\eta_r \in (0,R]$, $|\overline{\rho}_1 - \overline{\rho}_2|$ and $|\frac{ \overline{\rho}_1 + \overline{\rho}_2}{2} - \overline{\rho} |$.

\section*{Acknowledgements}
\noindent KFL gratefully acknowledges the support by the Research Grants Council of the Hong Kong Special Administrative Region, China [Project No.: HKBU 14303420].

\end{document}